\documentclass[12pt]{article}
\usepackage{amsmath}
\usepackage{amssymb}
\usepackage{amsthm}
\usepackage{amscd}
\usepackage{amsfonts}
\usepackage{graphicx}
\usepackage{fancyhdr}
\usepackage{color}
\usepackage{enumerate}
\usepackage{amsfonts}
\usepackage{psfrag}
\usepackage{epsfig}
\usepackage{subfigure}
\usepackage{epstopdf}
\usepackage{amsbsy}
\usepackage{latexsym}
\usepackage{moreverb}
\usepackage{textcomp}

\newtheorem{remark}{Remark}

\newtheorem{theorem}{Theorem}
\newtheorem{conjecture}{Conjecture}

\def\Real{\mathbb{R}}  
\newcommand{\ray}{\mathbf q}

\usepackage[top=2.8cm,bottom=2.8cm,left=2.5cm,right=2.5cm]{geometry}
\usepackage[T1]{fontenc}
\usepackage[utf8]{inputenc}

\usepackage[affil-it]{authblk}
\usepackage{lmodern}


\usepackage{sectsty}
\sectionfont{\fontsize{13}{15}\selectfont}
\subsectionfont{\fontsize{12}{15}\selectfont}

\begin{document}

\title{A Sparse Stochastic Collocation Technique for High-Frequency Wave Propagation with Uncertainty
}

\author[1]{G. Malenova\thanks{malenova@kth.se}}
\author[2]{M. Motamed\thanks{motamed@math.unm.edu}}
\author[1]{O. Runborg\thanks{olofr@nada.kth.se}}
\author[3]{R. Tempone\thanks{raul.tempone@kaust.edu.sa}}
\affil[1]{KTH Royal Institute of Technology}
\affil[2]{The University of New Mexico}
\affil[3]{King Abdullah University of Science and Technology}

\date{\small \today}


\maketitle

\begin{abstract}

We consider the wave equation with highly oscillatory initial data, where there is uncertainty in the wave speed, initial phase and/or initial amplitude.
To estimate quantities of interest related to the solution and their statistics, we combine a high-frequency method based on Gaussian beams with sparse stochastic collocation.
Although the wave solution, $u^\varepsilon$, is highly oscillatory in both physical and
stochastic spaces, we provide theoretical arguments and numerical evidence
that quantities of interest based on local averages of $|u^\varepsilon|^2$
are smooth, with derivatives in the stochastic space uniformly bounded
in $\varepsilon$, where $\varepsilon$ denotes the short wavelength. This observable related regularity makes the sparse stochastic collocation approach
more efficient than Monte Carlo methods.
We present numerical tests that demonstrate this advantage.
\end{abstract}

\section{Introduction}

The propagation of high-frequency waves is modeled by hyperbolic partial differential equations (PDEs) with highly oscillatory solutions. In these models, the wavelength $\varepsilon$ is very short compared to variations in the structure of the medium and the wave propagation distance, resulting in multiscale problems. In addition, these models are often subject to uncertainty due to both incomplete knowledge of the model's parameters (such as wave speed) and the intrinsic variability of the physical system (such as the location of an earthquake's hypocenter or the time frequency of volcanic forces). The problem therefore has two components: {\it multiple scales} and {\it uncertainty}. High levels of confidence in the predictions require understanding of the uncertainties in the model. This understanding can be obtained by a process called uncertainty quantification (UQ).
Furthermore, the systematic coupling and interaction of multiple scales and uncertainty must be considered.

Despite recent advances in the uncertainty quantification of a wide range of mathematical models, see, e.g.,\mbox{} \cite{Ghanem,LeMaitre_Knio,Motamed_etal:13,Motamed_etal:12, Xiu:10}, viable UQ methodologies for high-frequency waves are not well developed. Current methods are based on Monte Carlo sampling techniques \cite{Tsuji_etal:11}, which feature slow convergence rates. {The main reason that more efficient methods with fast spectral convergence rates are not developed for high-frequency waves is that the wave solution is highly oscillatory and consequently its derivatives in the stochastic space are not bounded uniformly with respect to $\varepsilon$.
Spectral techniques in the stochastic space have to
resolve these oscillations and will thus require very many sample
points to be accurate.
Here, we aim at developing efficient computational techniques that feature fast convergence rates independent of the oscillations.}

In this work, we are concerned mainly with the forward propagation of uncertainty in high-frequency waves, where the uncertainty in the input parameters, such as wave speed and initial data, propagates through the multiscale hyperbolic model to give information about uncertain output quantities of interest (QoIs).  As the prototype model equation for high-frequency waves, we consider the scalar wave equation with highly oscillatory initial data. The main sources of uncertainty are the wave speed and/or the initial data, which are here described by a finite number of independent random variables with known probability distributions. We propose a novel stochastic spectral asymptotic method, which combines two techniques:
(1) a Gaussian beam summation method for propagating high frequency waves;
and (2) a sparse stochastic collocation method for propagating uncertainty and approximating the statistics of QoIs.
An important property is the {\em stochastic regularity} of the QoIs.
By this, we mean the regularity of output QoIs with respect to input random parameters.
The fast spectral convergence of the proposed method depends
crucially on the presence of
high stochastic regularity,
{\em independent of the wave frequency $\varepsilon$}.
In particular, through both theoretical arguments for simplified problems and numerical experiments for more complicated problems, we show that although the derivatives of highly oscillatory wave solution $u^\varepsilon$ with respect to the random parameters are not bounded independently of $\varepsilon$, physically motivated
QoIs based on local averages of $|u^\varepsilon|^2$
are smooth with uniformly bounded derivatives in the stochastic space. Consequently, despite the slow algebraic convergence of the wave solution, the proposed method gives fast spectral convergence for those QoIs.

We note that the two methods composing the proposed method have already been employed and are not new. The stochastic collocation method has been widely used in forward propagation of uncertainty in many PDE models, see, e.g., \cite{Xiu_Hesthaven,BNT:10,NTW:08,NTW:08_2,BNTT:12,Motamed_etal:12,Ernst_Sprung:14,Gunzburger_etal:14,BMT:14}. The Gaussian beam summation method has also been successfully applied to deterministic high-frequency wave propagation problems, see, e.g., \cite{Runborg07,Tanushev08,MR_Cambridge:09,MotRun_Error:10,JinMarSpa:12,LRT:13,MotRun:15}. However, combining the two methods is not straightforward since the
convergence rate of the resulting method strongly depends on the
stochastic regularity of the QoI, which in general depends on $\varepsilon$ and in principle may be as bad as the regularity of $u^\varepsilon$.
The main contributions of this work include: (1) constructing a novel algorithm for the uncertainty propagation of high-frequency waves by systematic coupling of sparse stochastic collocation and Gaussian beam approximation; and (2) verifying the fast spectral convergence of the proposed algorithm for particular types of quadratic quantities, independent of frequency.

The rest of the work is organized as follows. In Section 2, we formulate the high-frequency wave propagation problem with stochastic parameters. The stochastic spectral asymptotic method for the forward propagation of uncertainty in the problem is presented in Section 3. In Section 4, we discuss the stochastic regularity of the highly oscillatory solution and a particular quantity of interest, through theoretical arguments and numerical experiments.
In Section 5, we provide numerical examples and demonstrate the spectral convergence of the proposed method for the quantity of interest. Finally, we summarize our conclusions in Section 6.

\section{Problem Statement}
\label{sec2}


Consider the Cauchy problem for the stochastic wave equation
\begin{subequations}\label{wave}
\begin{align}
&u^\varepsilon_{tt}(t,\bold{x},{\bf y})- c(\bold{x},{\bf y})^2 \, \Delta u^\varepsilon(t,\bold{x},{\bf y}) = 0, \qquad &&\text{in  } [0,T] \times {\mathbb R}^n \times \Gamma, \label{wave_1st}\\
&u^\varepsilon(0,\bold{x},{\bf y})=g^\varepsilon_1(\bold{x},{\bf y}), \hskip 0.5cm u^\varepsilon_t(0,\bold{x},{\bf y})=g^\varepsilon_2(\bold{x},{\bf y}), \qquad &&
\text{on } \{ t=0 \} \times {\mathbb R}^n \times \Gamma, \label{wave_2nd}
\end{align}
\end{subequations}
where $u^\varepsilon \in {\mathbb C}$ is the stochastic wave function, $t\in [0,T]$ is the time, ${\bf x}=(x_1,\dotsc,x_n) \in {\mathbb R}^n$ is the spatial variable, ${\bf y}= (y_1,\dotsc,y_N) \in \Gamma \subset \mathbb{R}^N$ is a random vector, and $\varepsilon \ll 1$ is a small positive number, representing a typical short wavelength defined below.
We use the convention that $\nabla$ represents the gradient with respect to the spatial variables ${\bf x}$.

\medskip
\noindent
{\bf Sources of uncertainty.}
The uncertainty in model \eqref{wave} is described by a random vector, ${\bf y}$, consisting of $N \in {\mathbb N}_+$ independent random variables, $y_1, \dotsc,y_N$, with a bounded joint probability density, $\rho ({\bf y}) = \prod_{n=1}^N \rho_n(y_n): \Gamma \rightarrow \mathbb{R}_+$.
There may be two sources of uncertainty:
uncertainty in the local wave speed, $c = c(\bold{x},{\bf y})$, and
uncertainty in the initial data, $g^\varepsilon_1=g^\varepsilon_1(\bold{x},{\bf y})$ and $g^\varepsilon_2=g^\varepsilon_2(\bold{x},{\bf y})$.
We make the following smoothness, uniform coercivity and boundedness assumptions on the wave speed:
\begin{gather}
c(\,\cdot\,,{\bf y}) \in C^{\infty}({\mathbb R}^n), \qquad \forall \, {\bf y} \in \Gamma, \label{assum_c1} \\
0<c_{\rm min} \le c(\bold{x},{\bf y}) \le c_{\rm max} < \infty, \qquad \forall \, {\bf x} \in {\mathbb R}^n, \quad \forall \, {\bf y} \in \Gamma.\label{assum_c2}
\end{gather}
In the case when $\Gamma$ is not compact (for example when random variables are Gaussian), in addition to \eqref{assum_c1}, we assume the spatial derivatives of the wave speed are uniformly bounded in $\Gamma$. If $\Gamma$ is compact (for example when random variables are uniform), this extra assumption is automatically satisfied.
We further assume that the random wave speed has bounded mixed ${\bf y}$-derivatives of any order, i.e., for a multi-index ${\bf k} {= (k_1,\dotsc,k_N)} \in {\mathbb Z}_+^N$ with $|{\bf k}| \ge 0$, we {define $\partial_{\bf y}^{\bf k} := \frac{\partial^{|{\bf k}|}}{\partial_{y_1}^{k_1} \, \dotsc \, \partial_{y_N}^{k_N}}$ and} assume
\begin{equation}\label{assum3}
\left\| \partial_{\bf y}^{\bf k} c(\,\cdot\,,{\bf y}) \right\|_{L^{\infty}({\mathbb R}^n)} < \infty, \qquad \forall \,  {\bf y} \in \Gamma.
\end{equation}

\medskip
\noindent
{\bf High-frequency waves.}
We consider highly oscillatory initial data of the form
\begin{equation} \label{waveinit}
  g^\varepsilon_1(\bold{x},{\bf y})   = A_0(\bold{x},{\bf y}) \, e^{i \, \Phi_0(\bold{x},{\bf y})/\varepsilon} \ , \qquad
  g^\varepsilon_2(\bold{x},{\bf y}) = \frac{1}{\varepsilon}B_0(\bold{x},{\bf y}) \, e^{i \, \Phi_0(\bold{x},{\bf y})/\varepsilon} \ ,
\end{equation}
where the short wavelength, $\varepsilon$, is assumed to be much smaller than the typical scale of the structure of the medium (variations in the wave speed) and the wave propagation distance (the size of the computational domain). Such initial data generate high-frequency waves propagating in low-frequency media. The functions $A_0$, $B_0$ and $\Phi_0$ are assumed to be real-valued and smooth, with $|\nabla \Phi|$ bounded away from zero, $\forall \, {\bf y} \in \Gamma$. Consequently, we consistently assume
\begin{equation}\label{assum4}
g^\varepsilon_1(\,\cdot\,,{\bf y}), \, \, g^\varepsilon_2(\,\cdot\,,{\bf y}) \, \in C_{\text{c}}^{\infty}({\mathbb R}^n), \qquad \forall \,  {\bf y} \in \Gamma, \quad \forall \, \varepsilon > 0.
\end{equation}

\medskip
\noindent
{\bf Well-posedness of the problem.} By formally partial-differentiating problem \eqref{wave} with respect to $t$ and $\bold{x}$ and a simple extension of the proof for deterministic problems \cite{Lions_Magenes1,Hormander,Evans}, we can show that problem \eqref{wave} is well-posed. In other words, with the random wave speed satisfying \eqref{assum_c1}-\eqref{assum_c2} and the initial data satisfying \eqref{assum4}, there exists a unique solution $u^\varepsilon(\,\cdot\,, {\bf y}) \in C^{\infty}([0,T] \times {\mathbb R}^n)$ for every ${\bf y} \in \Gamma$ to problem \eqref{wave}, which depends continuously on the data. Moreover, at a fixed $(t, \bold{x})$, the solution belongs to the Hilbert space of stochastic functions with bounded second moments, i.e.
$$
u^\varepsilon(t, \bold{x}, .) \in L_{\rho}^2(\Gamma) := \left\{ v: \Gamma \rightarrow \mathbb{R}, \, \int_{\Gamma} |v({\bf y})|^2 \, \rho({\bf y}) \, d{\bf y} < \infty \right\},
$$
where the space $L_{\rho}^2$ is endowed with its standard inner product.
We also note that by partial-differentiating the problem \eqref{wave} with respect to ${\bf y}$ and employing assumption \eqref{assum3}, we can further show $u^\varepsilon \in C^{\infty}([0,T] \times {\mathbb R}^n \times \Gamma)$. We refer to \cite{Motamed_etal:12,Motamed_etal:13} for a more detailed discussion on the well-posedness of stochastic wave propagation problems.

The ultimate goal is the prediction of statistical moments of the high oscillatory solution, $u^\varepsilon$, or the statistics of some given physically motivated QoIs. In particular, we consider the following QoI
\begin{equation}\label{Q}
{\mathcal Q}^\varepsilon ({\bf y}) = \int_{{\mathbb R}^n} |u^\varepsilon(T, {\bf x}, {\bf y})|^2 \, \psi({\bf x}) \, d{\bf x}, 
\end{equation}
where $\psi$ is a given real-valued function, referred to as the {\it QoI test function}. Throughout this paper, the QoI test functions we consider will always be smooth and compactly supported, $\psi \in C^{\infty}_{\text{c}}({\mathbb R}^n)$, unless otherwise stated.

\begin{remark}
In the present work, we concentrate only on the quadratic QoI \eqref{Q}, which  represents the average strength of the wave inside the support of $\psi$. Other types of QoIs will be considered and studied in future works. We note in particular that the {linear} quantity $\widetilde{\mathcal Q}^\varepsilon ({\bf y}) = \int_{{\mathbb R}^n} u^\varepsilon(T, {\bf x}, {\bf y}) \, \psi({\bf x}) \, d{\bf x}$ typically
vanishes as $\varepsilon \rightarrow 0$ and is therefore not of interest.
\end{remark}

\section{A Stochastic Spectral Asymptotic Method}
\label{sec3}

In this section, we present an efficient numerical method for solving the stochastic wave equation \eqref{wave} with highly oscillatory initial functions \eqref{waveinit}. The method combines two techniques: a Gaussian beam summation method for propagating high-frequency waves; and a sparse stochastic collocation method for approximating the statistics of QoIs.

\subsection{Gaussian beam approximation}
\label{sec:3_1}

The Gaussian beam method describes high-frequency waves in a way closely related to geometrical optics and ray tracing. Geometrical optics formally assumes the solution of \eqref{wave} to be of the form
            \begin{equation}\label{GOansatz}
                u^\varepsilon(t,\mathbf x, \mathbf y) = a(t,\mathbf x,\mathbf y) e^{i \phi(t,\mathbf x, \mathbf y)/\varepsilon}.
            \end{equation}
As the phase $\phi$ and amplitude $a$ are independent of frequency and vary on a much coarser scale than does the solution, they might be resolved with a computational cost independent of $\varepsilon$.
In the limit $\varepsilon \to \infty$, we arrive at the \emph{eikonal} and \emph{transport equations}  \cite{Runborg07}:
            \begin{align}
                \phi_t + c(\mathbf x, \mathbf y) |\nabla \phi| & = 0, \label{eiko} \\
                a_t  + \frac{c^2(\mathbf x, \mathbf y) \Delta \phi -\phi_{tt}}{2 c(\mathbf x, \mathbf y) |\nabla \phi|} a+ \frac{c(\mathbf x,\mathbf y) \nabla a \cdot \nabla \phi}{|\nabla \phi|} & = 0 ,\label{transp}
            \end{align}
where $(t,\mathbf x,\mathbf y) \in (0,T] \times \mathbb R^n \times \Gamma $. Initial data is given as
            \[
                a(0,\mathbf x,\mathbf y) = A_0(\mathbf x, \mathbf y), \qquad \phi(0,\mathbf x, \mathbf y) = \Phi_0(\mathbf x, \mathbf y),
            \]
with $A_0$ and $\Psi_0$ as in \eqref{waveinit}.

The bicharacteristics $({\bf q}(t,{\bf y}), {\bf p}(t,{\bf y}))$ of the eikonal equation \eqref{eiko}
satisfy
                \begin{equation}\label{qp}
                  \frac{d \ray}{dt} = c(\ray, \mathbf y)\frac{\mathbf p}{|\mathbf p|}, \qquad \frac{d \mathbf{p}}{dt} = -\nabla c(\ray, \mathbf y) |\mathbf p|,
                  \qquad t>0.
                \end{equation}
Here the vector ${\bf q}$ is usually called the {\em ray} and the vector ${\bf p}$ is called the {\em slowness}. For
initial values $\ray(0,\mathbf y) = {\bold z}_0(\mathbf y)$, $\mathbf{p}(0,\mathbf y) = \nabla \Phi_0({\bold z}_0(\mathbf y), \mathbf y)$
we have that $\frac{d}{dt}\phi(t,\ray(t,\mathbf y),\mathbf y) = 0$ and
thus the phase $\phi$ is given by the simple expression $\phi(t,\ray(t,\mathbf y),\mathbf y) = \Phi_0(\mathbf z_0(\mathbf y),\mathbf y)$.
Moreover, the slowness vector always points in the direction of
the phase gradient, ${\bf p}(t,{\bf y}) = \nabla\phi(t,\ray(t,\mathbf y),\mathbf y)$.
In this way, the eikonal equation always admits a smooth solution
locally in time, more precisely
up to the point where two rays cross.
If we denote the ray starting at ${\bf z}$ by
$\mathbf{q}(t,\mathbf y;\mathbf z)$, then {\em caustic} points are those where this {\em ray function} degenerates in ${\bf z}$,
i.e.,\mbox{} where $\det \partial\mathbf{q}(t,\mathbf y;\mathbf z)/\partial \mathbf z=0$.
At these points, the initial ray position fails to be a well-defined
function of the current ray position.
 A global smooth solution of \eqref{eiko} is in general not feasible as caustics develop when rays cross. The form of the solution in \eqref{GOansatz} is thus not correct. The solution becomes multi-phased and geometrical optics predicts an unbounded amplitude \cite{Runborg07}, cf.\ the discussion in Section \ref{sec4}.

\emph{Gaussian beams} present another type of high frequency approximation closely related to geometrical optics evaluated along rays. However, unlike geometrical optics, a Gaussian beam is well-defined for every $t$ and their superposition performs well even at caustics.

A Gaussian beam, $v^\varepsilon(t,\bold x,\bold y)$, has the same form as \eqref{GOansatz}:
            \begin{equation}\label{GB}
                v^\varepsilon(t,\mathbf x, \bold y) = A(t,\mathbf x,\bold y) e^{i \Phi(t,\mathbf x, \bold y)/\varepsilon},
            \end{equation}
where the phase, $\Phi$, and amplitude, $A$, are centered around the geometrical optics ray, $\ray(t,\bold y)$:
            \begin{align*}
                A(t,\mathbf x,\bold y) & = a(t,\mathbf x-\ray(t,\bold y), \bold y),\\
                \Phi(t,\mathbf x, \bold y) & = \phi(t, \mathbf x - \ray(t,\bold y), \bold y).
            \end{align*}
The two methods differ in their assumption on the phase, $\phi$. Geometrical optics assumes $\phi$ is real whereas the Gaussian beam method sets the phase to be complex away from the ray. The imaginary part is chosen such that the Gaussian
beam \eqref{GB} decreases exponentially away from the central ray.

For the simplest (first order) Gaussian beams, $a$ and $\phi$ read
            \begin{align}
                a(t,\mathbf x, \bold y) &= a_0(t,\bold y),\label{a}\\
                \phi(t,\mathbf x, \bold y) &= \phi_0(t, \bold y) + \mathbf x\cdot \mathbf p(t, \bold y) + \frac{1}{2} \mathbf x \cdot M(t, \bold y) \mathbf x, \label{phi}
            \end{align}
where $M$ is a symmetric matrix with a positive imaginary part.
Note that these expressions can be interpreted as a zeroth and a second order
Taylor expansion in ${\bf x}$ of the amplitude and the phase respectively;
the slowness ${\mathbf p}$ is then the phase gradient, as in geometrical
optics, while the $M$ matrix is the Hessian of the phase.
By choosing $\ray, \mathbf p, \phi_0, M$ and $a_0$ in the right way, the phase will satisfy the eikonal equation \eqref{eiko} up to ${O(|\mathbf x - \ray(t,\bold y)|^3)}$ and the amplitude will satisfy the transport equation \eqref{transp} up to ${O(|\mathbf x - \ray(t,\bold y)|)}$. 
The coefficients $\phi_0,M,a_0$ then obey the following ordinary differential equations (ODE) (see \cite{Ralston82}):
\begin{subequations}\label{phi0ma0}
    \begin{align}
    \dot \phi_0 & = 0,\\
    \dot M & = D + B^T M + M B + M C M, \\
    \dot a_0 & = \frac{1}{2 |\mathbf p|} \left(c(\bold q,\bold y) \: \text{Tr}(M) - \nabla c(\bold q,\bold y) \cdot \mathbf p - \frac{c(\bold q,\bold y)\: \mathbf p \cdot M \mathbf p}{|\mathbf p|^2}\right) a_0,
    \end{align}
\end{subequations}
with
       \begin{equation*}
        D= |\mathbf p| \nabla^2 c(\bold q,\bold y), \qquad B = \frac{\mathbf p \otimes \nabla c(\bold q,\bold y)}{|\mathbf p|}, \qquad
        C = \frac{c(\bold q,\bold y)}{|\mathbf p|} \text{I} - \frac{c(\bold q,\bold y)}{|\mathbf p|^3} \mathbf p \otimes \mathbf p,
       \end{equation*}
We will always assume that the solution to these ODEs exists for the time
intervals considered. A sufficient condition is, for instance, that
$c(\bold x,\bold y)$ and all its $\bold x$-derivatives up to order three are
uniformly bounded for all ${\bold x}\in\Real^n$ and ${\bold y}\in\Gamma$.
\begin{remark}
Geometrical optics also allows for an alternative to the eikonal and transport equation pair \eqref{eiko} and \eqref{transp}, namely:
            \begin{equation*}
                \phi_t - c(\mathbf x, \mathbf y) |\nabla \phi| = 0, \qquad
                a_t  - \frac{c^2(\mathbf x, \mathbf y) \Delta \phi -\phi_{tt}}{2 c(\mathbf x, \mathbf y) |\nabla \phi|} a- \frac{c(\mathbf x,\mathbf y) \nabla a \cdot \nabla \phi}{|\nabla \phi|} = 0.
            \end{equation*}
This alternative corresponds to waves moving in the opposite direction and leads to opposite signs in the ODEs \eqref{qp} and \eqref{phi0ma0}.
\end{remark}
The name of the Gaussian beam is substantiated by its Gaussian shape. The imaginary part of $\phi$ is determined by $\mathbf x \cdot \operatorname{Im}(M) \: \mathbf x$; thus
    \[|v^\varepsilon(t,\mathbf x, \bold y)| = a_0 \exp{\left(-\frac{1}{2 \varepsilon}(\mathbf x- \ray(t, \bold y))\cdot \operatorname{Im}(M) (\mathbf x- \ray(t, \bold y))\right)}.\]
Since $M$ has a positive imaginary part, $|v(t, \mathbf x, \mathbf y)|$ is Gaussian with a width $|\mathbf x-\ray(t, \bold y)| \sim \sqrt{\varepsilon}$.
The positivity of $\text{Im}(M)$ is key to the Gaussian beam approximation.
Indeed,
\cite{Ralston82} states that $M = M^T$ and $\text{Im}(M) > 0$ hold true at any time, provided they are valid for the initial data.

In general, to construct higher order Gaussian beams, higher order terms in the Taylor expansion of $\phi$ and $a$ in \eqref{a}, and \eqref{phi} respectively, must be employed. Analogously to the above derivation, $\Phi(t,\mathbf x)$ is required to solve the eikonal equation \eqref{eiko} to $O(|\mathbf x-\ray(t,\mathbf y)|^{k+2})$ for the $k$-th order beams and the amplitude terms $a_j$ to solve the transport equations to $O(|\mathbf x-\ray(t,\mathbf y)|^{k-2j})$, respectively. This again translates to a system of ODEs for the coefficients in the Taylor expansions.

To approximate more general solutions, we use a superposition of Gaussian beams
            \begin{equation}\label{GBsup}
                u^\varepsilon_{GB}(t,\mathbf x, \bold y) = \frac{1}{(2\pi \varepsilon)^{n/2}}\int_{K_0} v^\varepsilon(t,\mathbf x, \bold y;\mathbf z)\: d\mathbf z,
            \end{equation}
over a compact domain, $K_0 \subset \mathbb R^n$, where the initial data is supported, and
$v^\varepsilon(t,\mathbf x, \bold y;\mathbf z)$ is a Gaussian beam that
starts in the point $\mathbf z \in K_0$.

Due to linearity of the wave equation
the superposition in \eqref{GBsup}
is also an asymptotic solution, since each individual beam is such a solution.
The pre-factor $(2\pi \varepsilon)^{-n/2}$ normalizes the superposition such that $u^\varepsilon_{\text{GB}} = O(1)$ away from caustics.
The coefficients of the beams are similarly parameterized by $\bold z$
and denoted $\ray(t,\bold y;\mathbf z)$, $\bold p(t,\bold y;\mathbf z)$, etc.
Following \cite{Tanushev08}, we choose the initial values to be
\begin{subequations}\label{xpMinit}
        \begin{align}
            \ray(0,\bold y;\mathbf z) & = \mathbf z,\\
            \mathbf p(0,\bold y;\mathbf z) & = \nabla\Phi_0(\mathbf z, \bold y), \\
            M(0,\bold y;\mathbf z) & = \nabla^2 \Phi_0(\mathbf z, \bold y) + i \text{I},\\
            \phi(0,\bold y;\mathbf z) & = \Phi_0(\mathbf z, \bold y),\\
            a_0(0,\bold y;\mathbf z) & = A_0(\mathbf z, \bold y).
        \end{align}
\end{subequations}
It has been shown in \cite{LRT:13} that the initial data for $k$-th order beams can be chosen such that the error in an $\varepsilon$-scaled energy norm satisfies the estimate
        \[
            \sup_{t\in [0,T]}\| u^\varepsilon_{GB}(t, \cdot) - u^\varepsilon(t, \cdot) \|_{E} \leq C(T) \varepsilon^{\frac{k}{2}},
        \]
for some constant $C(T)$ independent of $\varepsilon$.
Note that this holds regardless of the presence of caustics in the
solution.
Numerically, integral \eqref{GBsup} is approximated by the trapezoidal rule:
        \begin{equation}\label{trapets}
            u^\varepsilon_{GB}(t,\bold x, \bold y) \approx \frac{1}{(2\pi \varepsilon )^{n/2}} \sum_{\{j:\:z_j \in K_0\}} v^\varepsilon(t,\mathbf x, \bold y; \bold z_j) \; \Delta z^n,
        \end{equation}
where $\Delta z\sim \sqrt{\varepsilon}$,
and the ODEs \eqref{qp} and \eqref{phi0ma0} are solved with a numerical ODE method. The computational cost of the Gaussian beam method is then much smaller than that of a direct solver.

\subsection{Sparse Stochastic Collocation}
\label{sec:3_2}

For stochastic PDE models, such as the wave model \eqref{wave} with stochastic inputs, there are in general two types of methods for propagating uncertainty: intrusive and non-intrusive. Intrusive methods, such as perturbation expansion and stochastic Galerkin  \cite{Ghanem,Xiu_Karniadakis:02,Matthies,BTZ:05,Todor_Schwab}, require extensive modifications of existing deterministic solvers. On the contrary, non-intrusive methods, such as Monte Carlo \cite{MC} and stochastic collocation \cite{Xiu_Hesthaven,BNT:10,Motamed_etal:12}, are sample-based approaches. They rely on a set of deterministic models corresponding to a set of realizations and hence require no modifications of existing deterministic solvers. In this work, we consider the stochastic collocation method on sparse grids and briefly review the method for the uncertainty propagation of high-frequency waves.

The stochastic collocation method consists of three main steps. First, the problem \eqref{wave} is discretized in the physical space, i.e.,\ in space and time, using a deterministic asymptotic solver such as the Gaussian beam
method described in Section \ref{sec:3_1}. We therefore obtain a semi-discrete solution $u^\varepsilon_{\text{GB}}(t,\bold{x},{\bf y})$, keeping the variable ${\bf y}$ in the stochastic space continuous. The semi-discrete problem is then collocated on a set of $\eta \in {\mathbb N}_+$ collocation points, $\{ {\bf y}^{(k)} \}_{k=1}^{\eta} \in \Gamma$, to compute $\eta$ approximate solutions, $u^\varepsilon_{\text{GB}}(t,\bold{x},{\bf y}^{(k)})$. Finally, a global polynomial approximation, $u^\varepsilon_{\text{GB}, \eta}$, is built upon those evaluations,
 \begin{equation}\label{SC_sol}
 u^\varepsilon_{\text{GB}, \eta} (t,\bold{x},{\bf y})= \sum_{k=1}^{\eta} u^\varepsilon_{\text{GB}}(t,\bold{x},{\bf y}^{(k)}) \, L_k({\bf y}),
 \end{equation}
 for suitable multivariate polynomials, $\{ L_k ({\bf y})\}_{k=1}^{\eta}$, such as Lagrange polynomials, in the stochastic space. 

A key point in the stochastic collocation method is the choice of the set of collocation points, $\{ {\bf y}^{(k)} \}_{k=1}^{\eta}$, i.e.,\mbox{} the type of computational grid in the $N$-dimensional stochastic space. A full tensor grid, based on the Cartesian product of mono-dimensional grids, can only be used when the dimension of the stochastic space, $N$, is small, since the computational cost grows exponentially fast with $N$ ({\it the curse of dimensionality}). To clarify this, let $\ell \in {\mathbb Z}_+$ be a non-negative integer, called the {\it level}. Moreover, for a given index, $j \in {\mathbb Z}_+$, let $p(j)$ be an increasing function, which relates index $j$ to the polynomial degree and hence to the number of interpolation nodes. Typical functions for the polynomial degree are given by the {\it linear rule}
\begin{equation}\label{pj1}
p(j) = j,
\end{equation}
and the {\it nested rule}
\begin{equation}\label{pj2}
p(j) = 2^j \, \, \text{for} \, \, j>0, \qquad p(0)=0.
\end{equation}
In the full tensor grid, we take all polynomials of degree at most $p(\ell)$ in each direction, and $\eta = (p(\ell) + 1)^N$ grid points (nodes or knots) are therefore needed.

Alternatively, sparse grids can reduce the curse of dimensionality. They were originally introduced by Smolyak for high-dimensional quadrature and interpolation computations \cite{Smolyak:1963}.
To understand the general sparse grid construction, let ${\bf j}= (j_1, \dotsc, j_N) \in {\mathbb Z}_+^N$ be a multi-index containing non-negative integers. In each direction, $\Gamma_n$, with $n=1, \dotsc, N$ and for a non-negative index $j_n$ in ${\bf j}$, we introduce the univariate polynomial interpolation operator
$$
{\mathcal U}^{j_n}: C^0(\Gamma_n) \rightarrow {\mathbb P}_{p(j_n)}(\Gamma_n),
$$
on $p(j_n)+1$ suitable knots. Here, ${\mathbb P}_p$ is the space of univariate polynomials of degree $p$. The univariate interpolation operator applied to a stochastic function of only one random variable, say $v(y_n)$, reads
$$
{\mathcal U}^{j_n} [v(y_n)] := \sum_{k=1}^{p(j_n)+1} v(y_n^{(k)}) \, L_k (y_n).
$$
Here, $L_k (y_n)$ is the univariate Lagrange polynomial of degree $k-1$. With ${\mathcal U}^{-1} = 0$, we then define the detail operator as
$$
\Delta^{j_n} := {\mathcal U}^{j_n} - {\mathcal U}^{j_n-1}.
$$
After the introduction of a sequence of index sets, ${\mathcal I}(\ell) \subset {\mathbb Z}_+^N$, the sparse grid approximation of the solution to \eqref{wave} at level $\ell$ reads
\begin{equation}\label{sparse_formula}
u^\varepsilon(\,\cdot\,,{\bf y}) \approx {\mathcal S}_{{\mathcal I}(\ell)} [u^\varepsilon_{\text{GB}}(\,\cdot\,,{\bf y})] = \sum_{{\bf j} \in {\mathcal I}(\ell)} \bigotimes_{n=1}^N \Delta^{j_n} \, [u^\varepsilon_{\text{GB}} (\,\cdot\,,{\bf y})].
\end{equation}
Equivalently, we can rewrite the sparse approximation \eqref{sparse_formula} as
\begin{equation}\label{sparse_formula2}
{\mathcal S}_{{\mathcal I}(\ell)} [u^\varepsilon_{\text{GB}}(\,\cdot\,,{\bf y})] = \sum_{{\bf j} \in {\mathcal I}(\ell)} C({\bf j}) \, \bigotimes_{n=1}^N {\mathcal U}^{j_n} \, [u^\varepsilon_{\text{GB}} (\,\cdot\,,{\bf y})], \qquad C({\bf j}) = \sum_{\substack{{\bf i} = \{ 0, 1 \}^N \\ {\bf i} + {\bf j} \in {\mathcal I}(\ell)}} (-1)^{|{\bf i}|},
\end{equation}
where the multivariate interpolation operator reads
\begin{equation}\label{multivar_interp}
\bigotimes_{n=1}^N {\mathcal U}^{j_n} \, [u^\varepsilon_{\text{GB}} (\,\cdot\,,{\bf y})] := \sum_{k_1 = 1}^{p(j_1)+1} \dotsc \sum_{k_N = 1}^{p(j_N)+1}  u^\varepsilon_{\text{GB}}(\,\cdot\,,(y_1^{(k_1)}, \dotsc, y_N^{(k_N)})) \, \widetilde{L}_k({\bf y}),
\end{equation}
and $\widetilde{L}_k({\bf y}) = \prod_{n=1}^N L_{k_n}(y_n)$.
 This operator is given by the tensor product of $N$ univariate interpolation operators, constructed on a full tensor grid corresponding to the multi-index ${\bf j}$. The second formulation \eqref{sparse_formula2} shows that a sparse grid is a linear combination of a few tensor product grids, each with a relatively small number of grid points.

To characterize the sparse approximation operator in \eqref{sparse_formula} and \eqref{sparse_formula2} fully, we need to provide the following:
\begin{itemize}
\item[(1)] A level $\ell \in {\mathbb Z}_+$ and a function $p(j)$ given by either \eqref{pj1} or \eqref{pj2}.

\item[(2)] A sequence of sets ${\mathcal I}(\ell)$. Typical examples of index sets include:
\begin{itemize}
\item[$\circ$] Total degree index set: $\, {\mathcal I}_{\text{TD}}(\ell) = \{ {\bf j}: \, \sum_{n=1}^N j_n \le \ell \}$.

\item[$\circ$] Hyperbolic cross index set: $\, {\mathcal I}_{\text{HC}}(\ell) = \{ {\bf j}: \, \prod_{n=1}^N (j_n+1) \le \ell+1 \}$.
\end{itemize}

\item[(3)] The family of grid points to be used. Typical choices include:
\begin{itemize}
\item[$\circ$] Gauss abscissas, that are the zeros of $\rho$-orthogonal polynomials; see e.g. \cite{Trefethen:08}.

\item[$\circ$] Clenshaw-Curtis abscissas, that are the extrema of Chebyshev polynomials; see e.g. \cite{Trefethen:08}.

\item[$\circ$] Leja abscissas; see e.g. \cite{Leja:14}.
\end{itemize}
\end{itemize}

\begin{remark} (Full tensor and Smolyak grids) The full tensor grid corresponds to the full tensor index set, ${\mathcal I}_{\text{FT}}(\ell) = \{ {\bf j}: \, \max_{n} j_n \le \ell \}$. The Smolyak sparse grid is a particular type of sparse grid, where we use the nested rule \eqref{pj2} and the total degree index set, based on either Gauss or Clenshaw-Curtis abscissas. In the latter case, we obtain a nested grid, in which grids in different levels are embedded.
\end{remark}

In many practical applications, the main goal is the computation of the statistical moments of the solution or some QoIs. In such cases, we can directly compute the statistical moments using Gauss or Clenshaw-Curtis quadrature formulas, without explicitly constructing the solution or other quantities by sparse approximation formulas such as \eqref{sparse_formula2}. Suppose we want to compute the statistics of an operator applied to the solution, say $F(u^{\varepsilon})$. This can for instance be the solution itself, the QoI \eqref{Q}, or different powers of these quantities when computing higher moments. We write
$$
{\mathbb E}[F(u^{\varepsilon}(\,\cdot\,,{\bf y}))] \approx {\mathbb E}[
{\mathcal S}_{{\mathcal I}(\ell)} [F(u^\varepsilon_{\text{GB}}(\,\cdot\,,{\bf y}))]
]
= \int_{\Gamma}  {\mathcal S}_{{\mathcal I}(\ell)} [F(u^\varepsilon_{\text{GB}}(\,\cdot\,,{\bf y}))]  \, \rho({\bf y}) \, d{\bf y}
\approx \sum_{k=1}^{\eta} \theta_k \, F(u^\varepsilon_{\text{GB}}(\,\cdot\,,{\bf y}^{(k)})),
$$
where, by \eqref{sparse_formula2}-\eqref{multivar_interp}, the quadrature weights read
$$
\theta_k = c_k \int_{\Gamma} \widetilde{L}_k({\bf y}) \, \rho({\bf y}) \, d{\bf y}.
$$
Here, each global index, $k=1, \dotsc, \eta$, corresponds to a local multi-index $[k_1, \dotsc, k_N]$ in the formulas \eqref{sparse_formula2}-\eqref{multivar_interp}. Moreover, the coefficients $c_k$ correspond to the coefficients $C$ in \eqref{sparse_formula2}.

\medskip
\noindent
{\bf Convergence of stochastic collocation.}
{It is well known that the rate of convergence of stochastic collocation for a stochastic function, say $v({\bf y})$, in general depends on the stochastic regularity of the function, i.e., the regularity of the mapping $v: \Gamma \rightarrow {\mathbb R}$. Fast convergence is attained in the presence of high stochastic regularity. For instance, suppose that $v({\bf y})$ is continuous and admits an analytic extension in the complex region
$$
\tilde\Gamma_\tau = \left\{
\bold z=(z_1,\ldots,z_N)\in{\mathbb C}^N\,|\, \text{\rm for some $j$: dist}(z_j,\Gamma_j)\leq \tau\ {\rm and}\ z_k\in \Gamma_k\ {\rm for}\ k\neq j
\right\},
$$
whose size is characterized by the radius of analyticity, $\tau >0$.
Then, the maximum error, $\epsilon_{\text{max}}$, in the sparse approximation, ${\mathcal S}_{{\mathcal I}(\ell)} [v({\bf y})]$, in \eqref{sparse_formula2} on the Smolyak grid based on Clenshaw-Curtis abcissas satisfies (see Theorem 3.10 in \cite{NTW:08})
$$
\epsilon_{\text{max}} \le C(\tau,N,\Gamma) \, M(v,\tau,\Gamma)\, \eta^{- \sigma / (1 + \log(2 \, N))},\qquad M(v,\tau,\Gamma) = \max_{\bold z\in\tilde\Gamma_\tau}|v(\bold z)|,
$$
where $\sigma > 0$ is directly proportional to the radius of analyticity, $\tau$, and $C$, $M$ are independent of $\eta$. Therefore, the larger the radius of analyticity of the map, the faster the convergence rate in $\eta$; see also \cite{NTT:14,BNTT:14,BNTT:14b}.
In the case of high-frequency waves with highly oscillatory solutions, $\tau$, $\sigma$ and $M$ in general all depend on the
wavelength $\varepsilon$. To maintain fast convergence, it is therefore important that $M$ and $1/\tau$ remain bounded as $\varepsilon\to 0$.
Since the size of $M$ depends on the size of the $\bold y$-derivatives of $v$ on $\Gamma$,
we need uniform bounds for ${\bf y}$-derivatives independent of $\varepsilon$ to obtain fast convergence for all $\varepsilon$. In other words, we require high regularity of the mapping ${\mathcal Q}^\varepsilon: \Gamma \rightarrow {\mathbb R}$, uniformly in $\varepsilon$. This poses extra challenges in employing stochastic spectral techniques for high-frequency waves. The stochastic regularity and the convergence of stochastic collocation for high frequency waves is discussed in more detail in Section \ref{sec4}.
}

\section{Stochastic Regularity of High-Frequency Waves}
\label{sec4}

For fast convergence, regular dependence of the quantity of interest on the input random variables is required, as discussed in
Section \ref{sec:3_2}, and, ideally,
${\mathcal Q}^\varepsilon({\bold y})\in C^\infty(\Gamma)$.
For high-frequency problems, we would like to have an even stronger bound,
namely
        \begin{equation}\label{est}
       \sup_{{\bold y}\in\Gamma} \left|\frac{\partial^{\bold k} \mathcal Q^\varepsilon({\bold y})}{\partial {\bold y}^{\bold k}}\right| \leq C_{\bold k},   \quad \forall {\bold k} \in \mathbb Z_+^N,
        \end{equation}
        where the constants $C_{\bold k}$ can be taken independent
        of the wavelength, $\varepsilon$. If this is not true,
        the spectral convergence rate of
        a stochastic collocation method could deteriorate
        for small wavelengths, $\varepsilon$; with the bound (\ref{est}),
        we expect that there is a rate that is uniform in $\varepsilon$.

        We note that in general
        $u^\varepsilon(t,{\bf x},{\bf y})$ will be
        oscillating with period $\sim\varepsilon$ in both ${\bf x}$ and ${\bf y}$.
        If the corresponding  quantity of interest
        would inherit this property, a bound like (\ref{est}) could not hold.
        Nevertheless, in this section, we conjecture that
        for a Gaussian beam superposition approximation of $u^\varepsilon$
        with the initial data given by \eqref{waveinit}
        and \eqref{phi0ma0} respectively,
        the
        estimate \eqref{est} holds.
\begin{conjecture}
Bound (\ref{est}) holds for the following quantity of interest computed with the Gaussian beam approximation (\ref{GBsup})
$$
 {\mathcal Q}^\varepsilon({\bold y}) = \int_{{\mathbb R}^n} |u_{GB}^\varepsilon(T,{\bf x},{\bf y})|^2\psi({\bf x})d{\bf x},
$$
if $\psi\in C^\infty_c({\mathbb R}^n)$.
\end{conjecture}

  This is based on formal theoretical arguments and several numerical
  experiments, which are presented in the subsections below.
The conjecture will be proved in an upcoming work.

\subsection{Motivation away from caustics}
\label{sec:4_1}

In this section, we give a non-rigorous argument for why \eqref{est} should
hold, at least in the case when there are no caustics in the support
of the quantity of interest.
We thus consider a high-frequency solution
for \eqref{wave} and we assume $(t_0,{\bold x}_0,{\bold y}_0)\in[0,T]\times\Real^n\times\Gamma$ is a non-caustic point.
By the theory of Maslov \cite{maslov:65,maslov:81} and Fourier
Integral Operators \cite{hormander83,Hormander,Duistermaat:74}
there is then at least a neighborhood, $U$, of $(t_0,{\bf x}_0,{\bf y}_0)$ in which
the solution has the form
\begin{equation}\label{multiphase0}
u^\varepsilon(t,{\bold x},{\bold y}) = \sum_{j=1}^{K} A_j(t,{\bold x},{\bold y})e^{i\phi_j(t,{\bold x},{\bold y})/\varepsilon+im_j\pi/4} + O(\varepsilon),\qquad \forall (t,{\bold x},{\bold y})\in U,
\end{equation}
for some integers $K$, which
 represents the number of waves passing through
$(t_0,{\bold x}_0)$, and
$m_j$, which is the Keller--Maslov index
of wave $j$, i.e.,\mbox{} the
number of times it has passed through a caustic.
Note that in general $K>1$ when $t_0>0$,
even if the initial data is of the single wave form \eqref{waveinit}
since new terms appear when
caustics develop; see \cite{SparberEtAl:03} for further explanations.
The amplitudes
$A_j(t,{\bf x},{\bf y})$ and phases
$\phi_j(t,{\bf x},{\bf y})$ satisfy the geometrical optics
equations, i.e.,\mbox{} the transport and eikonal equations \eqref{eiko}, and \eqref{transp} respectively.
They are independent of $\varepsilon$ and
smooth in all variables,
given the assumptions of a smooth wave speed and
initial data.
Moreover, the phases, $\phi_j$, are related to the initial phase, $\Phi_0$,
and to ${\bold p}(t,{\bf y};{\bf z})$ via
the ray function, ${\bf q}(t,{\bf y};{\bf z})$, as
follows. For each $(t,{\bf x},{\bf y})\in U$ and each $1\leq j\leq K$,
there is a ${\bf z}_j\in Z_j(t,{\bf y})$ such that
\begin{equation}\label{multiphase}
{\bf x}={\bf q}(t,{\bf y};{\bf z}_j), \qquad
  \phi_j(t,{\bf x},{\bf y})=\Phi_0({\bf z}_j,{\bf y}),\qquad
  \nabla\phi_j(t,{\bf x},{\bf y})={\bold p}(t,{\bf y};{\bf z}_j),
\end{equation}
where $Z_j(t,{\bf y})$ are disjoint open subsets of supp($A_0(\,\cdot\,,{\bf y})$).


Based on this, let us assume that in the
support of the QoI test function, $\psi$,
there are no caustic points at $t=T$,
and that the solution satisfies \eqref{multiphase0} and
\eqref{multiphase}.
The quantity of interest \eqref{Q} can then be written
\begin{align*}
 {\mathcal Q}^\varepsilon({\bf y}) &= \int_{\Real^n} |u^\varepsilon(T,{\bf x},{\bf y})|^2\psi({\bf x})d{\bf x}
 \\
 &\approx
 \int_{\Real^n}
 \left|\sum_{j=1}^{K} A_j(T,{\bf x},{\bf y})e^{i\phi_j(T,{\bf x},{\bf y})/\varepsilon+im_j\pi/4}\right|^2
\psi({\bf x})d{\bf x}
\\
&=\Re \sum_{i=1}^{K}\sum_{j=1}^{K}\widetilde{Q}^\varepsilon_{ij}({\bf y})
e^{i(m_i-m_j)\pi/4},
\end{align*}
where $\Re$ denotes the real part and
\begin{equation}\label{Qij}
\widetilde{Q}^\varepsilon_{ij}({\bf y}) =
\int_{\Real^n}
 A_i(T,{\bf x},{\bf y})A_j(T,{\bf x},{\bf y})
 e^{i(\phi_i(T,{\bf x},{\bf y})-\phi_j(T,{\bf x},{\bf y}))/\varepsilon}
\psi({\bf x}) d{\bf x}.
\end{equation}
We note that \eqref{est} holds if the same estimate holds for
each $\widetilde{Q}_{ij}^\varepsilon(\bold y)$.
Clearly the diagonal terms,
$$
\widetilde{Q}^\varepsilon_{jj}({\bf y}) =
\int_{\Real^n}
  A_j(T,{\bf x},{\bf y})^2
\psi({\bf x}) d{\bf x},
$$
are smooth and independent of $\varepsilon$, hence satisfying
\eqref{est}. For the off-diagonal entries, $i\neq j$, however,
we have for a multi-index, ${\bold k}\in\mathbb Z_+^N$,
$$
\frac{\partial^{\bold k} \widetilde{Q}_{ij}^{\varepsilon}({\bf y})}{\partial{\bf y}^{\bold k}}
=\sum_{\ell=1}^{|{\bold k}|} {\frac{1}{\varepsilon^{\ell}}}
\int_{\Real^n}  r_{\ell,{\bold k}}({\bf x},{\bf y})\psi({\bf x}) \exp\left(i\frac{\phi_i(T,{\bf x},{\bf y})-\phi_j(T,{\bf x},{\bf y})}{\varepsilon}\right) d{\bf x},\qquad  i\neq j,
$$
where $r_{\ell,{\bold k}}$ are smooth, $\varepsilon$-independent functions built from sums and products of $\bold y$-derivatives
of
$A_i$, $A_j$, $\phi_i$ and $\phi_j$ evaluated at $t=T$.
Because of the pre-factors, $\varepsilon^{-\ell}$, it follows that we cannot
bound $\widetilde{Q}_{ij}^{\varepsilon}$ as in \eqref{est} unless
\begin{equation}\label{eq:oscbound}
\int_{\Real^n}  r_{\ell,{\bold k}}({\bf x},{\bf y})\psi({\bf x}) \exp\left(i\frac{\phi_i(T,{\bf x},{\bf y})-\phi_j(T,{\bf x},{\bf y})}{\varepsilon}\right) d{\bf x}=O(\varepsilon^{m}),
\qquad \forall m\in{\mathbb Z}_+.
\end{equation}
To verify this, we use stationary phase arguments.
More specifically, the following theorem formalizes the dependence of oscillatory integrals on critical points of the phase; see \cite{hormander83}.
\begin{theorem}[Principle of (non)-stationary phase]
Let $\psi \in C_c^\infty(\mathbb R^n)$ and $\varphi \in C^\infty(\mathbb R^n; \mathbb R)$ such that $\nabla \varphi \ne 0$ on $supp(\psi)$. Then, for all $m\in \mathbb Z_+$,
there exist constants, $C_m$, independent of $\varepsilon$, such that
        \[
           \left| \int_{\mathbb R^n} e^{i \varphi({\bf x})/\varepsilon}\: \psi({\bf x}) d{\bf x}\right| \leq C_M \varepsilon^{m}.
        \]
\end{theorem}\label{th:stat_phase}
When applying this theorem to the integral in
\eqref{eq:oscbound}, we thus obtain the
required estimate provided two conditions are fulfilled.
First, the QoI test function should be {\em smooth} with compact support, i.e.,
  $\psi\in C_c^{\infty}(\Real^n)$.
  Since the functions, $r_{\ell,{\bold k}}$, are
  smooth, the theorem applies. However, we can, for instance, not expect \eqref{est} to hold if we define the
  QoI just as an integral over a domain, $D\subset\Real^n$,
  \begin{equation}\label{nonsmoothQ}
     {\mathcal Q}^\varepsilon(y) = \int_{D} |u^\varepsilon(T,{\bf x},{\bf y})|^2d{\bf x},
  \end{equation}
which corresponds to choosing $\psi$ as the (non-smooth) characteristic function of $D$.
Second, there should be no stationary points in the support of $\psi$.
However, under our assumptions, this is not possible. Indeed, suppose
there is a stationary point, i.e.,
$$
  \nabla\phi_i(T,{\bf x},{\bf y})=\nabla\phi_j(T,{\bf x},{\bf y})
$$
for some ${\bold x}\in$ supp($\psi$). Then, by \eqref{multiphase},
there are two geometrical optics rays starting at ${\bf z}_i
\in Z_i(T,{\bf y})$
and ${\bf z}_j\in Z_j(T,{\bf y})$, respectively,
such that
\begin{align*}
{\bf q}(T,{\bf y};{\bf z}_i)&= {\bf q}(T,{\bf y};{\bf z}_j)={\bf x},\\
{\bf p}(T,{\bf y};{\bf z}_i) &=\nabla\phi_i(T,{\bf x},{\bf y})
=\nabla\phi_j(T,{\bf x},{\bf y})=
{\bf p}(T,{\bf y};{\bf z}_j).
\end{align*}
But, by uniqueness of the solution to the Hamiltonian system \eqref{qp}, we must therefore have ${\bf z}_i={\bf z}_j$. Hence, the two rays are the same
and since the sets $\{Z_j(t,{\bf y})\}$ are disjoint, we
get $i=j$, a contradiction.

In conclusion, the formal arguments above suggest that
if the QoI test function is smooth,
then, at least away from caustics, the bound \eqref{est}  holds.

\subsection{Numerical justification and preliminaries}
\label{sec:4_2}

In this section, we consider a number of problems with high-frequency solutions to \eqref{wave}, to justify the
bound \eqref{est}. We recognize that the theoretical arguments
for \eqref{est}
in the previous section were only
formal.
In general high-frequency solutions, there are caustics,
which we did not account for.
Furthermore, we did not consider the regularity of the
approximation error, $O(\varepsilon)$, in \eqref{multiphase0}.

Here we will show three numerical examples to demonstrate that
the quantity of interest, ${\mathcal Q}^\varepsilon$,
can also be smooth in more general cases.
First, we show a simple one-dimensional example with constant coefficients,
where the arguments in Section \ref{sec:4_1} are essentially
rigorous; in one dimension, there are no caustics, and for constant
coefficients, geometrical optics is exact. Our second example
is still one-dimensional with no caustics,
but the speed of propagation varies in space. Geometrical optics is therefore
accurate only up to $O(\varepsilon)$. The last example is a two-dimensional
problem where there are caustics in the solution. As we will see,
for all these cases, the QoI and its derivatives
remain non-oscillatory when the QoI test function is smooth.

All three examples use an initial condition
of the single wave type,
        \begin{equation}\label{init}
            u^\varepsilon(0,{\bold x},{\bold y}) = A_0(\mathbf x,{\bold y}) e^{i \Phi_0(\mathbf x,{\bold y})/\varepsilon},
        \end{equation}
where the amplitude, $A_0$, is a sum of two pulses moving towards
each other. The integration via the
trapezoidal rule \eqref{trapets} is carried out to obtain the approximated Gaussian beam solution, $u_{GB}^\varepsilon$. The spacing between the starting points is set to $\Delta z = \sqrt{\varepsilon}$ in all examples.
The QoI test function in \eqref{Q} is based on
the smooth function
\begin{equation}\label{testfun}
        \widetilde{\psi}(\mathbf x) = \left\{
            \begin{array}{ll}
                e^{-\frac{|\mathbf x|^2}{1-|\mathbf x|^2}}, & |\mathbf x| < 1,\\
                0, & \text{otherwise}.
            \end{array}
        \right.
    \end{equation}
To show the regularity of ${\mathcal Q}^\varepsilon$
at high frequencies, it is computed for
a sequence of small wavelengths, $\varepsilon = (1/40,1/80,1/160)$.
The final time in \eqref{Q} is fixed at $T = 1$ for all three examples below.

The integration over variable $\mathbf x$ in \eqref{Q} is again carried out by the trapezoidal rule with a spatial step that is uniform in all dimensions, $\Delta x = \frac{2 \pi \varepsilon}{10}$, unless stated otherwise. This corresponds to ten points per wavelength, which gives
a very high accuracy when $\psi\in C^\infty_c(\Real^n)$ due
to the spectral convergence of the trapezoidal rule for such functions.

\subsubsection{Example 1: Constant speed of propagation in 1D}

\begin{figure}
\centering
    \begin{minipage}{0.43\textwidth}
            \subfigure[$t=0.25$]{
\includegraphics[width=\textwidth]{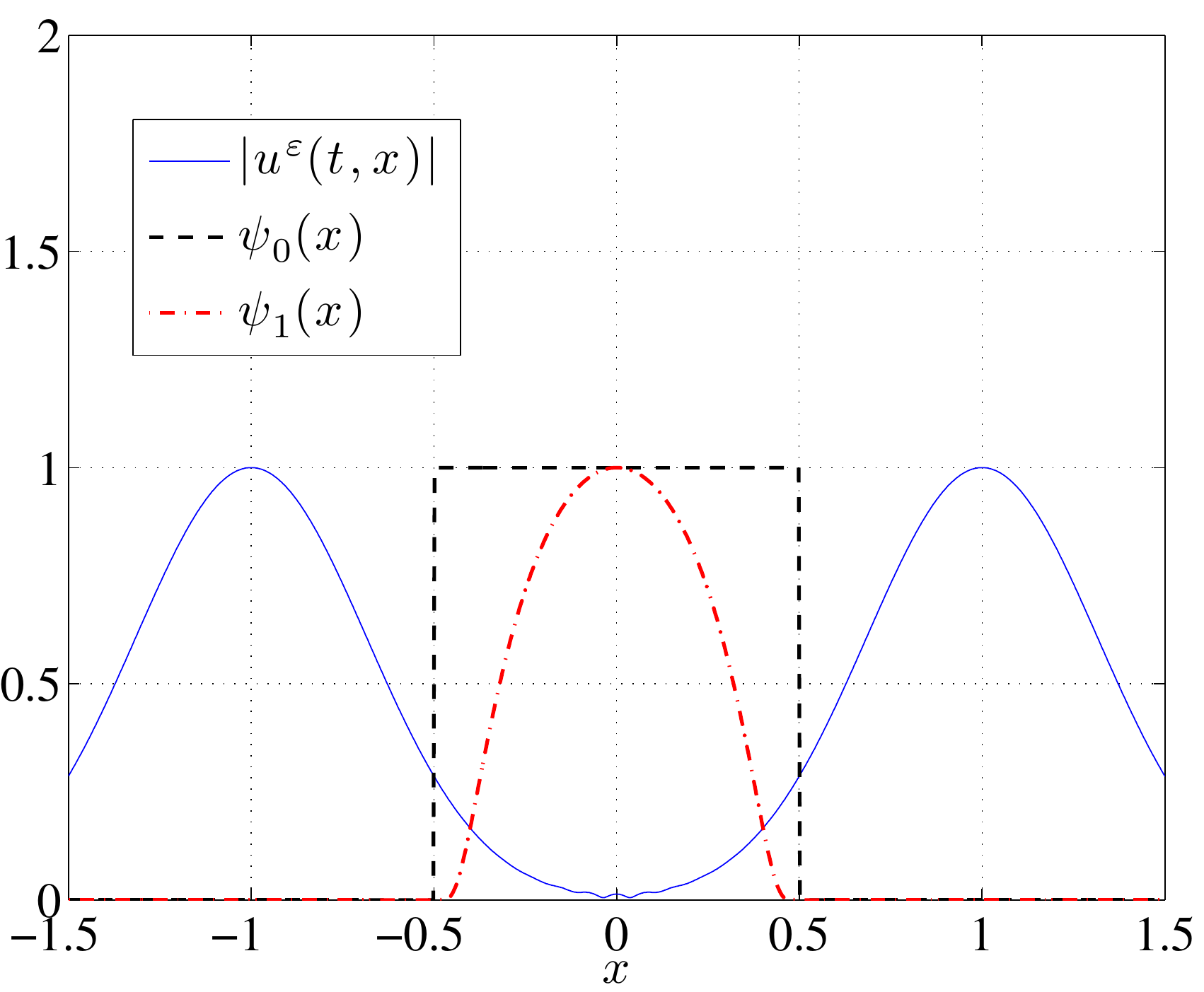}}
    \end{minipage}
    \begin{minipage}{0.43\textwidth}
            \subfigure[$t=0.5$]{
\includegraphics[width=\textwidth]{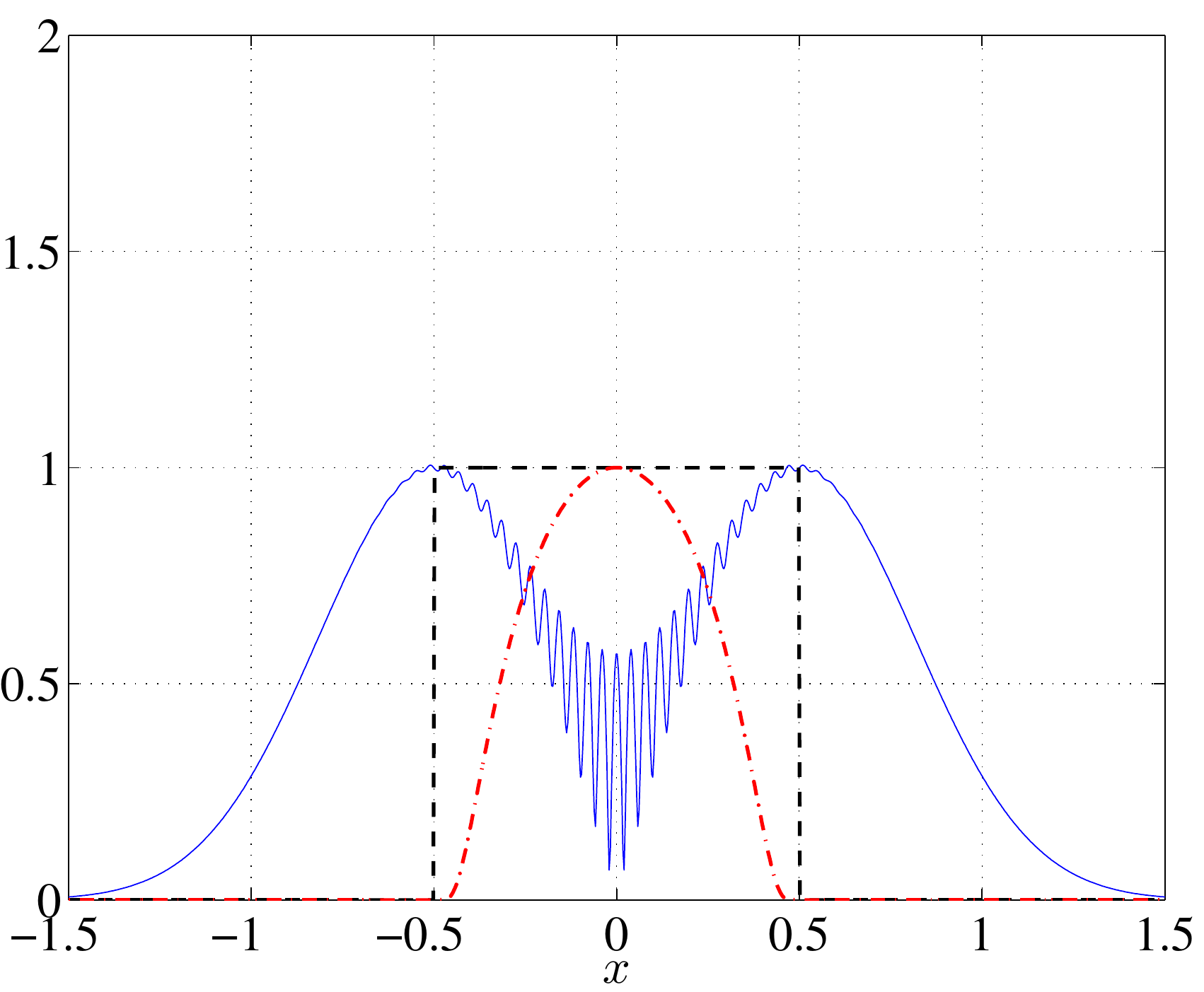}}
    \end{minipage}
    \begin{minipage}{0.43\textwidth}
            \subfigure[$t=0.75$]{
\includegraphics[width=\textwidth]{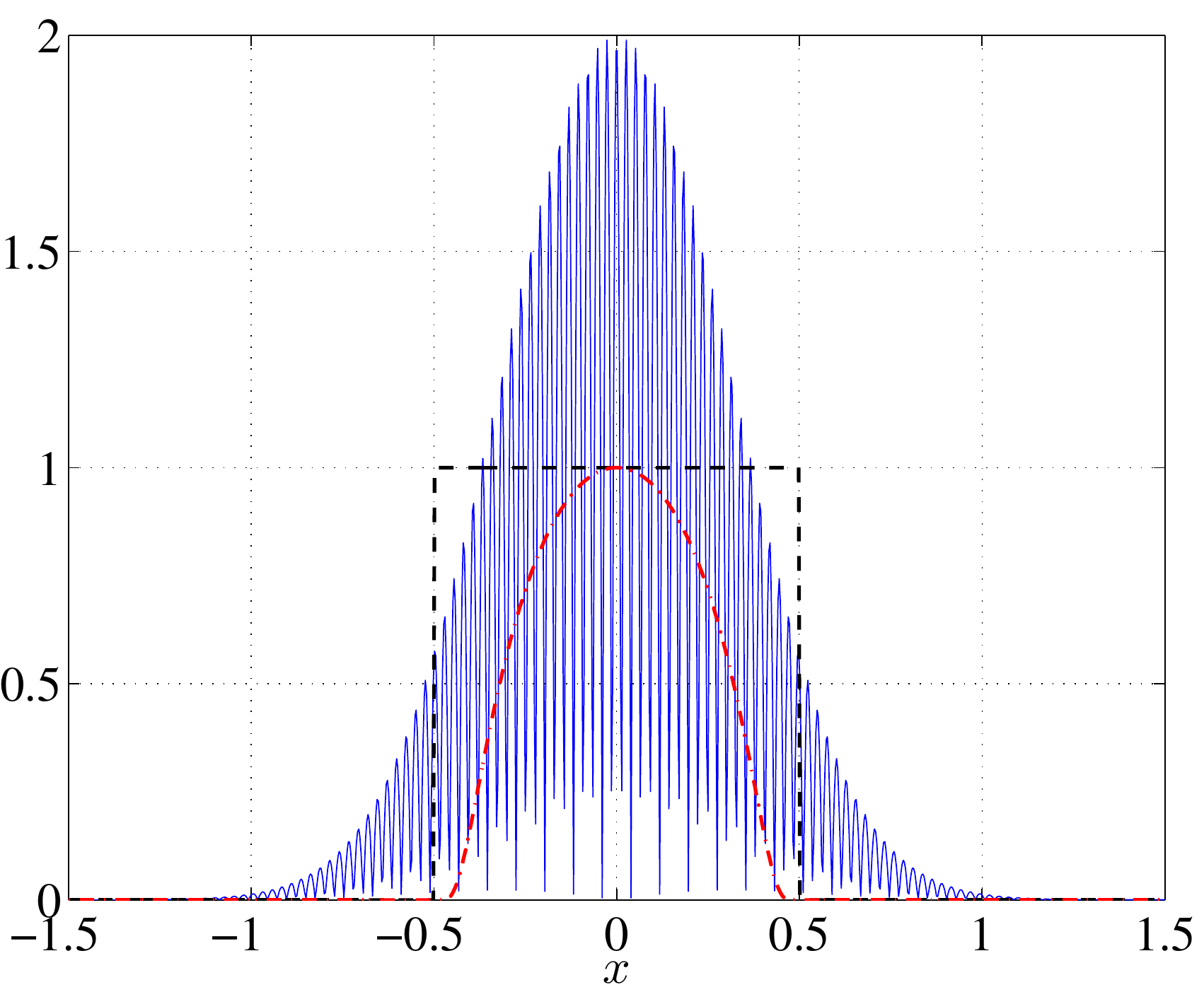}}
    \end{minipage}
    \begin{minipage}{0.43\textwidth}
            \subfigure[$t=1$]{
\includegraphics[width=\textwidth]{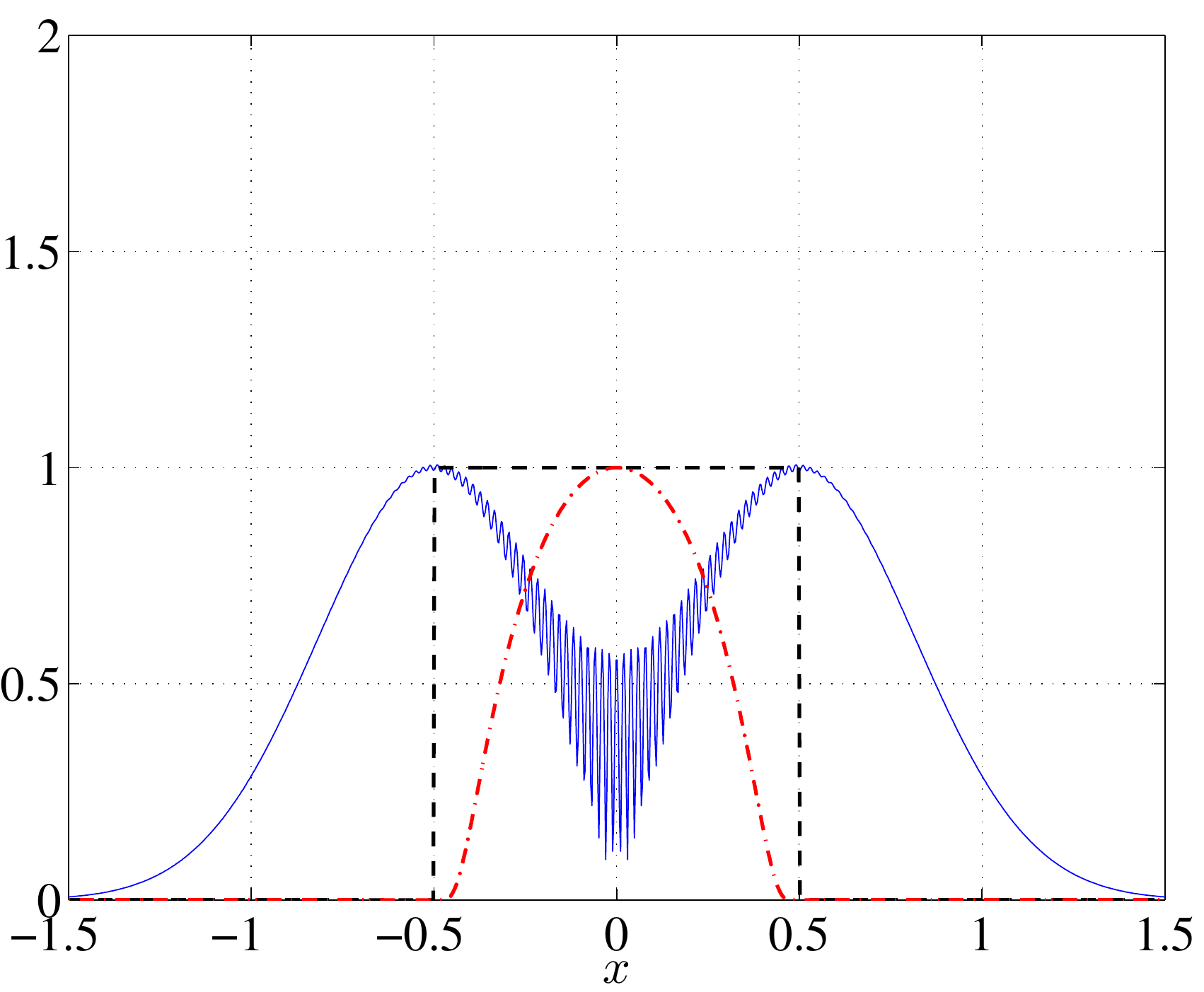}}
    \end{minipage}
\caption{Example 1. Absolute value of solution for various times, $t$, when $y = 2$ and $\varepsilon=1/40$.
The two QoI functions, $\psi_0$ and $\psi_1$, are overlaid.}
\label{fig:Ex1}
\end{figure}

We start with the simple case of the
one-dimensional wave equation \eqref{wave} with constant speed of propagation,
where we have an explicit expression for the
solution $u(t,x,y)$.
We investigate a case when two high-frequency pulses are moving towards each other,
$$
  A_0(x) = g(x-s_1) + g(x-s_2), \quad g(x) = e^{-5x^2},\qquad \Phi_0(x) = x^2,
$$
where $s_1 = -s_2 = 1.5$.
The initial data has no uncertainty, i.e.,\mbox{} it does not depend on $y$. However, the constant (in $x$) speed is uncertain, $c=c(y)$.
The initial time-derivative of $u$ is taken such that the pulses
propagate towards each other,
    \begin{align*}
        u^\varepsilon_t(0,x,y)
        &= -c(y)\left(g'(x-s_1) + \frac{i \Phi'_0(x)}{\varepsilon} \: g(x-s_1)\right) e^{i\Phi_0(x)/\varepsilon} +\\
         &+ c(y) \left(g'(x-s_2) + \frac{i \Phi'_0(x)}{\varepsilon} \: g(x-s_2)\right) e^{i\Phi_0(x)/\varepsilon}.
    \end{align*}
    Moreover, we let $c(y)=y>0$.
    The solution is then given by d'Alembert as
\begin{equation}\label{eq:solex1}
    u^\varepsilon(t,x,y) =
    g(x-s_1-yt) e^{i\Phi_0(x-yt)/\varepsilon}+
    g(x-s_2+yt) e^{i\Phi_0(x+yt)/\varepsilon}.
\end{equation}
The absolute value of a particular solution with $\varepsilon = 1/40$ and $y = 2$ is shown at various times
in Figure \ref{fig:Ex1}.

We consider two different quantities of interest,
$$
{\mathcal Q}^\varepsilon_j(y) = \int |u^\varepsilon(T,x,y)|^2\psi_j(x) dx,\qquad
T=1,\quad j=0,1,
$$
where $\psi_1(x)=\widetilde{\psi}(2x)$ is smooth, based on \eqref{testfun},
and $\psi_0(x)$ is the (non-smooth)
characteristic function for the interval $[-1/2,1/2]$, which reduces
the form of ${\mathcal Q}^\varepsilon_0$
to \eqref{nonsmoothQ}.
The QoI test functions, $\psi_j(x)$, are plotted together with the solution in Figure \ref{fig:Ex1}.

\begin{figure}
            \centering
    \begin{minipage}{0.32\textwidth}
            \subfigure[${\mathcal Q}_0^\varepsilon$]{
                    \includegraphics[width=\textwidth]{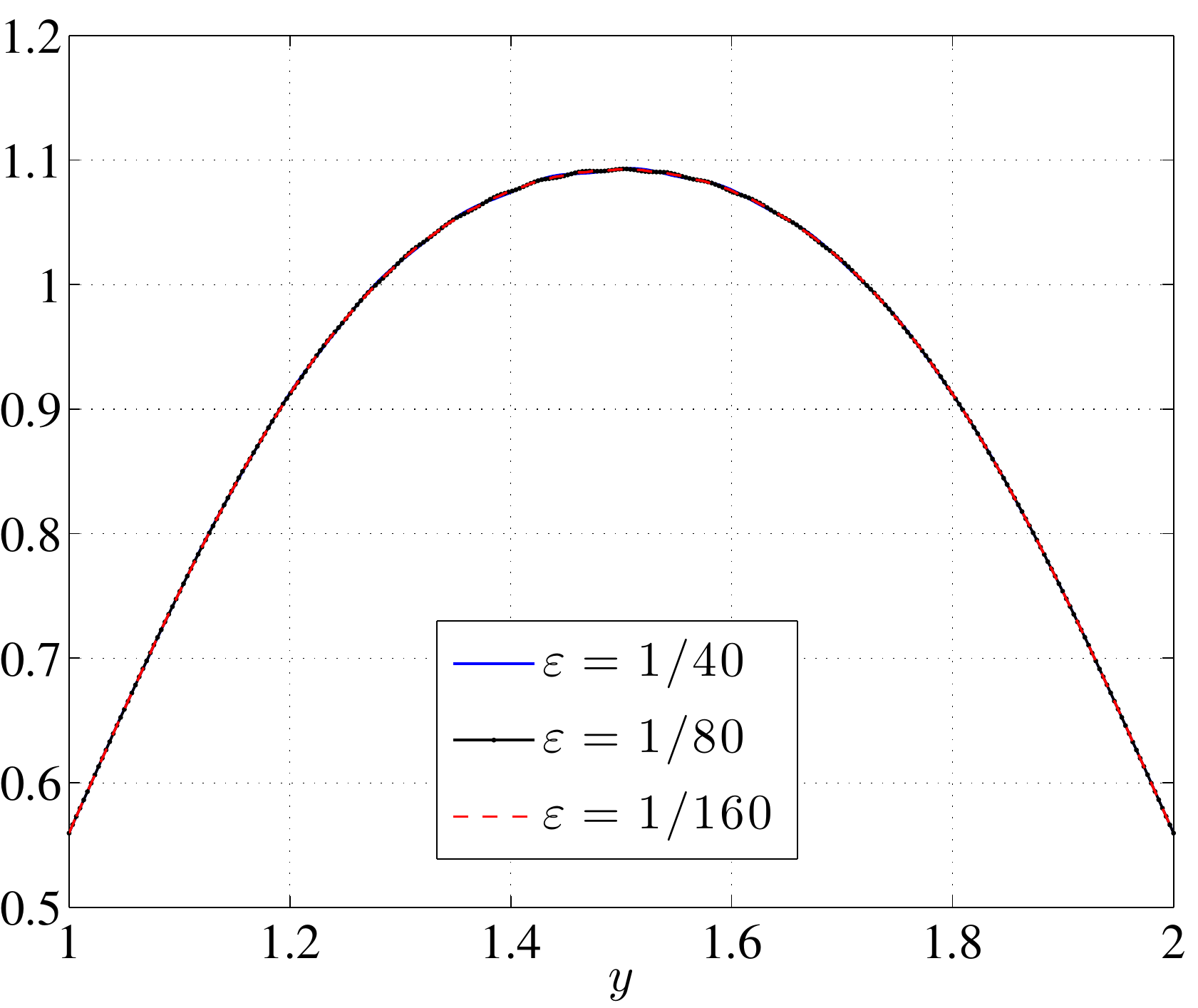}}
    \end{minipage}
    \begin{minipage}{0.32\textwidth}
            \subfigure[$\frac{d}{dy}{\mathcal Q}_0^\varepsilon$]{
                    \includegraphics[width=\textwidth]{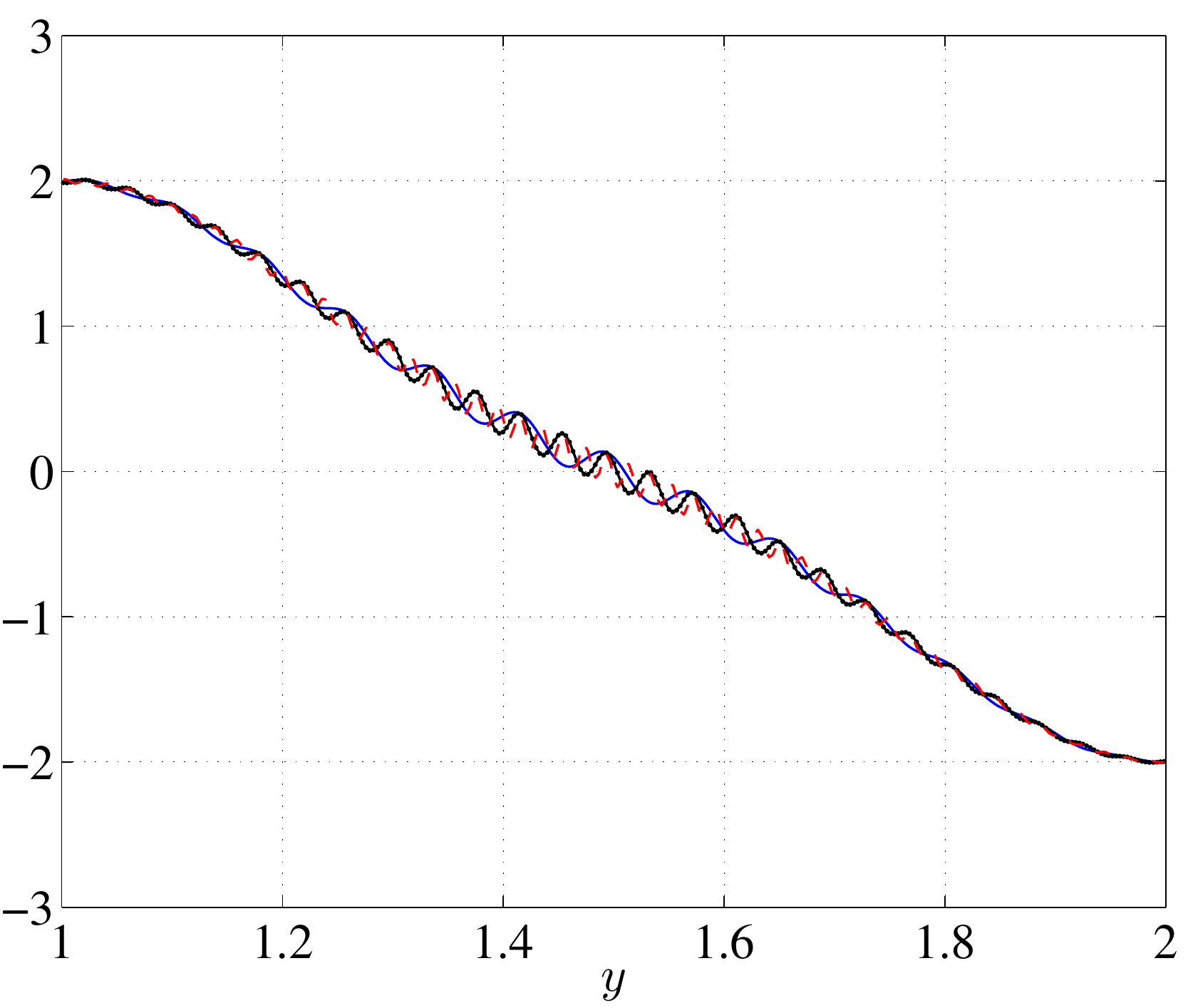}}
    \end{minipage}
    \begin{minipage}{0.32\textwidth}
            \subfigure[$\frac{d^2}{dy^2}{\mathcal Q}_0^\varepsilon$]{
                    \includegraphics[width=\textwidth]{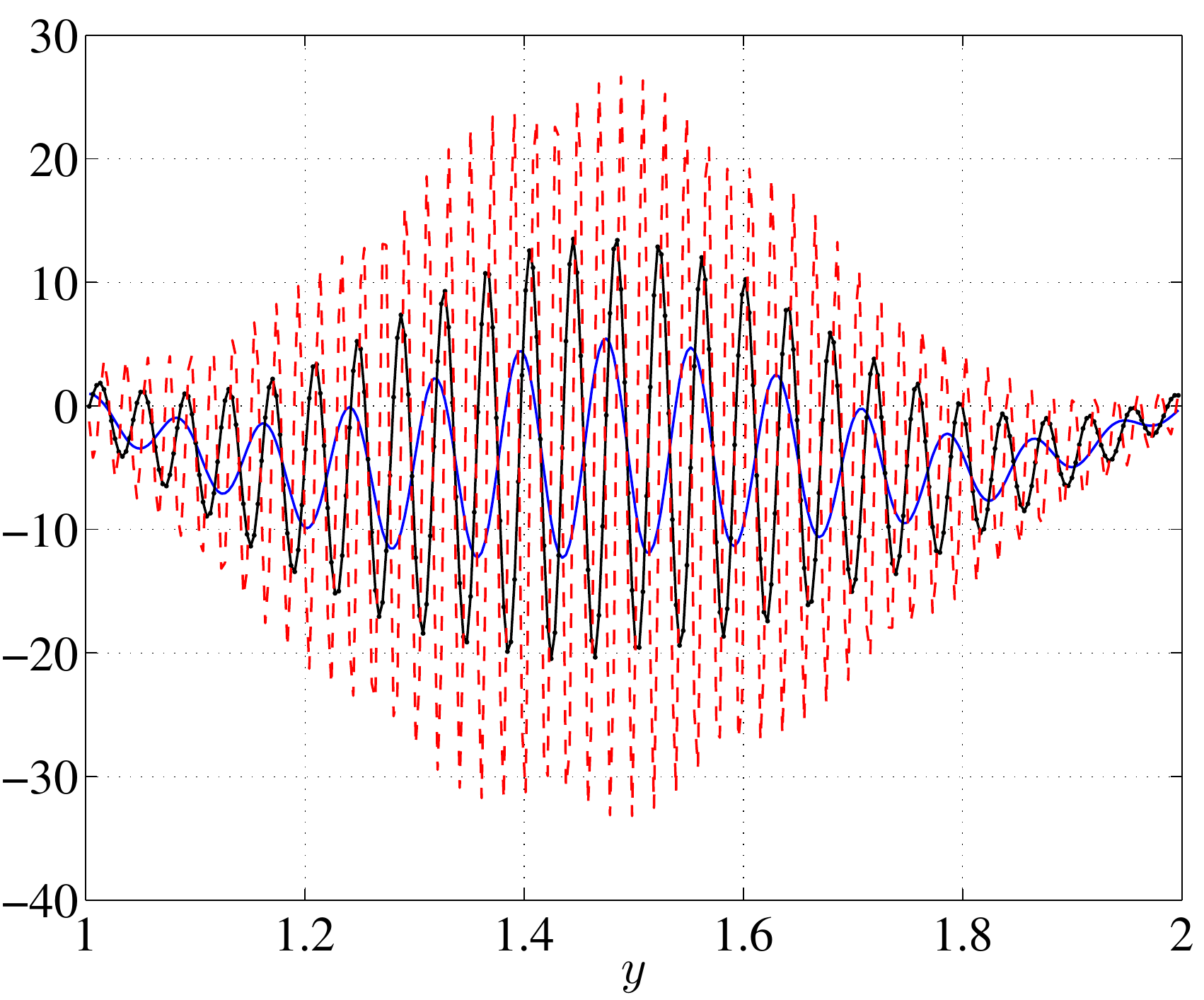}}
    \end{minipage}
    \begin{minipage}{0.32\textwidth}
            \subfigure[${\mathcal Q}_1^\varepsilon$]{
                    \includegraphics[width=\textwidth]{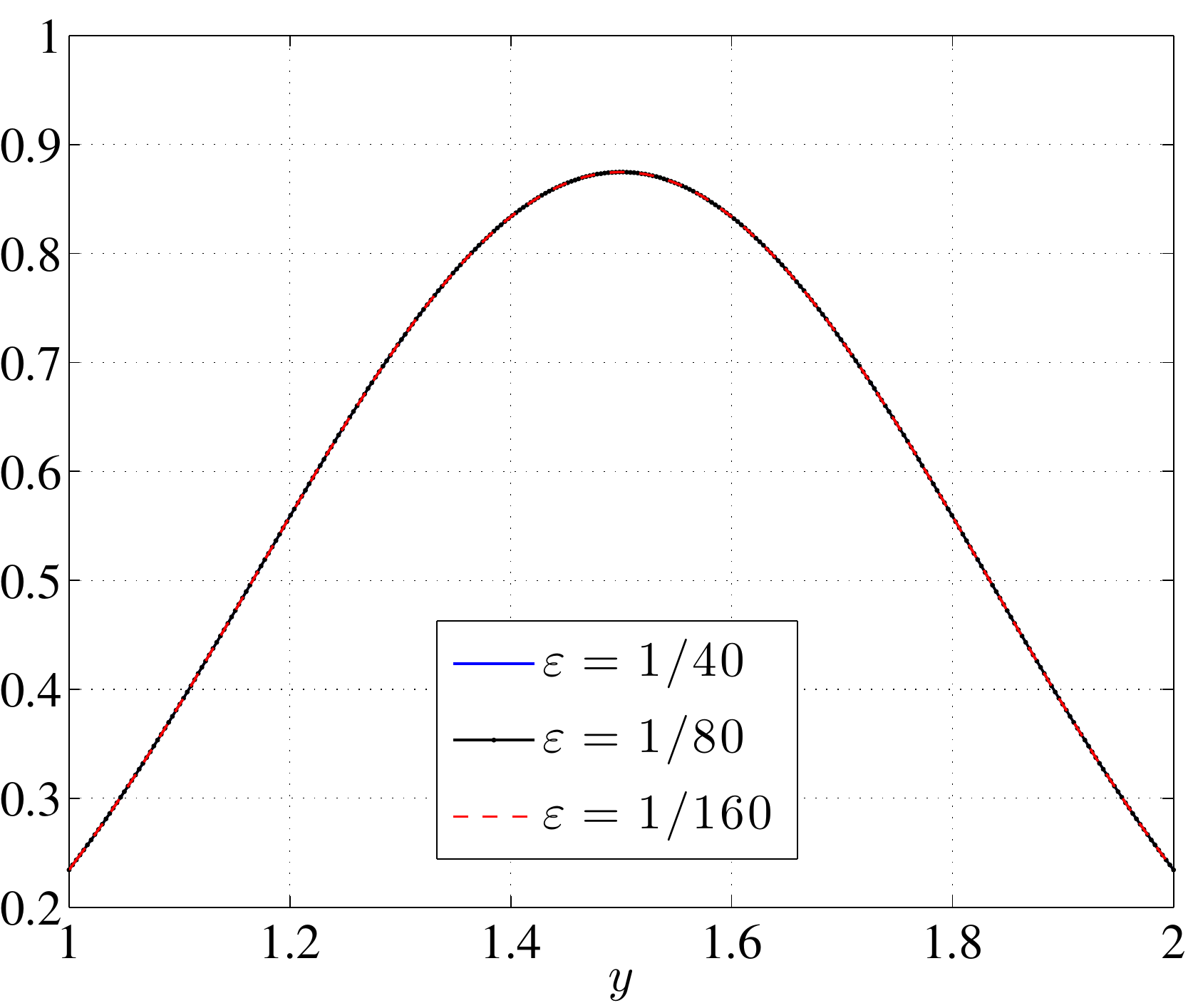}}
    \end{minipage}
    \begin{minipage}{0.32\textwidth}
            \subfigure[$\frac{d}{dy}{\mathcal Q}_1^\varepsilon$]{
                    \includegraphics[width=\textwidth]{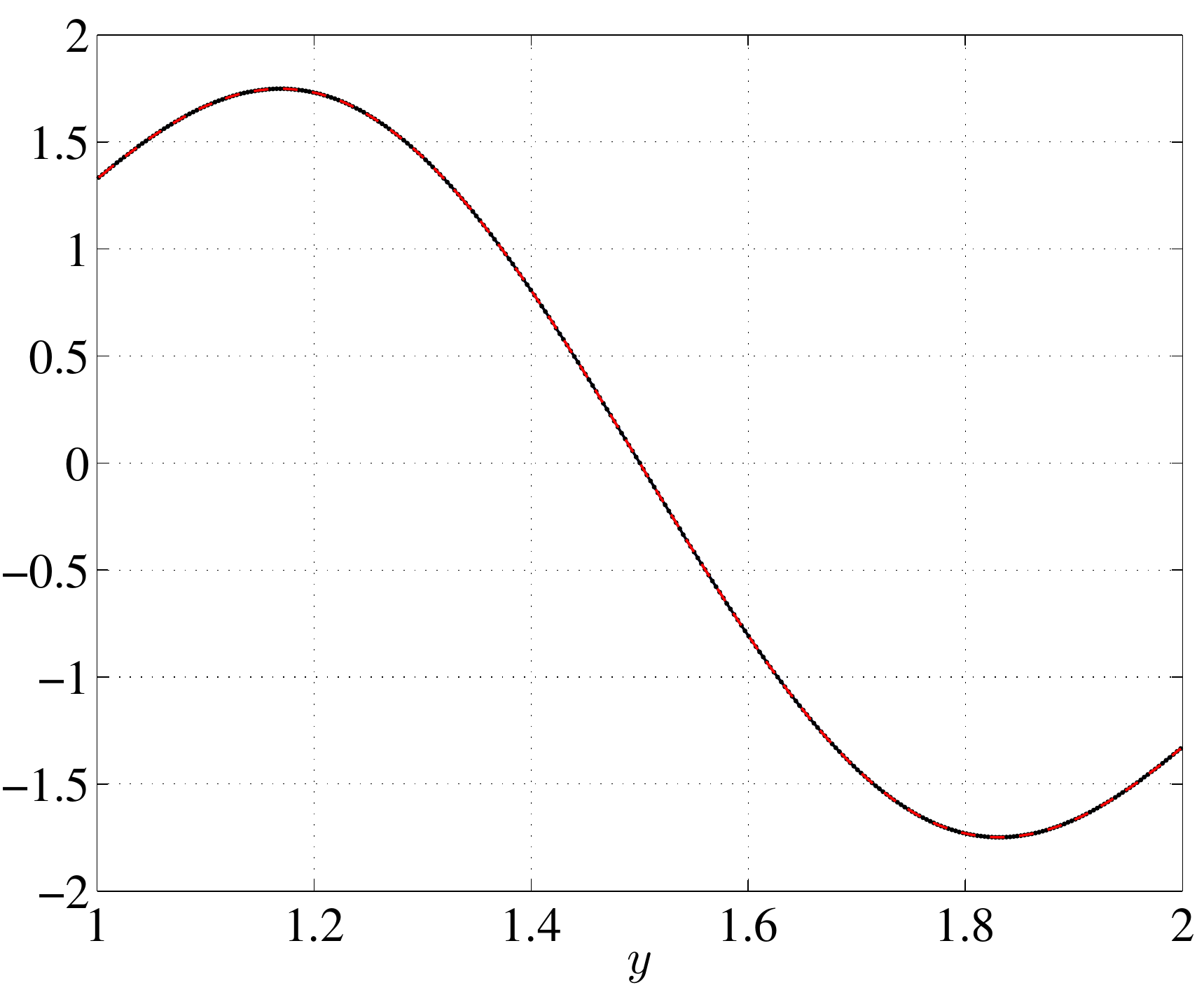}}
    \end{minipage}
    \begin{minipage}{0.32\textwidth}
            \subfigure[$\frac{d^2}{dy^2}{\mathcal Q}_1^\varepsilon$]{
                    \includegraphics[width=\textwidth]{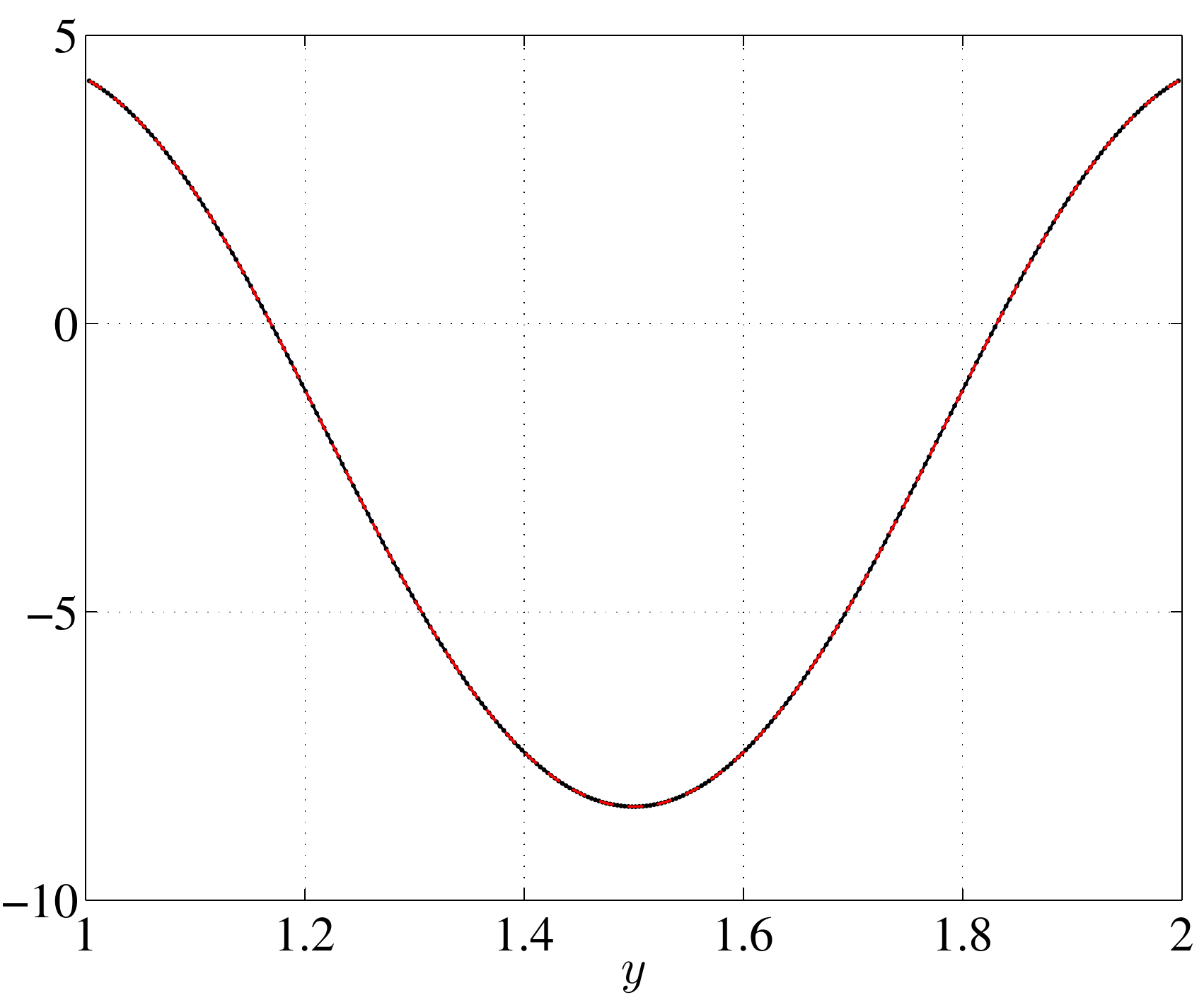}}
    \end{minipage}
            \caption{Example 1:
            QoIs,
${\mathcal Q}^\varepsilon_0$ (top row) and
${\mathcal Q}^\varepsilon_1$ (bottom row),
and their first and second derivatives,
for different wavelengths, $\varepsilon$.
For ${\mathcal Q}^\varepsilon_1$ the curves
are almost independent of $\varepsilon$ while the non-smooth $\mathcal Q_0^\varepsilon$ exhibits a clear loss of regularity for decreasing $\varepsilon$.
}
\label{fig:1D_qoi}
\end{figure}

By using the explicit solution \eqref{eq:solex1}, we can write ${\mathcal Q}^\varepsilon_j(y)$ as a sum in the same way as
 in Section \ref{sec:4_1},
 $$
 {\mathcal Q}^\varepsilon_j(y)  = \Re \Bigl(\widetilde{Q}^\varepsilon_{11}(y)
 +\widetilde{Q}^\varepsilon_{12}(y)+\widetilde{Q}^\varepsilon_{21}(y)+\widetilde{Q}^\varepsilon_{22}(y)\Bigr),
 $$
where $\widetilde{Q}^\varepsilon_{11}(y)$ and $\widetilde{Q}^\varepsilon_{22}(y)$
are independent of $\varepsilon$ and
    \begin{align*}
    \widetilde Q_{12}^\varepsilon(y) &=\overline{ \widetilde Q_{21}^\varepsilon(y)}=
    \int g(x-s_1-y)\: g(x-s_2+y)
    e^{i\frac{\Phi_0(x-y)-\Phi_0(x+y)}{\varepsilon}}\psi_j(x)dx\\
    &=
    \int g(x-s_1-y)\: g(x-s_2+y)
    \exp\left(\frac{-4xyi}{\varepsilon}\right)\psi_j(x)dx.
    \end{align*}
    From this expression, one can derive that
for $\psi_0$, the derivatives behave as $\frac{\partial^k {\mathcal Q}_0^\varepsilon}{\partial y^k} \sim \frac{1}{\varepsilon^{k-1}}$. Hence, condition \eqref{eq:oscbound} is not fulfilled and we should observe oscillatory behavior.
On the other hand, the function $\psi_1 \in C_c^\infty$ suppresses the oscillations
when $y>0$
according to Theorem \ref{th:stat_phase} and we should expect
smooth higher derivatives of ${\mathcal Q}_1^\varepsilon$.

To demonstrate this numerically, we
plot in Figure \ref{fig:1D_qoi} the quantities of interest,
${\mathcal Q}_0^\varepsilon$ and
${\mathcal Q}_1^\varepsilon$, with their first and second derivatives.
Indeed, the derivatives of ${\mathcal Q}_1^\varepsilon$ are smooth,
while those of
${\mathcal Q}_0^\varepsilon$ oscillate; in particular,
the second derivative of ${\mathcal Q}_0^\varepsilon$ grows as $\frac{1}{\varepsilon}$
when $\varepsilon\to 0$, which means that \eqref{est} is violated. For accurate approximation of $\mathcal Q_0^\varepsilon$, we use 50 points per wavelength to discretize the spatial grid $x$ which makes the quadrature error negligible.

\subsubsection{Example 2: Variable speed of propagation in 1D}
 Problem \eqref{wave} is only analytically solvable in the d'Alembert form \eqref{eq:solex1} if the speed of propagation $c$ is constant. In more general cases, the Gaussian beam approximation has to be implemented.
Analogously to the previous example, we choose the initial data as
$$
  A_0(x,\mathbf y) = g(x-s_1(\mathbf y)) + g(x-s_2(\mathbf y)), \quad g(x) = e^{-10x^2},\qquad \Phi_0(x) = |x|,
$$
which represents two high-frequency pulses moving towards each other. We set the Gaussian beam initial data as in \eqref{xpMinit}. The initial position of the pulses and the non-constant speed of propagation depend on two stochastic variables, $\mathbf y = (y_1,y_2)$. More precisely,
\[
    s_1(\mathbf y) = - s_2(\mathbf y) = y_1, \qquad c(x,{\bf y}) = 1 + \frac{1}{2}\left(e^{-(x-1)^2} + y_2 e^{(x + 1)^2}\right).
\]

\begin{figure}
\centering
    \begin{minipage}{0.43\textwidth}
            \subfigure[$t=0.25$]{
\includegraphics[width=\textwidth]{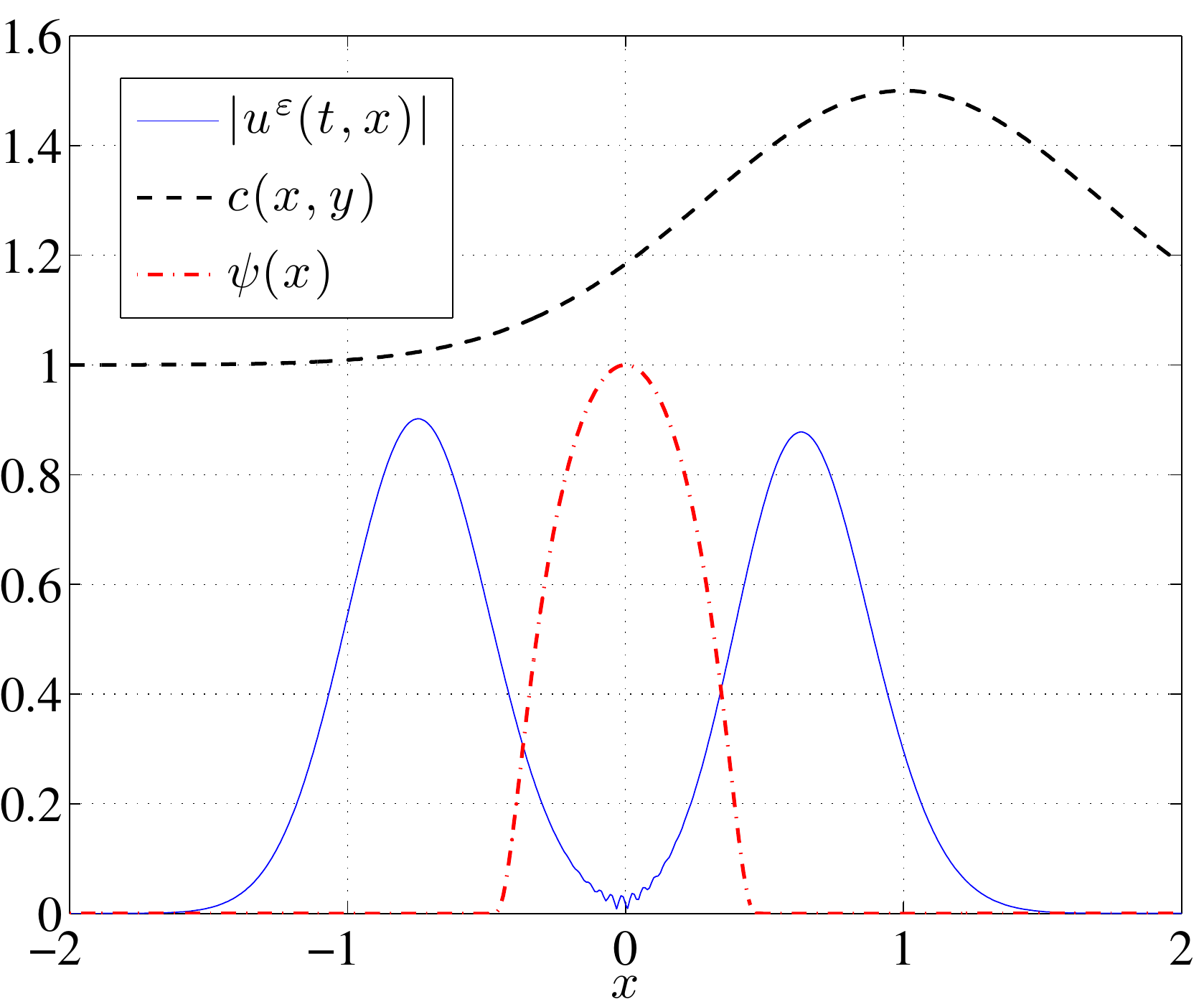}}
    \end{minipage}
    \begin{minipage}{0.43\textwidth}
            \subfigure[$t=0.5$]{
\includegraphics[width=\textwidth]{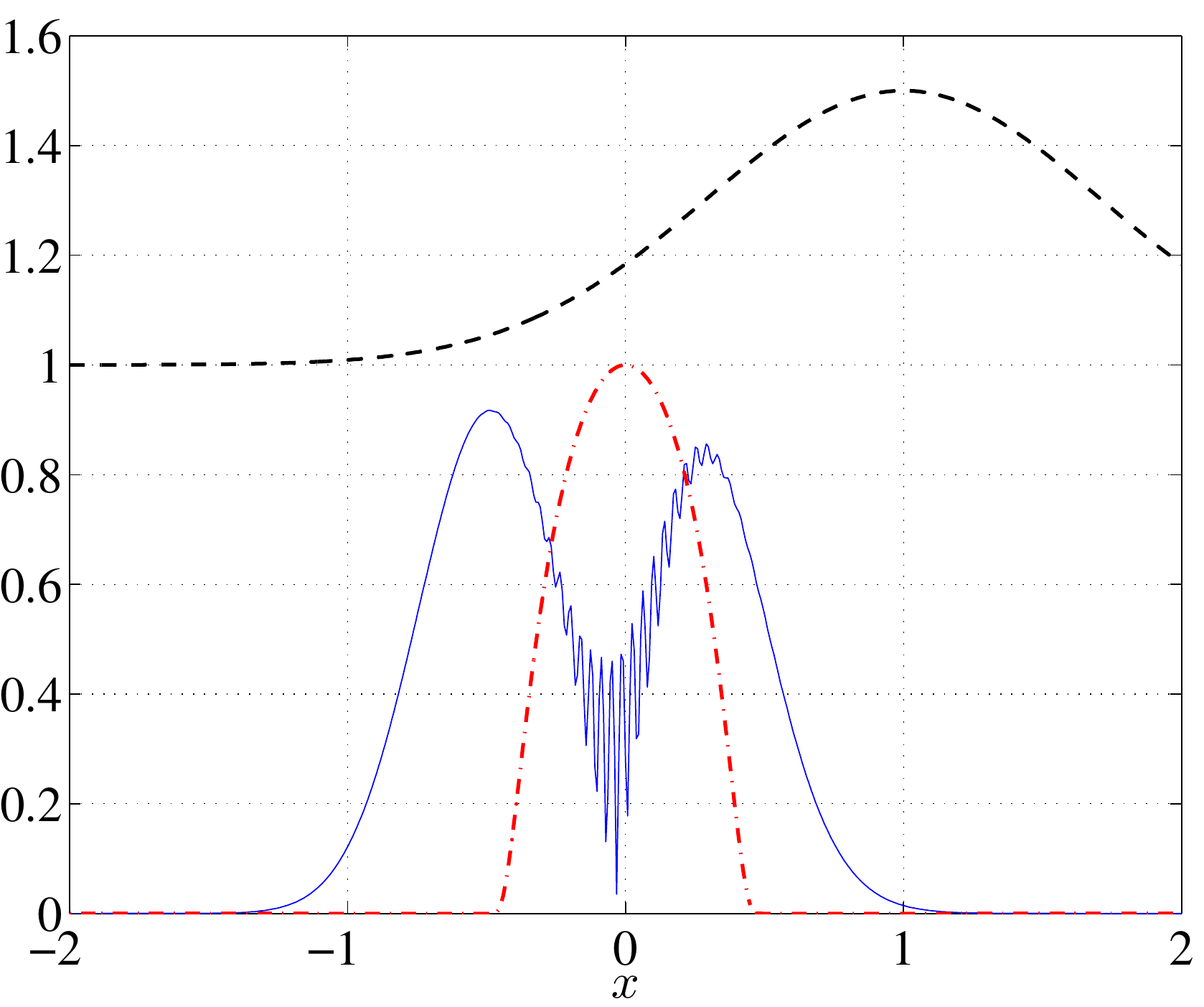}}
    \end{minipage}
    \begin{minipage}{0.43\textwidth}
            \subfigure[$t=0.75$]{
\includegraphics[width=\textwidth]{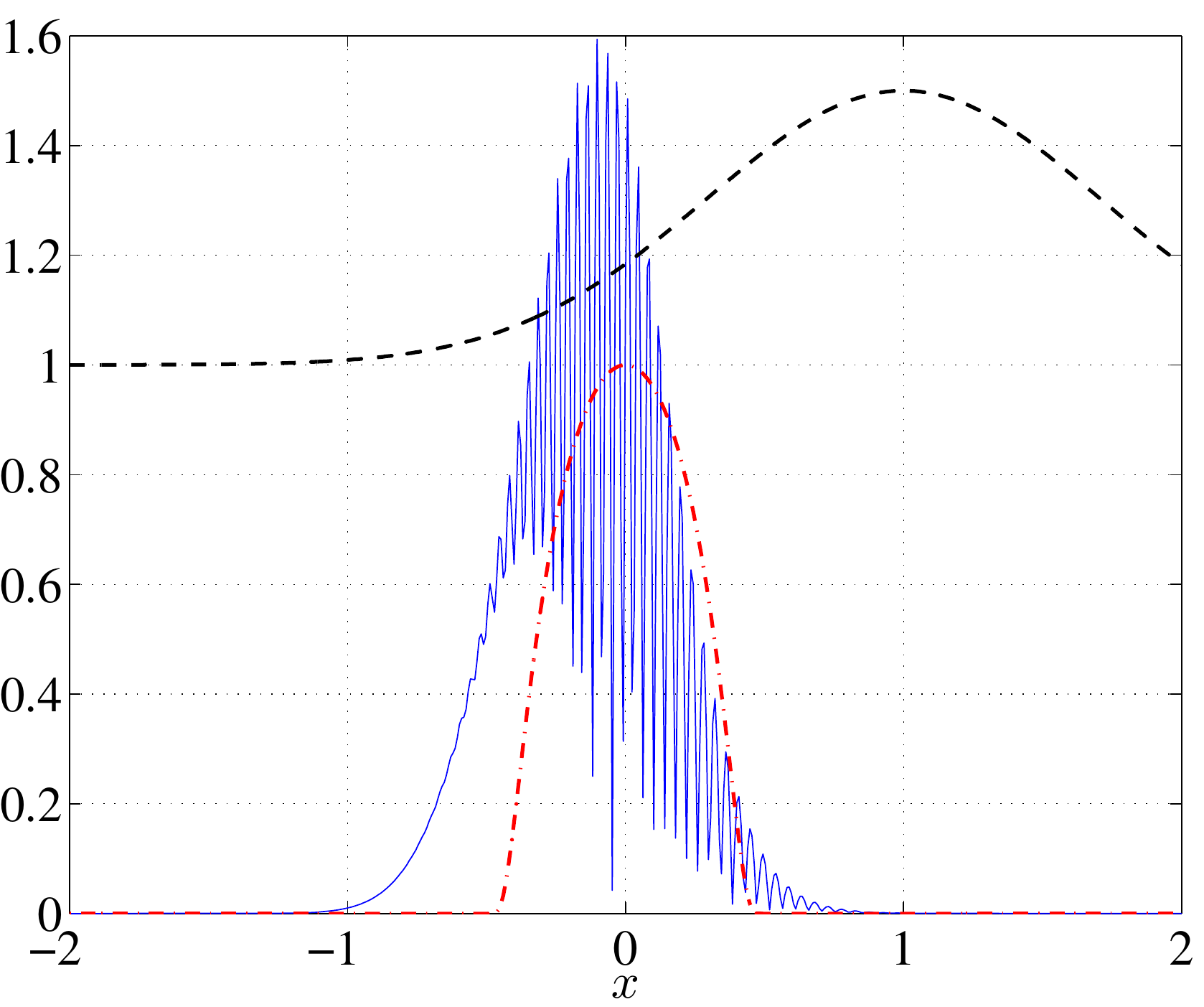}}
    \end{minipage}
    \begin{minipage}{0.43\textwidth}
            \subfigure[$t=1$]{
\includegraphics[width=\textwidth]{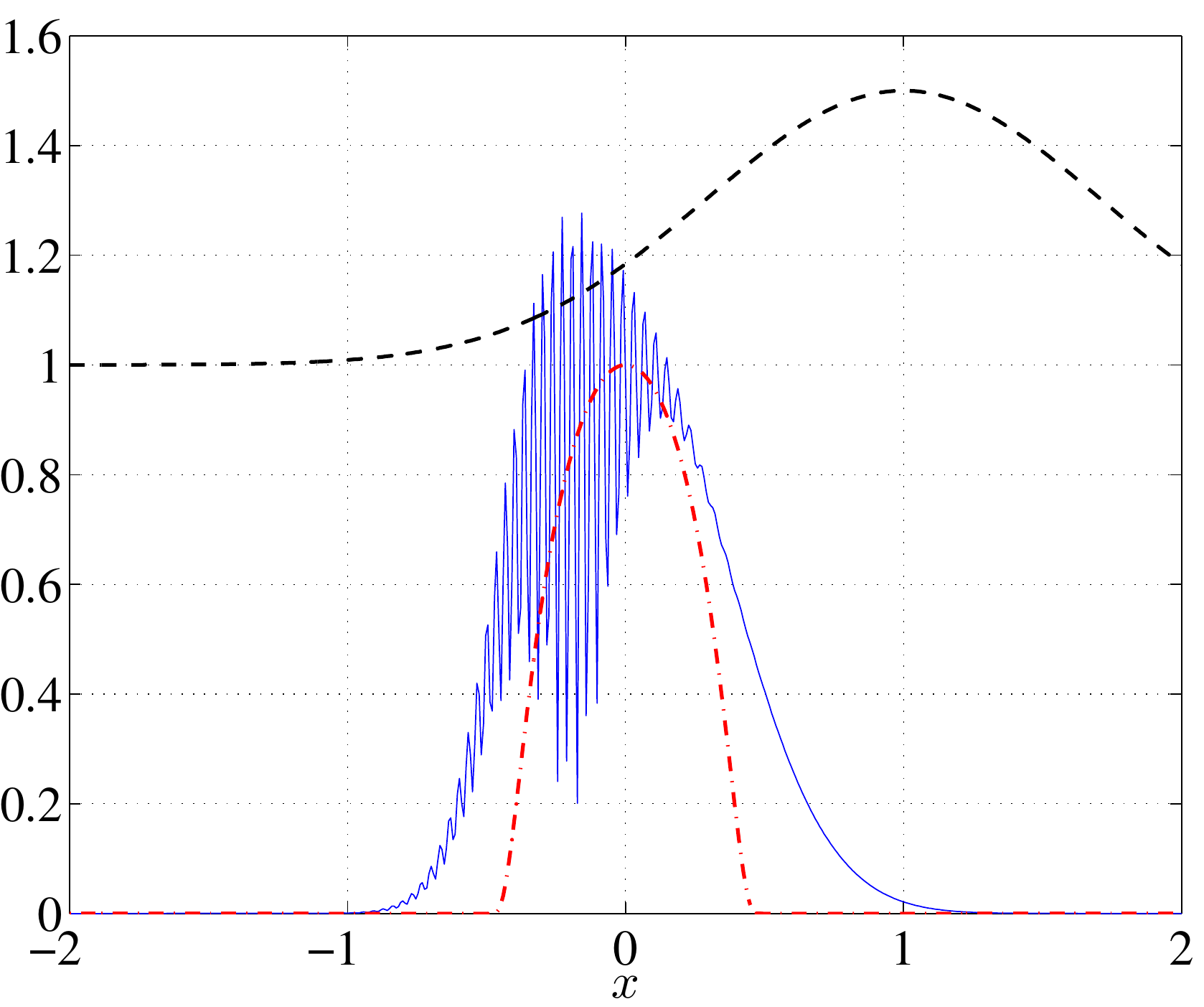}}
    \end{minipage}
\caption{Example 2. Absolute value of solution for $t = 0.25, \: 0.5, \: 0,75, \: 1$ when $\varepsilon = 1/80$ and ${\bold y}=(1,0)$. The QoI test function, $\psi$, and speed, $c$, are overlaid.}
\label{fig:ex2_1}
\end{figure}

\begin{figure}
\centering
\begin{minipage}{.43\textwidth}
            \subfigure[${\bold y}=(1,0.25)$]{
\includegraphics[width=\textwidth]{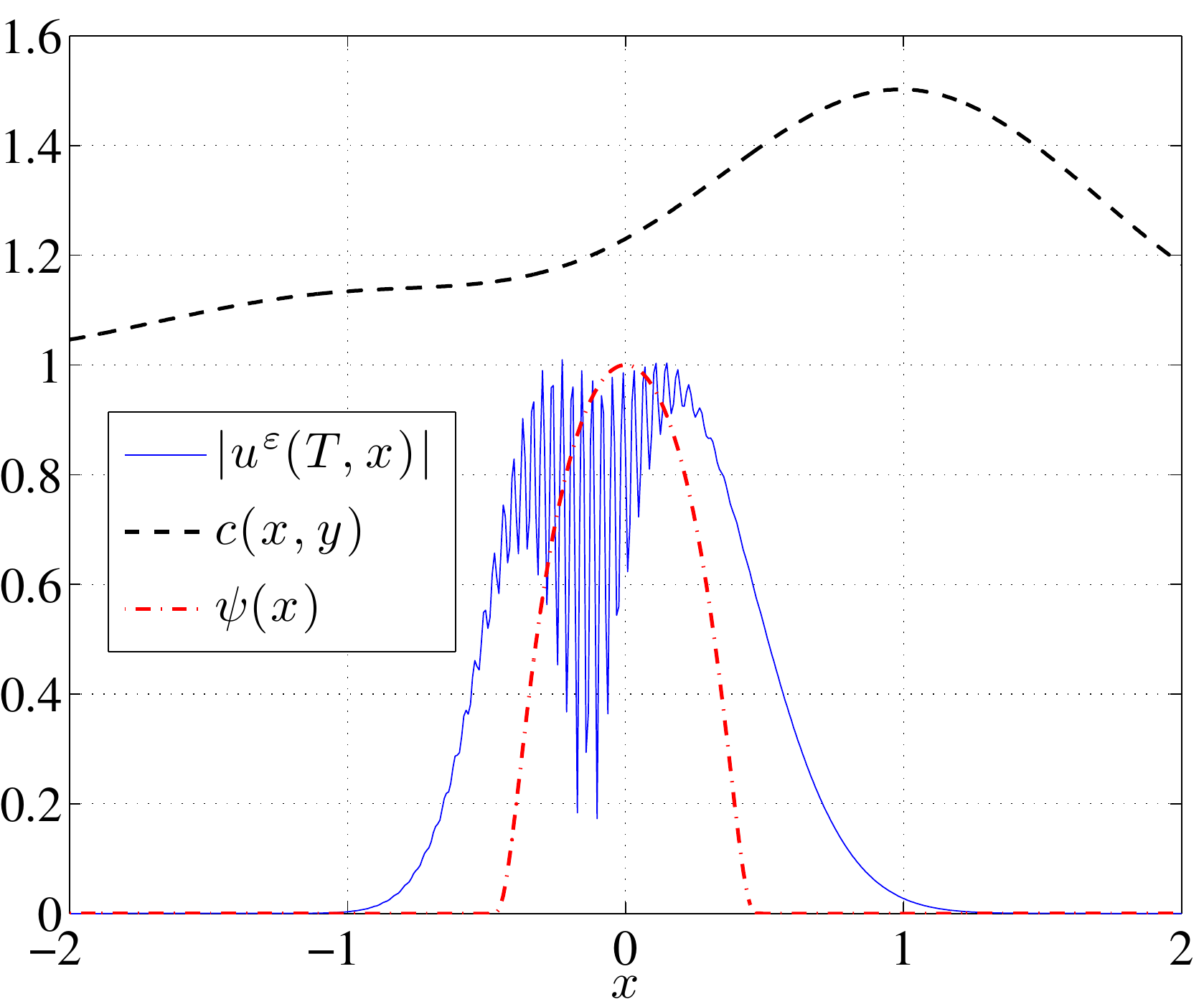}}
\end{minipage}
\begin{minipage}{.43\textwidth}
            \subfigure[${\bold y}=(1,0.5)$]{
\includegraphics[width=\textwidth]{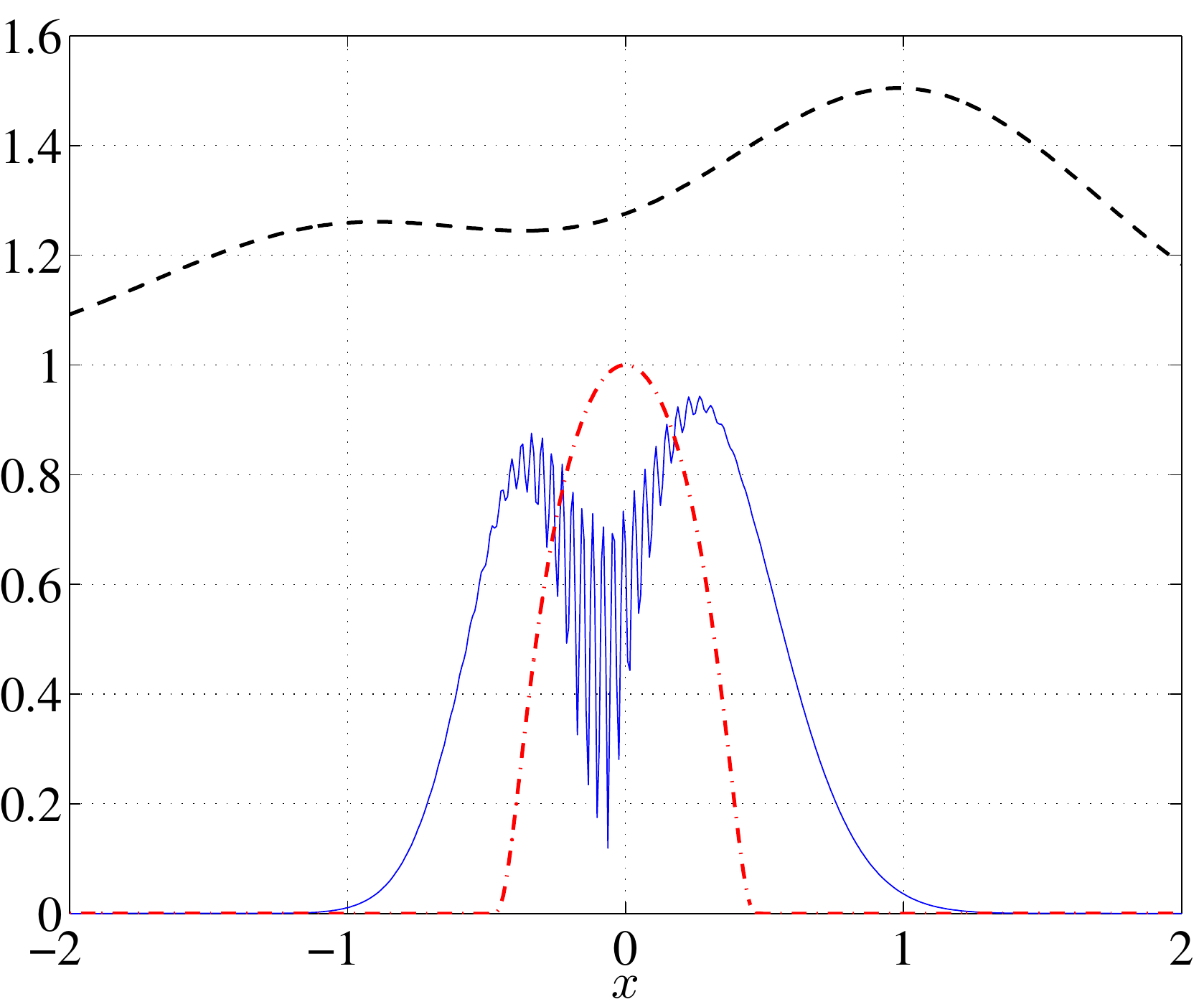}}
\end{minipage}
\begin{minipage}{.43\textwidth}
            \subfigure[${\bold y}=(1,0.75)$]{
\includegraphics[width=\textwidth]{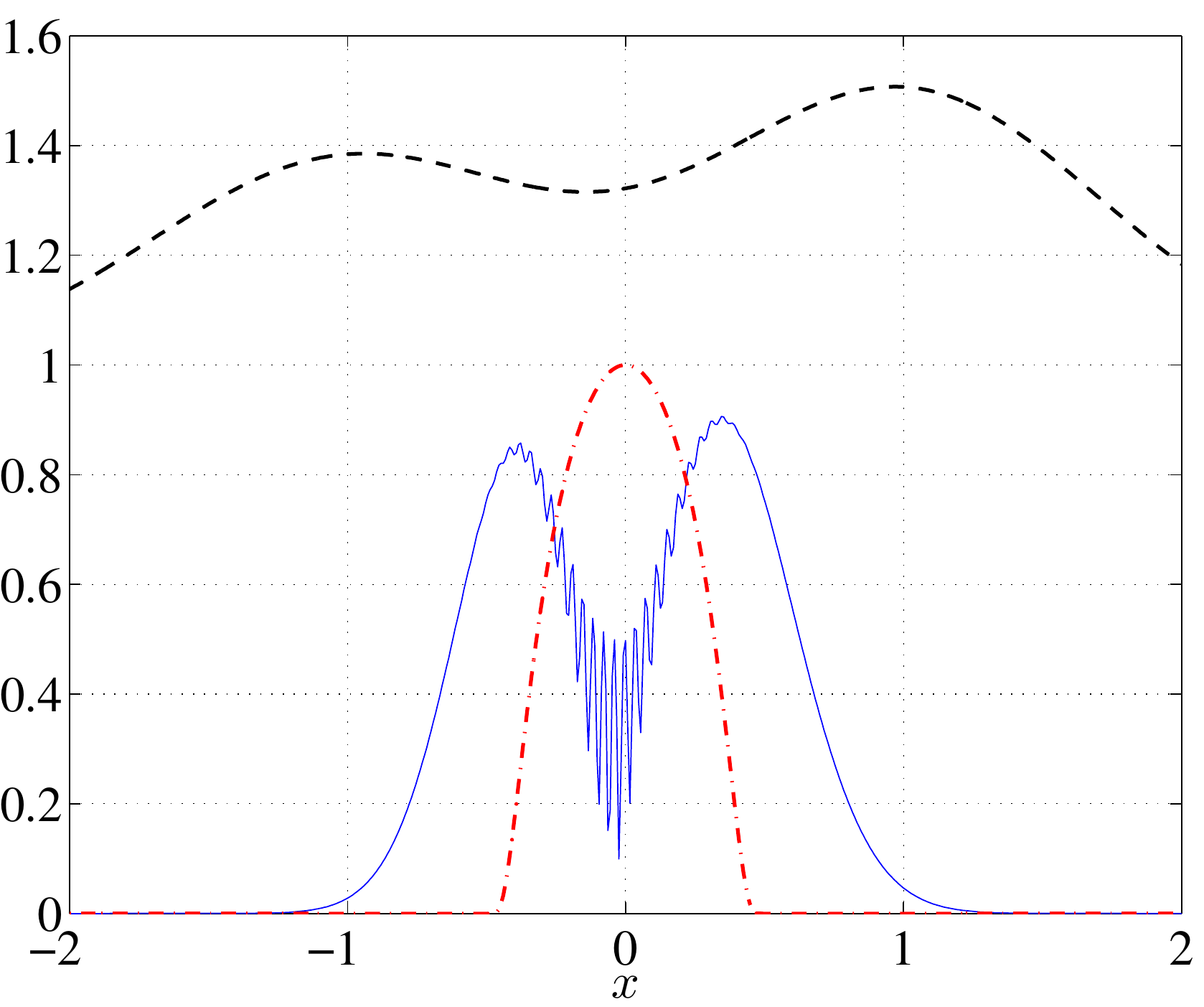}}
\end{minipage}
\begin{minipage}{.43\textwidth}
            \subfigure[${\bold y}=(1,1)$]{
\includegraphics[width=\textwidth]{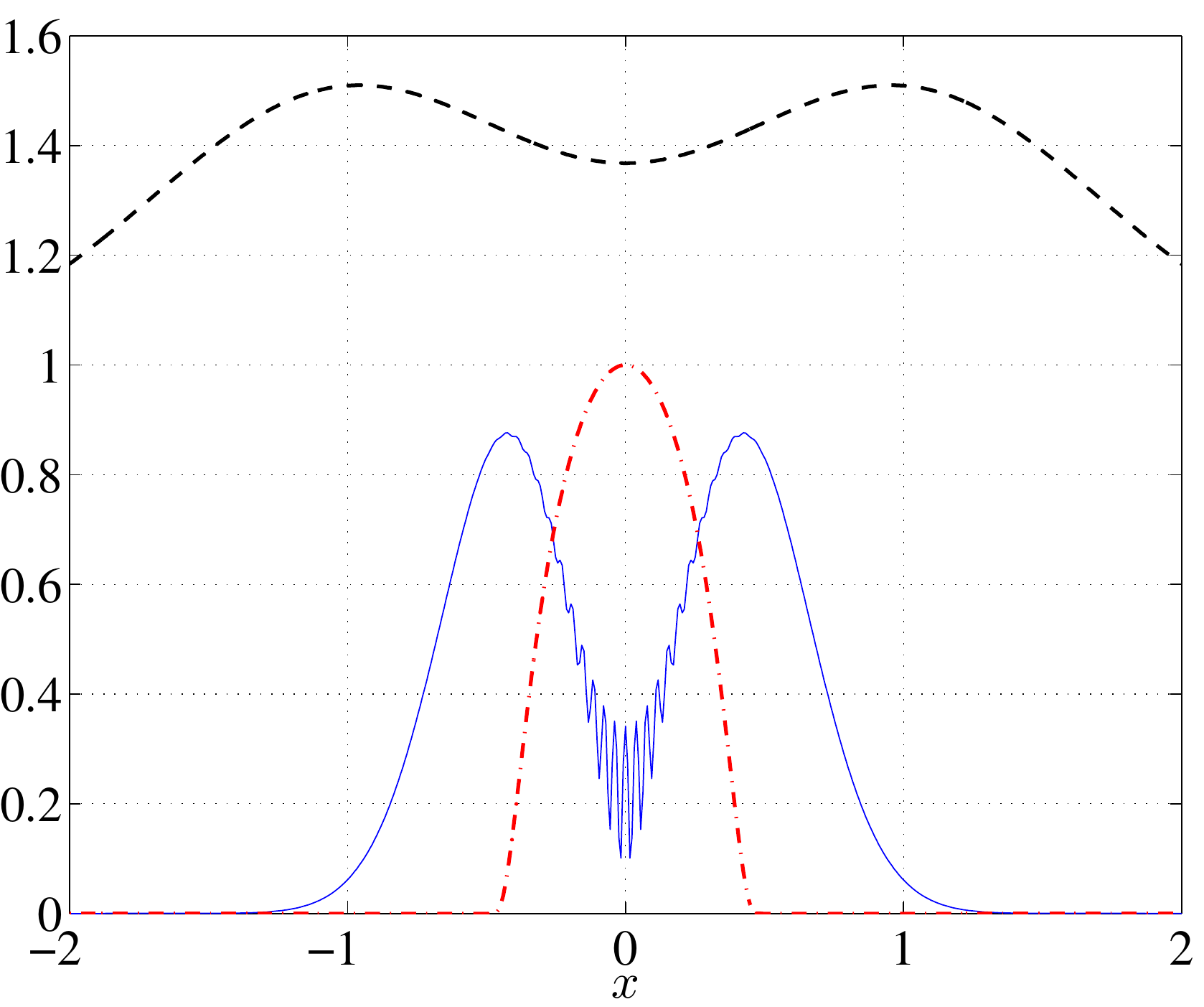}}
\end{minipage}
\caption{Example 2:
Four realizations of the absolute value of the solution at time $T=1$ with $\varepsilon = 1/80$. The
QoI test function, $\psi$,
and speed, $c$, are overlaid.
}
\label{fig:ex2_2}
\end{figure}

The absolute value of the solution with $\varepsilon = 1/80$ for $\mathbf y = (1,0)$ at
different times
is plotted in Figure \ref{fig:ex2_1}.
Figure \ref{fig:ex2_2} shows four different realizations of the absolute value of the solution with $\varepsilon = 1/80$ at time $T=1$. In all realizations, we keep $y_1=1$ fixed and vary $y_2$. This corresponds to different realizations of the wave speed.
The QoI test function, $\psi(x) = \widetilde{\psi}(2x)$, and the wave speed, $c(x,{\bf y})$, corresponding to different realizations are also shown in the figure.

\begin{figure}
            \centering
\begin{minipage}{.32\textwidth}
            \subfigure[${\mathcal Q}^\varepsilon({\bold y}(r))$]{
                    \includegraphics[width=\textwidth]{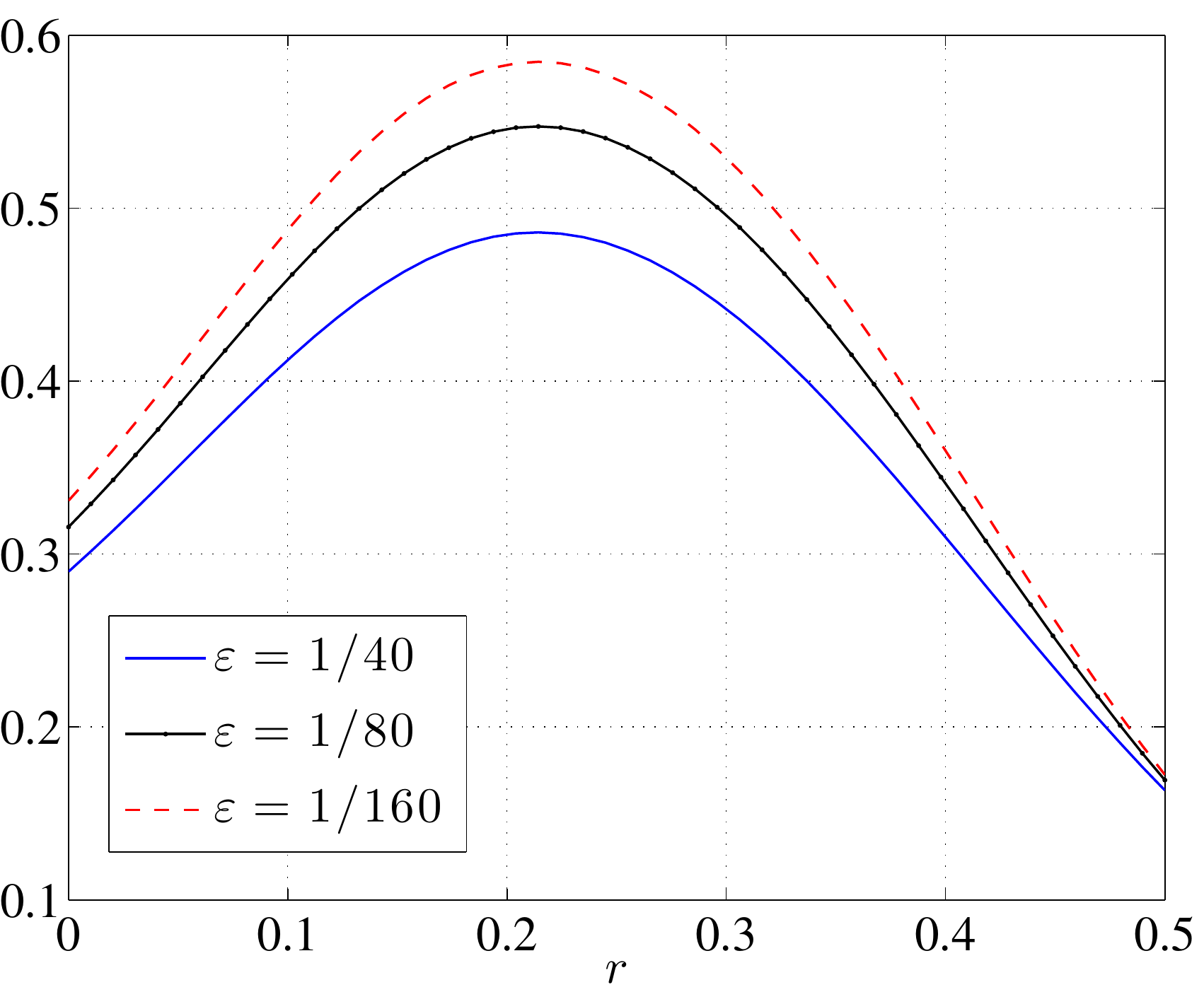}}
\end{minipage}
\begin{minipage}{.32\textwidth}
            \subfigure[$\frac{d}{dr}{\mathcal Q}^\varepsilon({\bold y}(r))$]{
                    \includegraphics[width=\textwidth]{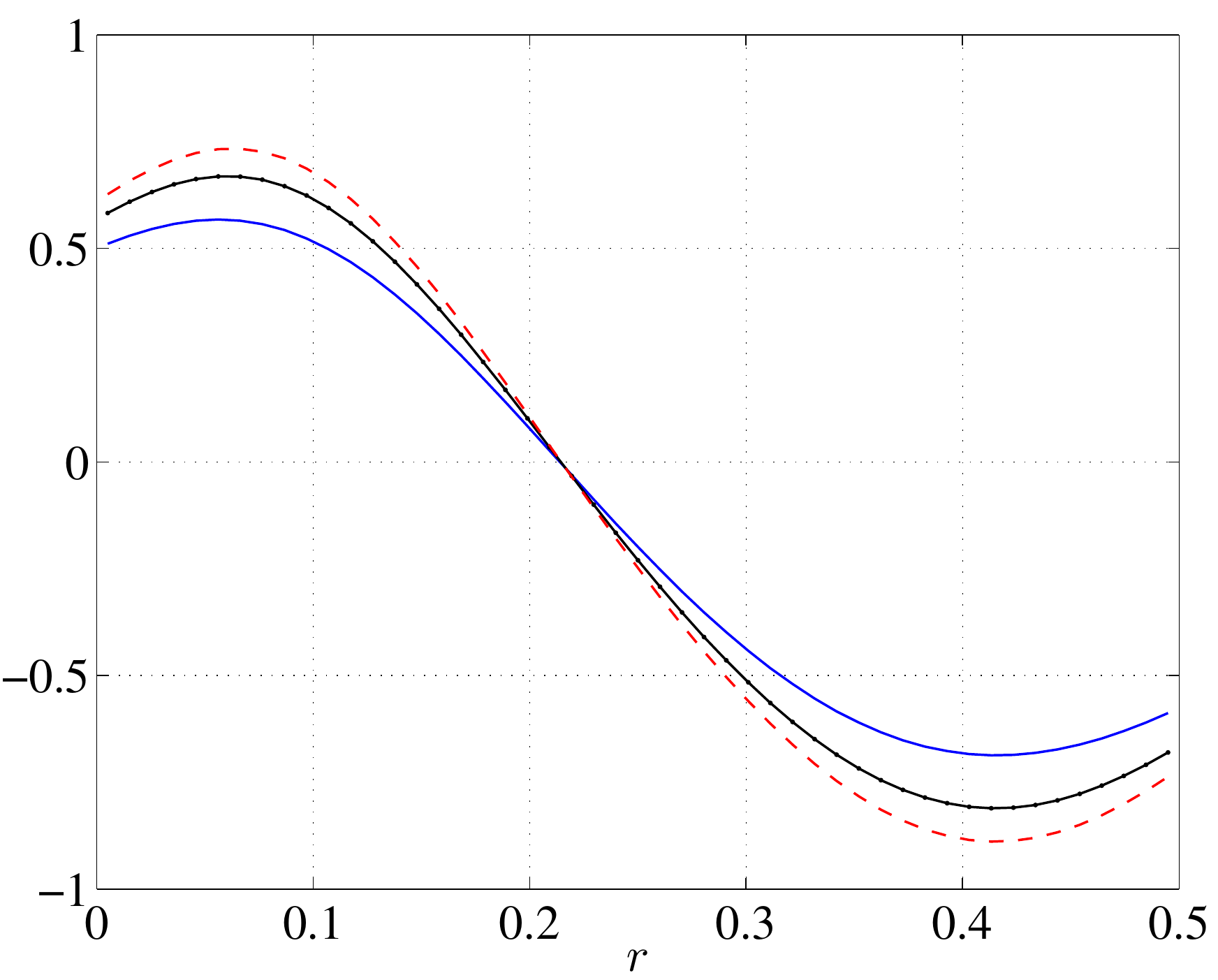}}
\end{minipage}
\begin{minipage}{.32\textwidth}
            \subfigure[$\frac{d^2}{dr^2}{\mathcal Q}^\varepsilon({\bold y}(r))$]{
                    \includegraphics[width=\textwidth]{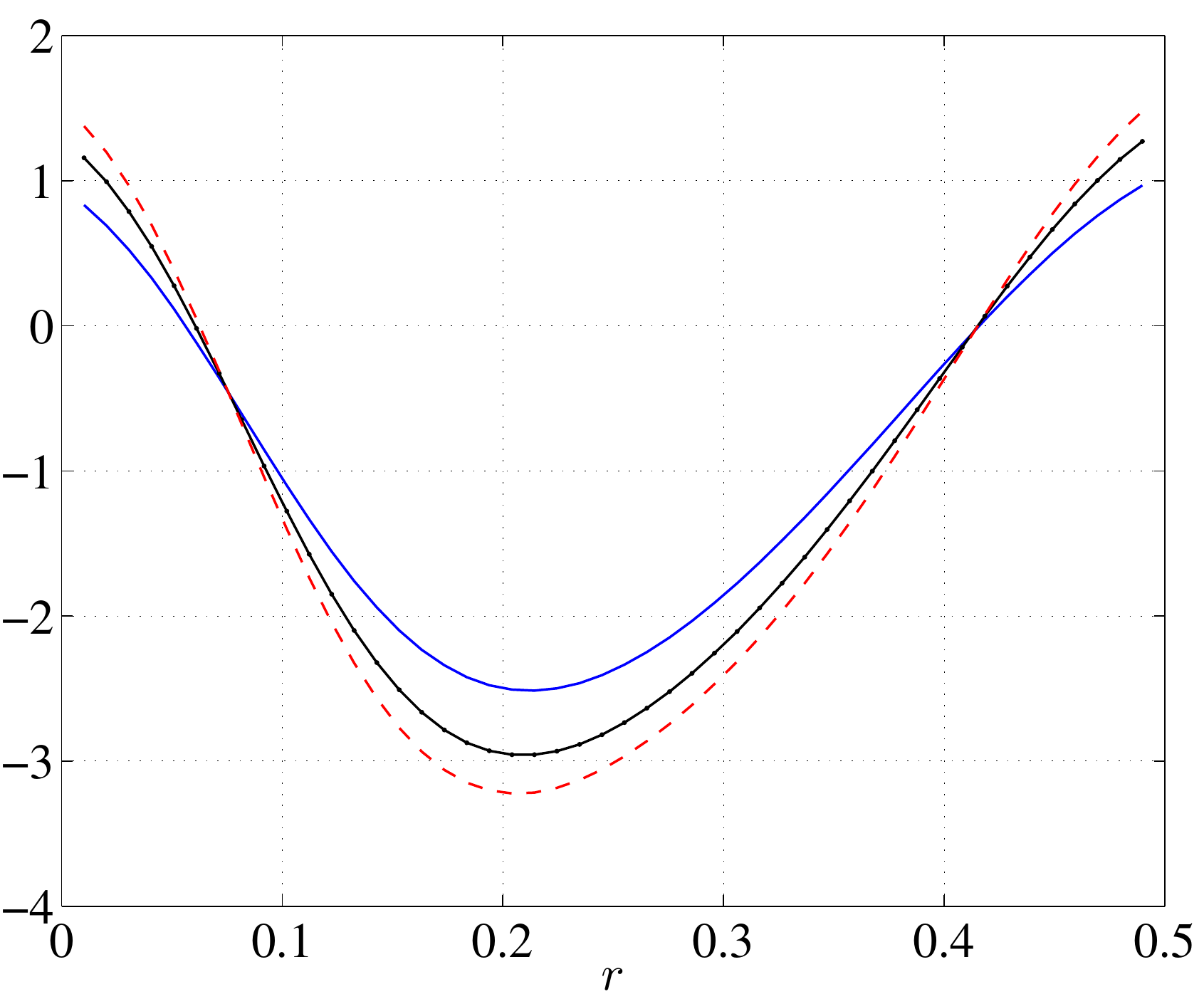}}
\end{minipage}
            \caption{Example 2:
            QoI
${\mathcal Q}^\varepsilon$ and its first and second derivatives
along the line ${\bold y}(r) = (1+r,1-2r)$ for $r\in[0,0.5]$
and different wavelengths $\varepsilon$.}
\label{fig:1D_qoi_nonc}
\end{figure}

In Figure \ref{fig:1D_qoi_nonc}, we plot the quantity of interest \eqref{Q} and its first and second derivatives along the line ${\bold y}(r) = (1 + r,1-2r)$.
We observe that ${\mathcal Q}^\varepsilon$ and its derivatives are smooth
and do not oscillate with $\varepsilon$. Hence, the estimate \eqref{est} is fulfilled.

\subsubsection{Example 3: Caustics in 2D}

In our final example, we consider a two-dimensional case with
two initial wave pulses,
$$
  A_0({\bold x},{\bold y}) = g({\bold x}-{\bold s}_1({\bold y})) + g({\bold x}-{\bold s}_2({\bold y})), \qquad
  g({\bold x})=e^{-5|{\bold x}|^2}.
$$
The deterministic initial phase
$$
\Phi_0(\mathbf x) = |x_1| + x_2^2,\qquad {\bold x} = (x_1,x_2),
$$
is chosen such that a cusp caustic develops at $t=0.5$,
and two fold caustics form at $t>0.5$.

Both the initial location of the pulses and the constant speed of
propagation are random, depending on two stochastic variables,
${\bold y}=(y_1,y_2)$, as
$$
{\bold s}_1({\bold y})=-{\bold s}_2({\bold y})=(y_1,0),\qquad
  c({\bold y}) = y_2.
$$

We consider the QoI \eqref{Q} with $\psi(\mathbf x)=\widetilde{\psi}(2\mathbf x)$.

\begin{figure}
\centering
\begin{minipage}{.32\textwidth}
\subfigure[$t=0$]{
\includegraphics[width=\textwidth]{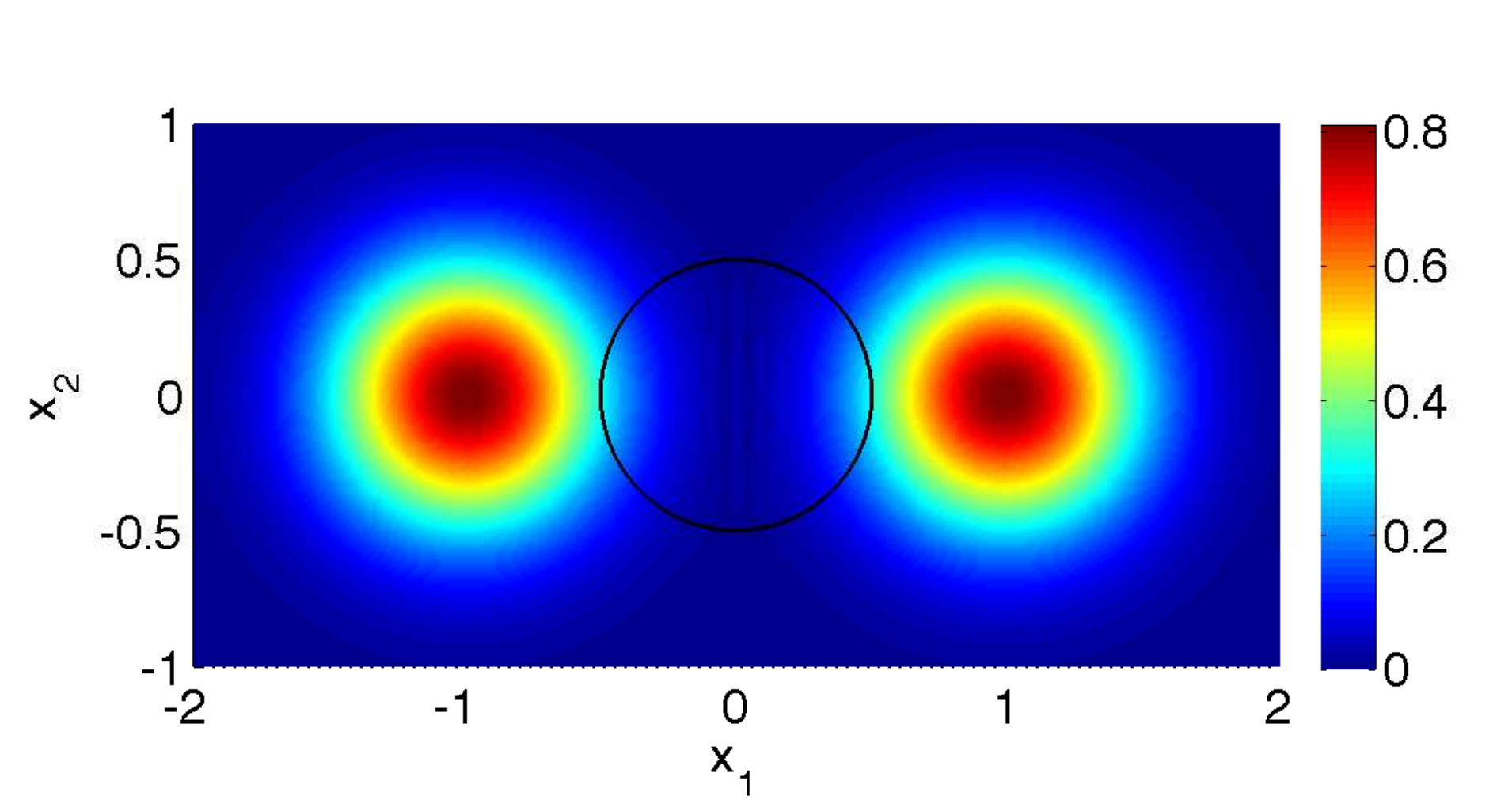}}
\end{minipage}
\begin{minipage}{.32\textwidth}
\subfigure[$t=0.2$]{
\includegraphics[width=\textwidth]{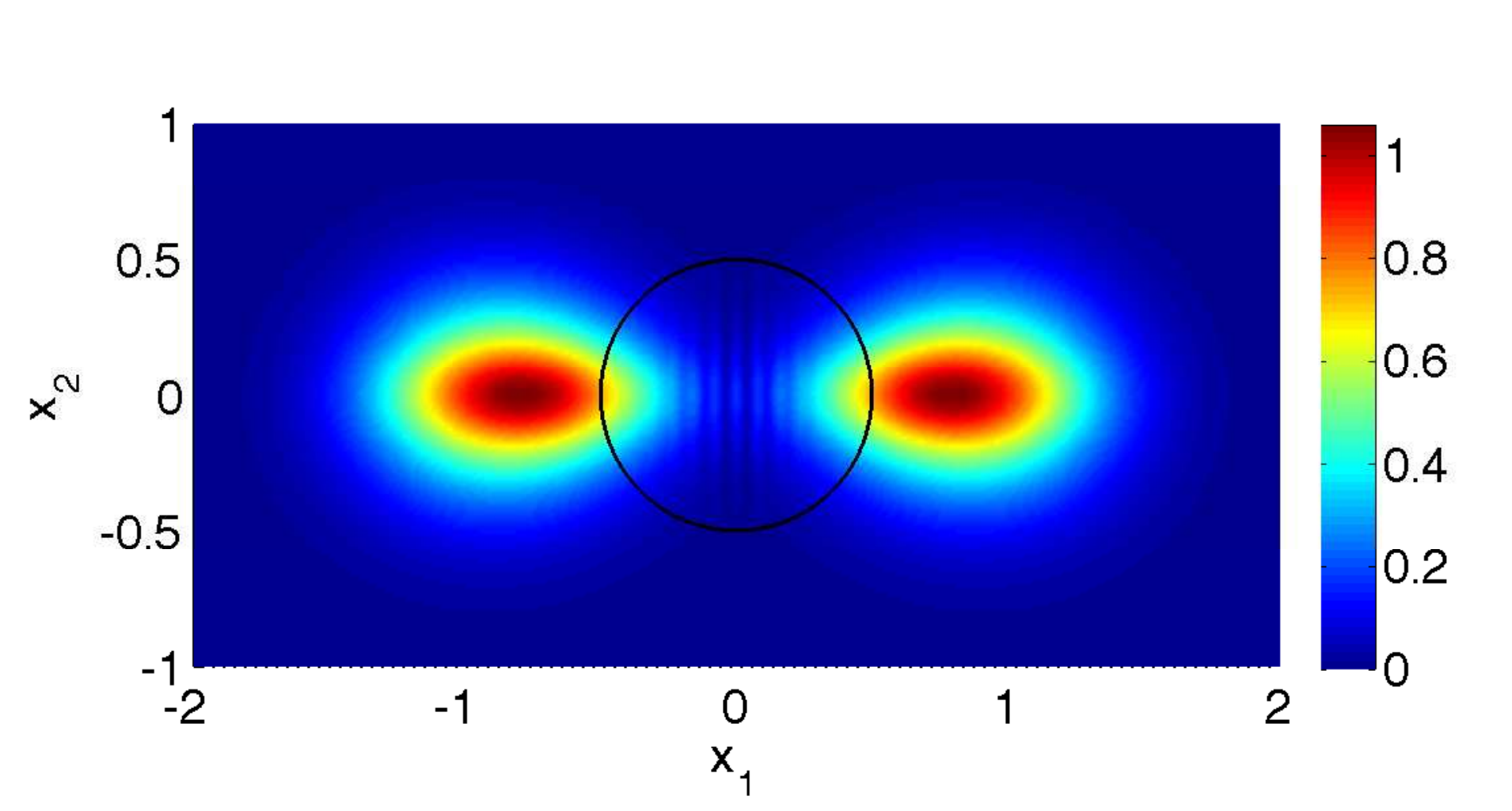}}
\end{minipage}
\begin{minipage}{.32\textwidth}
\subfigure[$t=0.4$]{
\includegraphics[width=\textwidth]{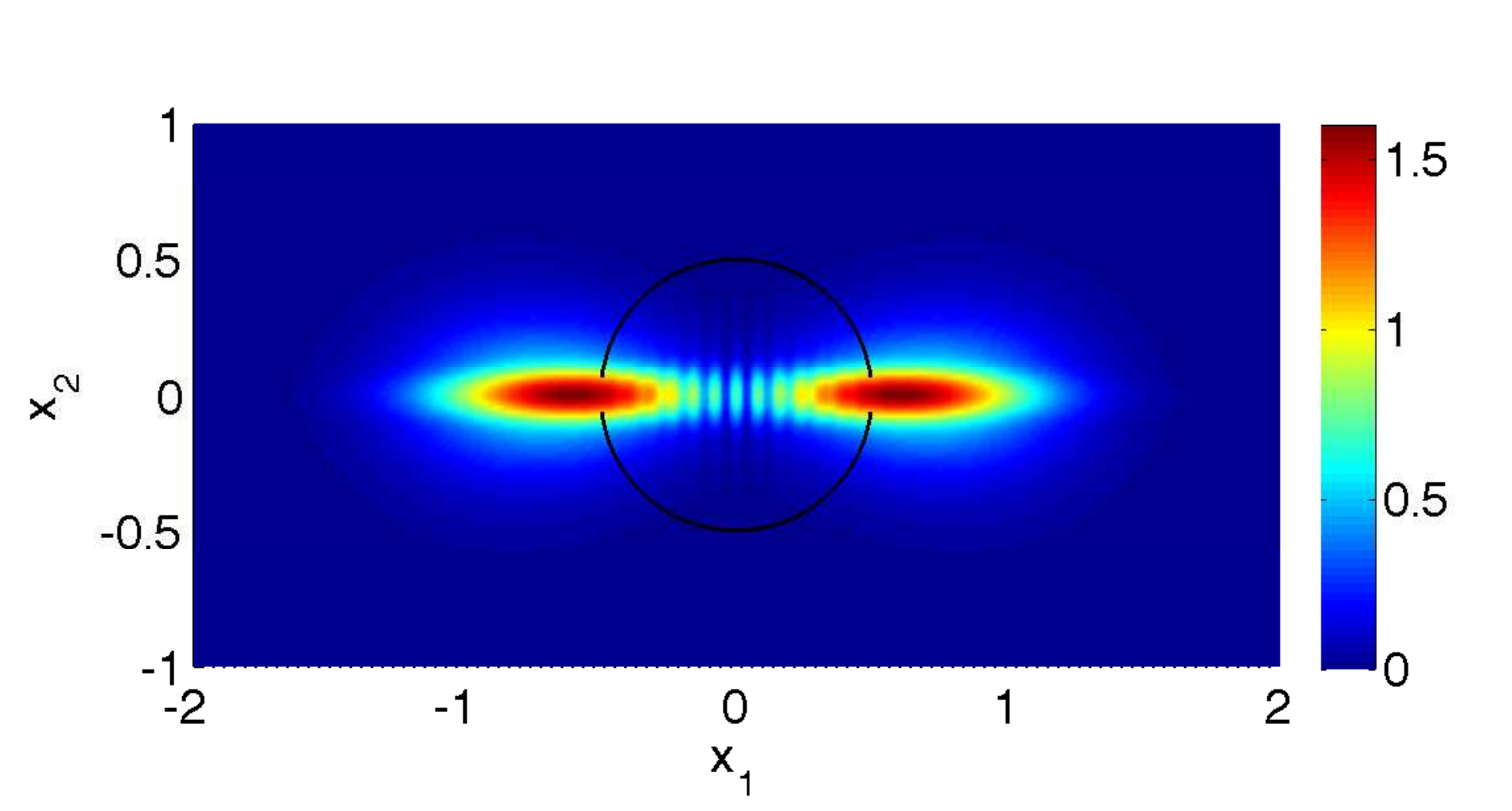}}
\end{minipage}
\begin{minipage}{.32\textwidth}
\subfigure[$t=0.6$]{
\includegraphics[width=\textwidth]{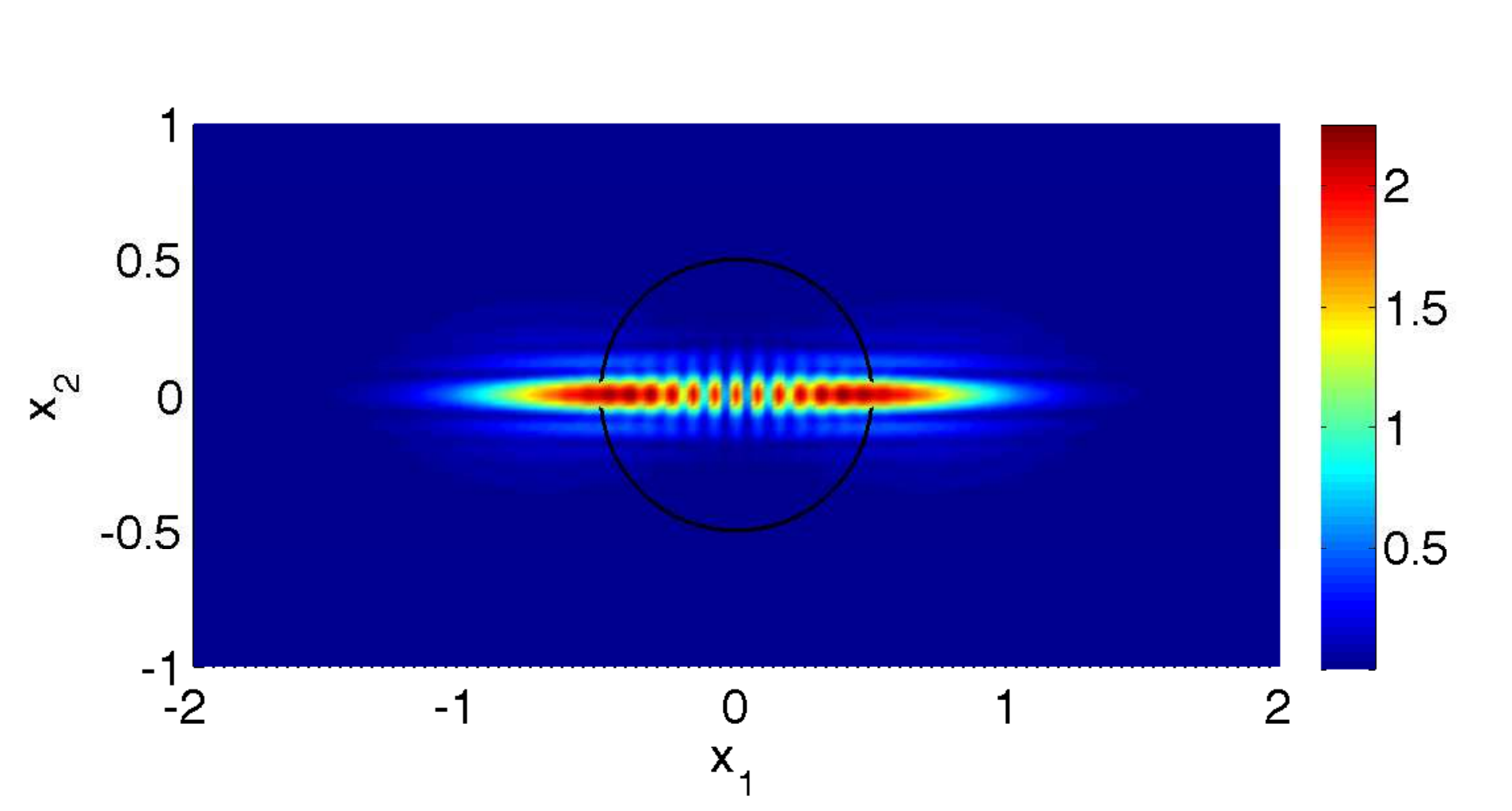}}
\end{minipage}
\begin{minipage}{.32\textwidth}
\subfigure[$t=0.8$]{
\includegraphics[width=\textwidth]{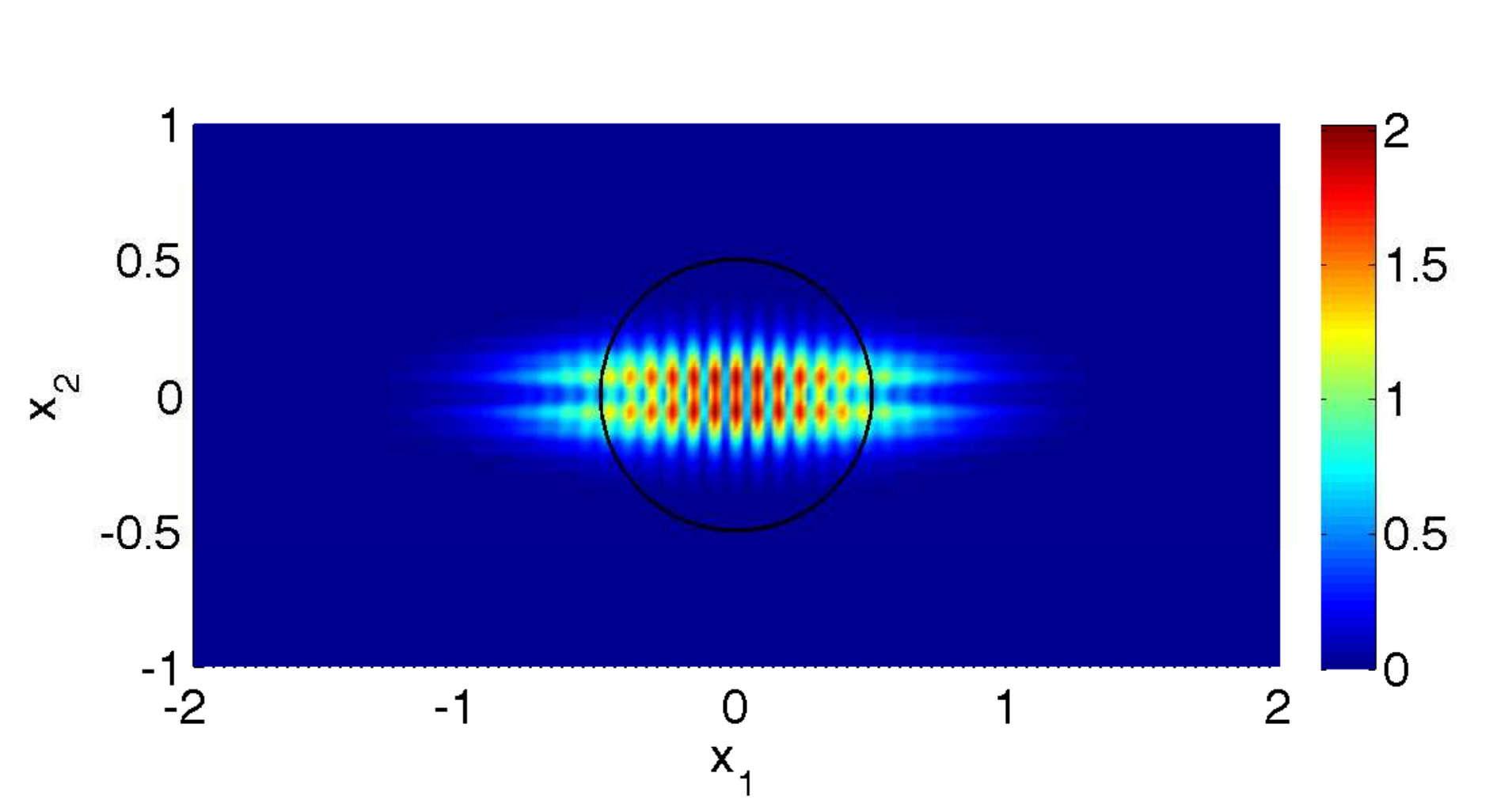}}
\end{minipage}
\begin{minipage}{.32\textwidth}
\subfigure[$t=1.0$]{
\includegraphics[width=\textwidth]{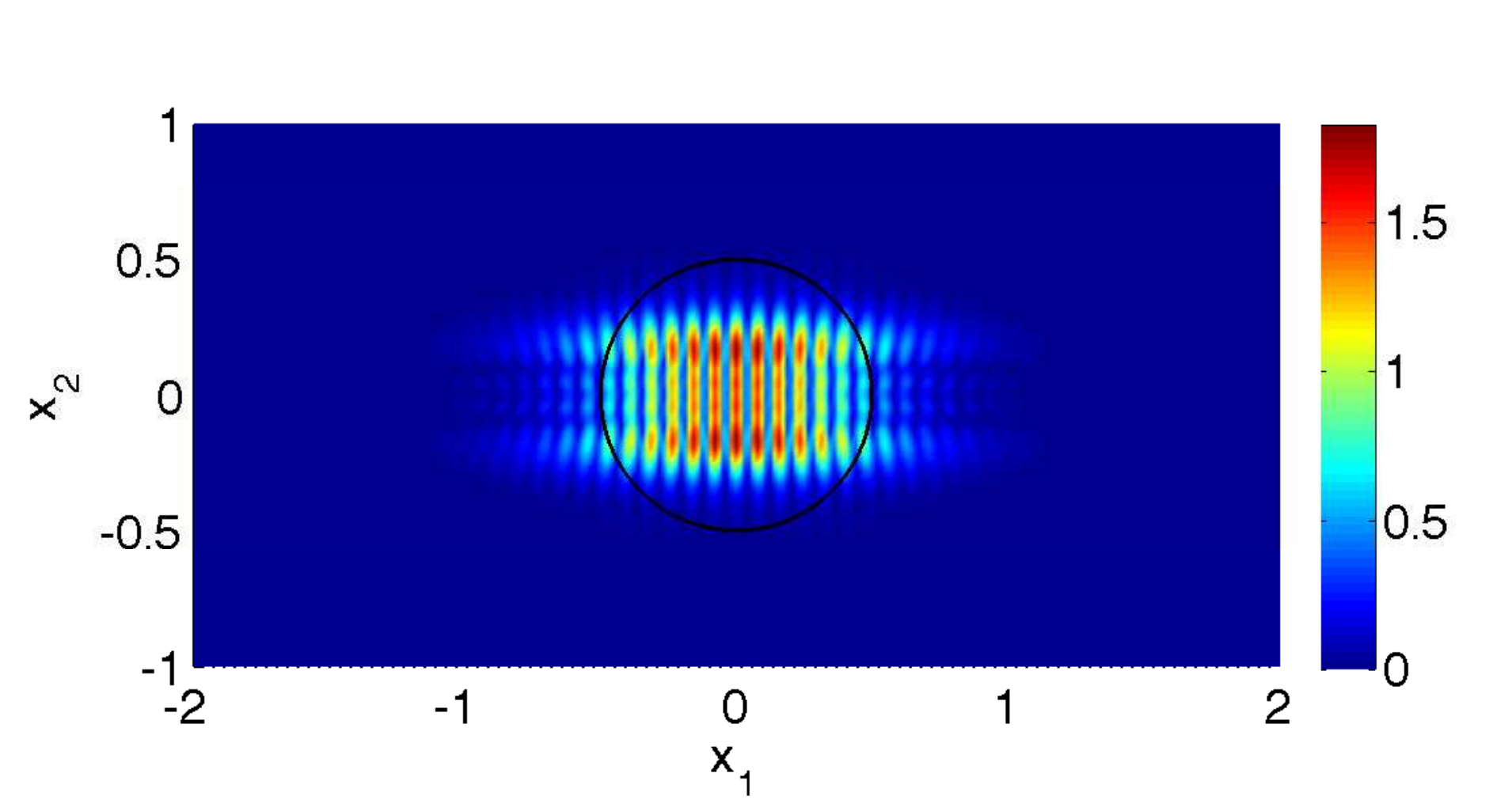}}
\end{minipage}
\caption{Example 3:
Absolute value of solution for various times, $t$, when $\mathbf y = (1,1)$.
Caustics appear for $t\geq 0.5$.
Circle indicates the support of the QoI test function.}
\label{fig:2D_sol_caust2}
\end{figure}

\begin{figure}
            \centering
\begin{minipage}{.32\textwidth}
            \subfigure[${\mathcal Q}^\varepsilon({\bold y}(r))$]{
                    \includegraphics[width=\textwidth]{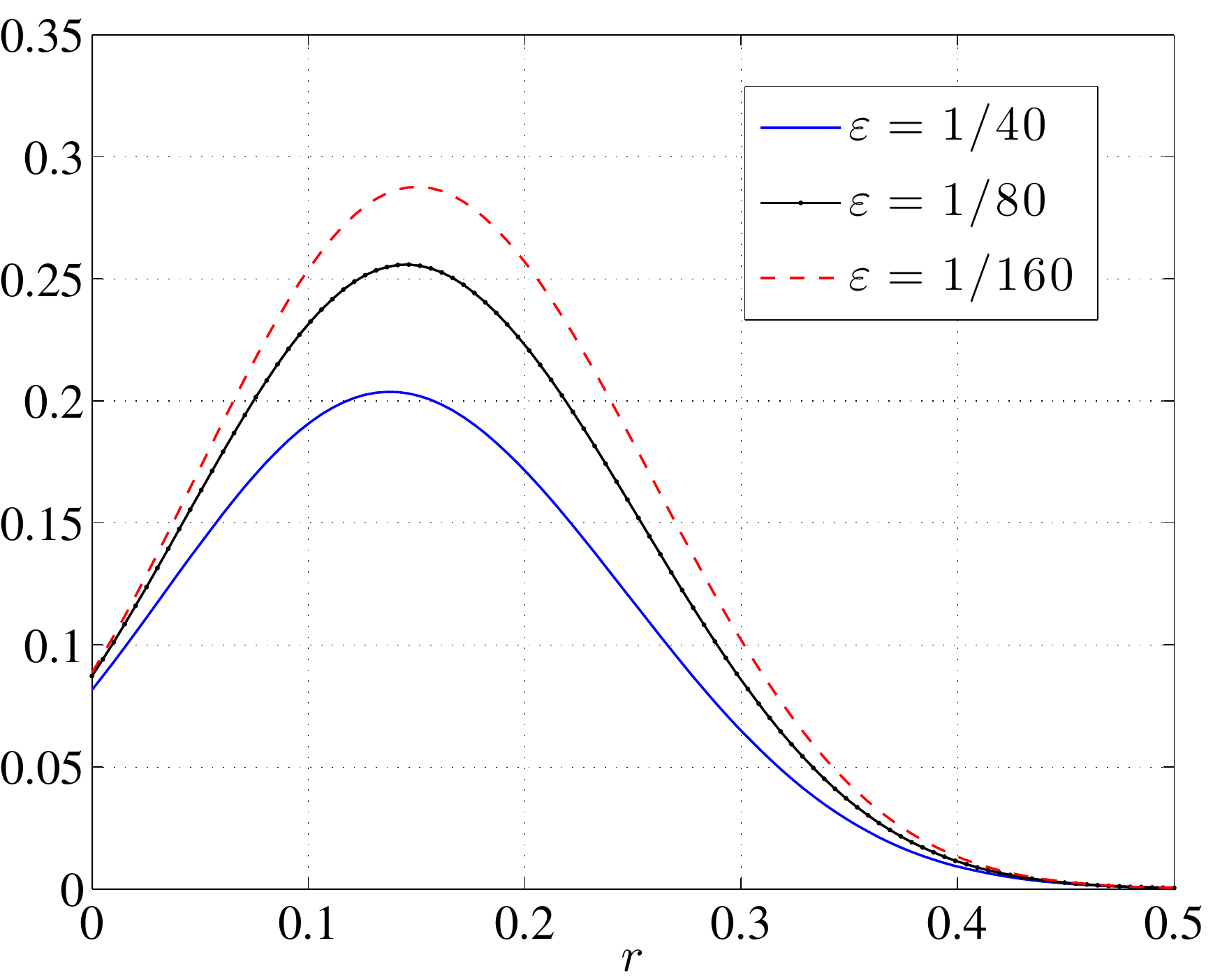}}
\end{minipage}
\begin{minipage}{.32\textwidth}
            \subfigure[$\frac{d}{dr}{\mathcal Q}^\varepsilon({\bold y}(r))$]{
                    \includegraphics[width=\textwidth]{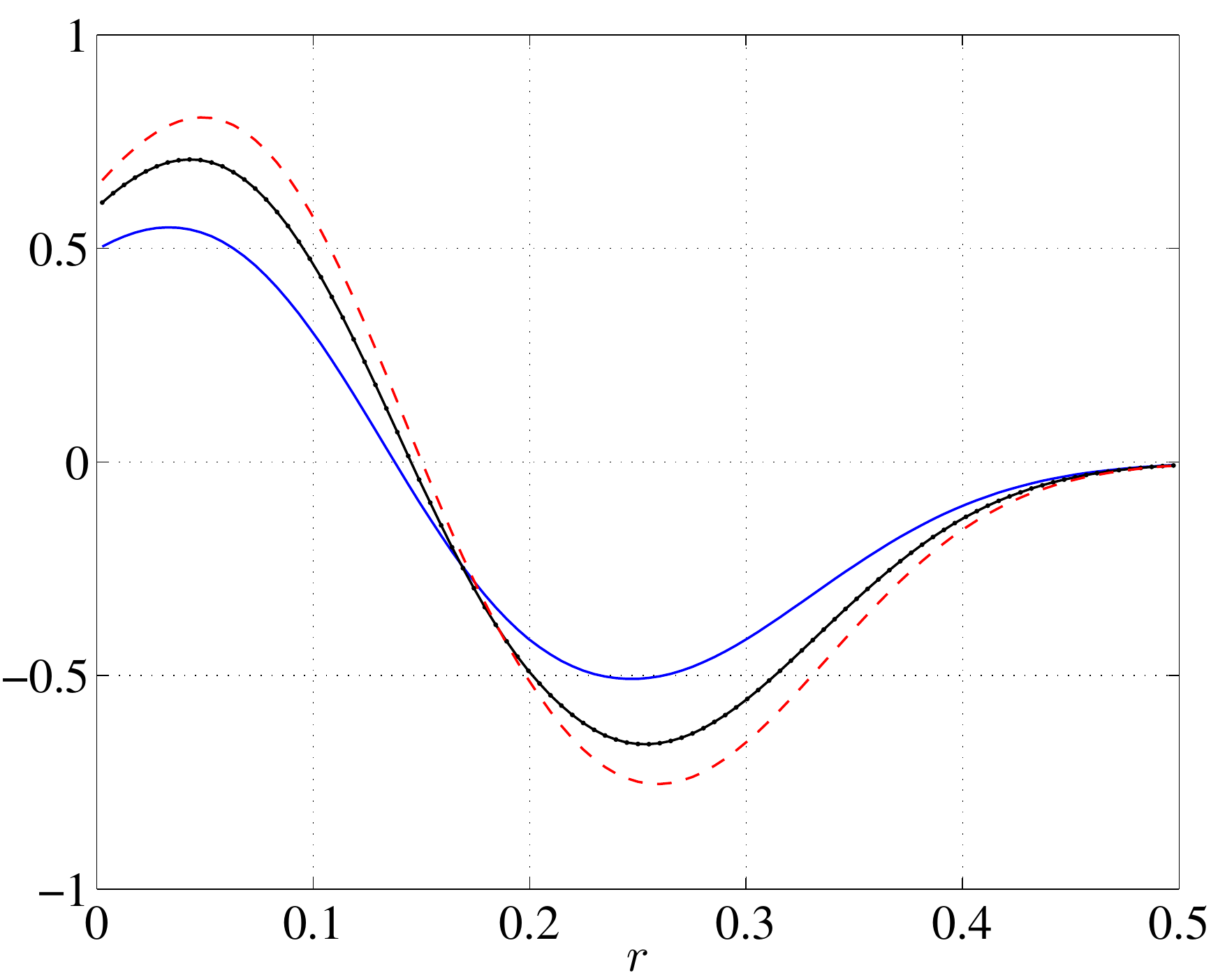}}
\end{minipage}
\begin{minipage}{.32\textwidth}
            \subfigure[$\frac{d^2}{dr^2}{\mathcal Q}^\varepsilon({\bold y}(r))$]{
                    \includegraphics[width=\textwidth]{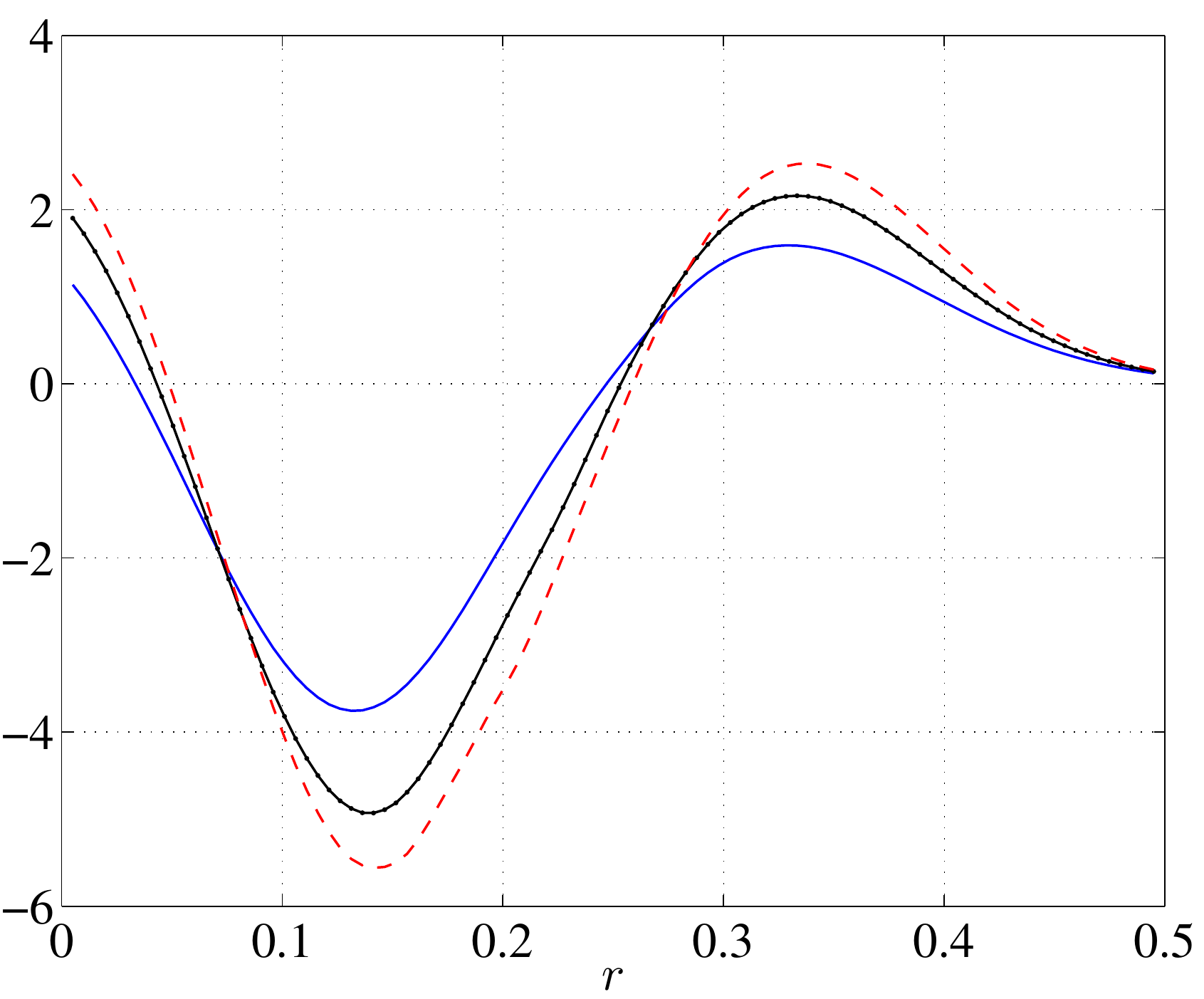}}
\end{minipage}
            \caption{Example 3:
            QoI ${\mathcal Q}^\varepsilon$, and its first and second derivatives
along the line ${\bold y}(r) = (1+r,1+2r)$ where $r \in [0,0.5]$,
for different wavelengths, $\varepsilon$.}
            \label{fig:2D_qoi_caust2}
\end{figure}

In Figure \ref{fig:2D_sol_caust2} the absolute value of the solution
at various times is shown. In this simulation,
the wavelength is $\varepsilon = 1/ 40$, the pulse centers are
${\bold s}_1=-{\bold s}_2=(1,0)$,
and the speed of propagation
$c \equiv 1$. The bold line represents the support of the QoI test function, $\psi$.
Figure \ref{fig:2D_qoi_caust2} shows the QoI
and its derivatives along the line
${\bold y}(r) = (1+r,1+2r)$. Note that for most of these ${\bold y}$-values,
the two pulses overlap at the final time, $T=1$.

\section{Numerical Examples}
\label{sec5}

In this section, we present two numerical examples to demonstrate the efficiency and applicability of the method proposed in Section \ref{sec3}.

We consider the Cauchy problem \eqref{wave} in a two-dimensional spatial space and let ${\bf x}=(x_1,x_2) \in {\mathbb R}^2$. We employ the proposed stochastic spectral asymptotic method to approximate the solution, $u^{\varepsilon}$, and the expected value of the QoI in \eqref{Q}. The QoI test functions are given in terms of the smooth function, $\widetilde{\psi} \in C_{\text{c}}^{\infty}({\mathbb R}^2)$, in \eqref{testfun} as $\psi(x_1,x_2) = \widetilde{\psi}(2 x_1, 2 x_2)$ and $\psi(x_1,x_2) = \widetilde{\psi}(x_1-1,x_2)$ for the first and second numerical examples, respectively. In both examples, we use the Smolyak sparse grid based on Clenshaw-Curtis abscissas and the nested rule \eqref{pj2}. We show that fast convergence rates are obtained as predicted in Section \ref{sec4}.
For each example, we also compare the
convergence rate with the average rate obtained from
ten independent Monte Carlo simulations.

As in Section \ref{sec:4_2}, the step size used in the quadrature approximation of \eqref{Q} is chosen as $\Delta x = \frac{2 \pi \varepsilon}{10}$ and the space between the Gaussian beams in \eqref{trapets} as $\Delta z = \sqrt{\varepsilon}$.

\subsection{Numerical test 1: Two pulses}

In this example, both wave speed and initial data are uncertain and described by a random vector, ${\bf y}=(y_1, \dotsc, y_5)$, containing $N = 5$ independent uniformly distributed random variables. The constant random wave speed is given by
$$
c({\bold y}) = y_1 \sim {\mathcal U}(0.8,1),
$$
and the random initial data are given by \eqref{waveinit} with
$$
\Phi_0({\bf x}) = |x_1|, \qquad  A_0({\bf x},{\bf y}) = g\Bigl(\mathbf x-\mathbf s_1,\mathbf d_1({\bf y})\Bigr) + g\Bigl(\mathbf x-\mathbf s_2({\bf y}),\mathbf d_2({\bf y})\Bigr),
$$
where
$$
g(\mathbf x,\mathbf d) = e^{-(d_{1} \, x_1^2 + d_{2} \, x_2^2)},
$$
and
$$
\mathbf s_1 = (-1,0), \qquad \mathbf s_2({\bf y}) = (y_2,y_3), \qquad \mathbf d_1({\bf y}) = (y_4,5), \qquad \mathbf d_2({\bf y}) = (y_5,y_5),
$$
and
$$
y_2 \sim {\mathcal U}(1,1.5), \qquad y_3 \sim {\mathcal U}(0,0.5), \qquad y_4 \sim {\mathcal U}(5,10), \qquad  y_5 \sim {\mathcal U}(5,10).
$$
Hence, the initial solution consists of two Gaussian wave pulses, and the vectors $\mathbf s_j$ and $\mathbf d_j$, with $j=1,2$, represent the position and shape of the pulses, respectively.

\begin{figure}
\vskip -.2cm
\centering
\begin{minipage}{.32\textwidth}
      \subfigure[${\bf y}=(0.8,1,0.5,5,5)$]{\includegraphics[width=\textwidth]{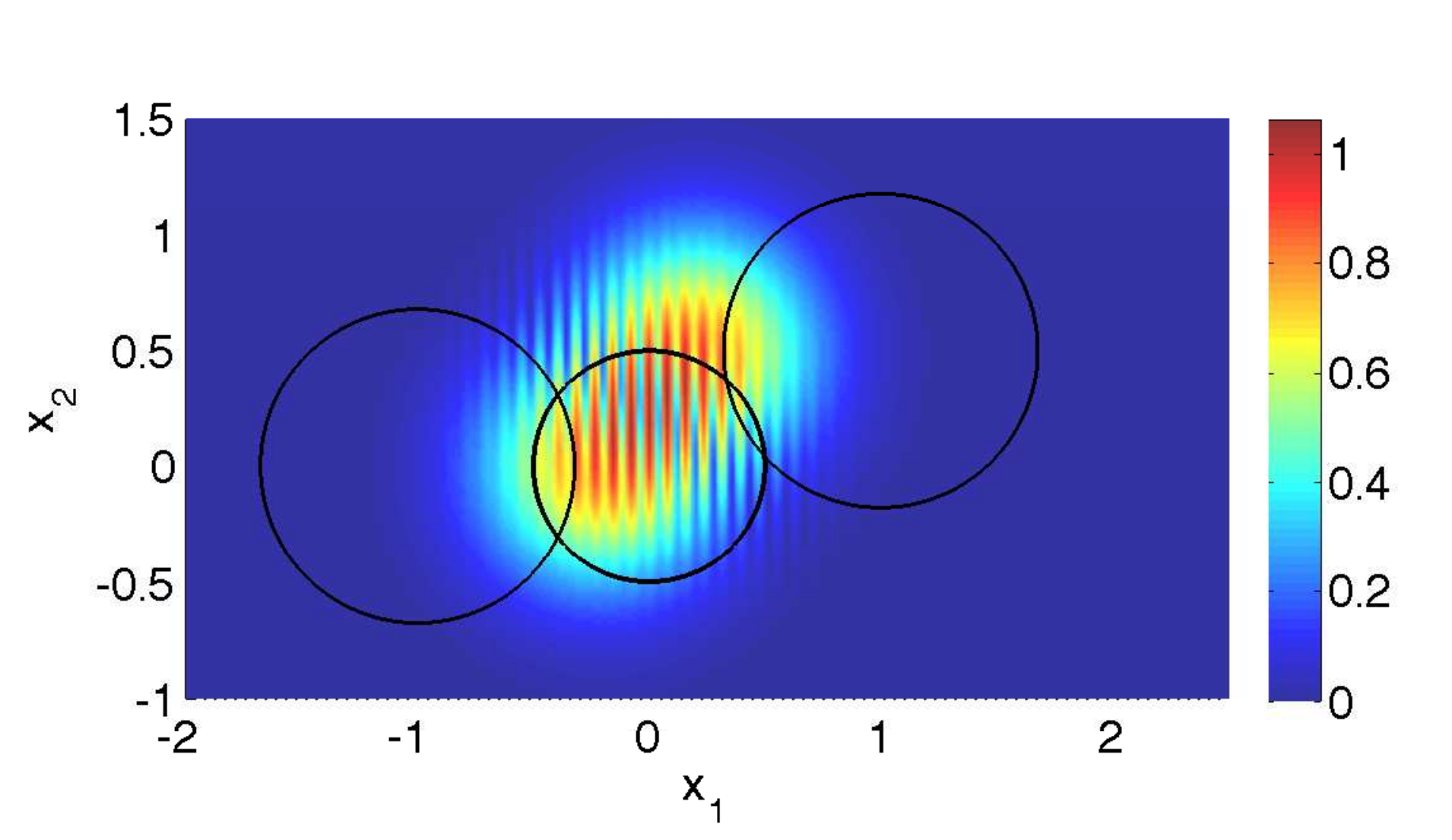}}
\end{minipage}
\begin{minipage}{.32\textwidth}
      \subfigure[${\bf y}=(0.8,1,0.5,10,10)$]{\includegraphics[width=\textwidth]{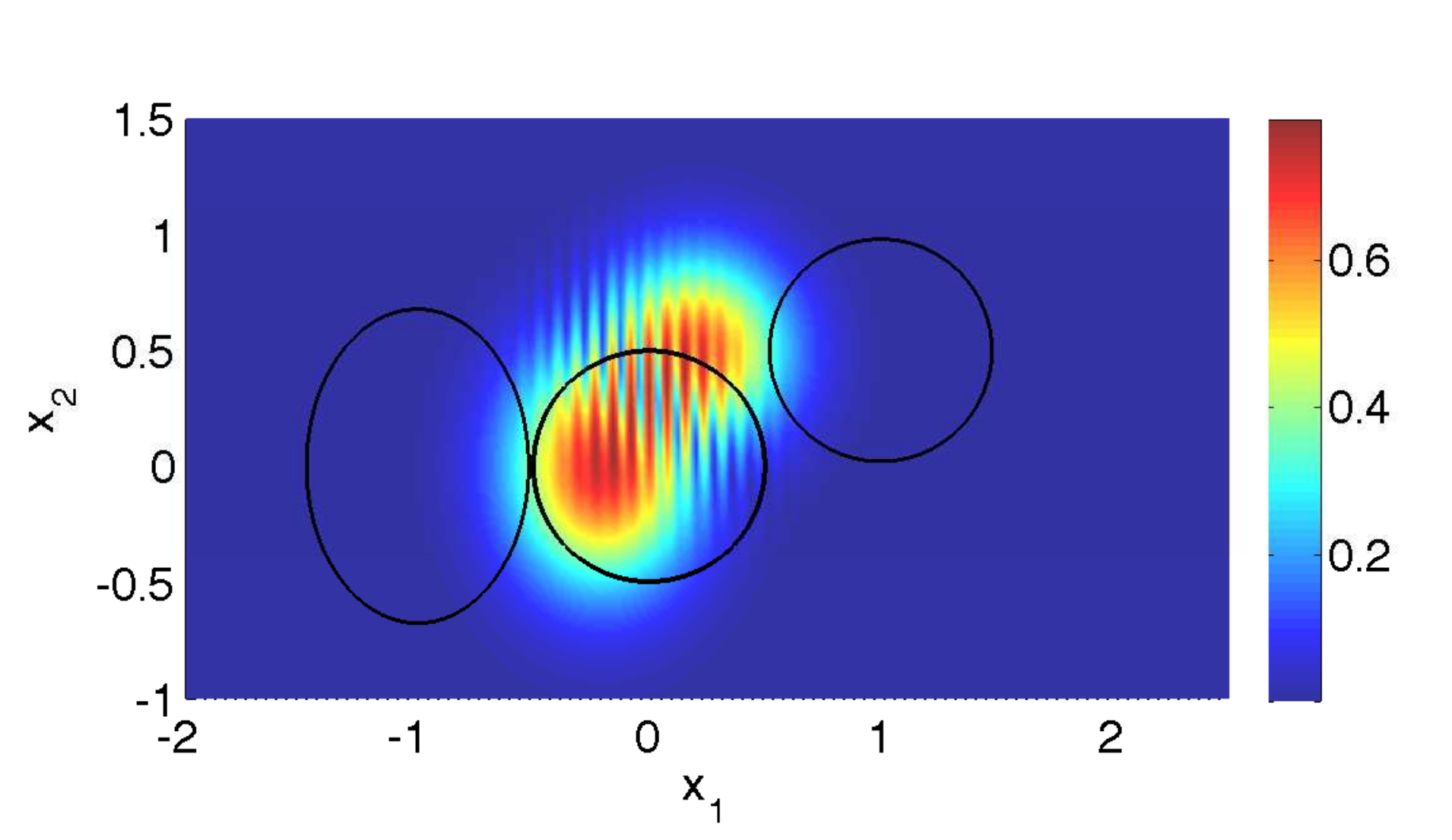}}
\end{minipage}
\begin{minipage}{.32\textwidth}
      \subfigure[${\bf y}=(1,1.5,0,5,5)$]{\includegraphics[width=\textwidth]{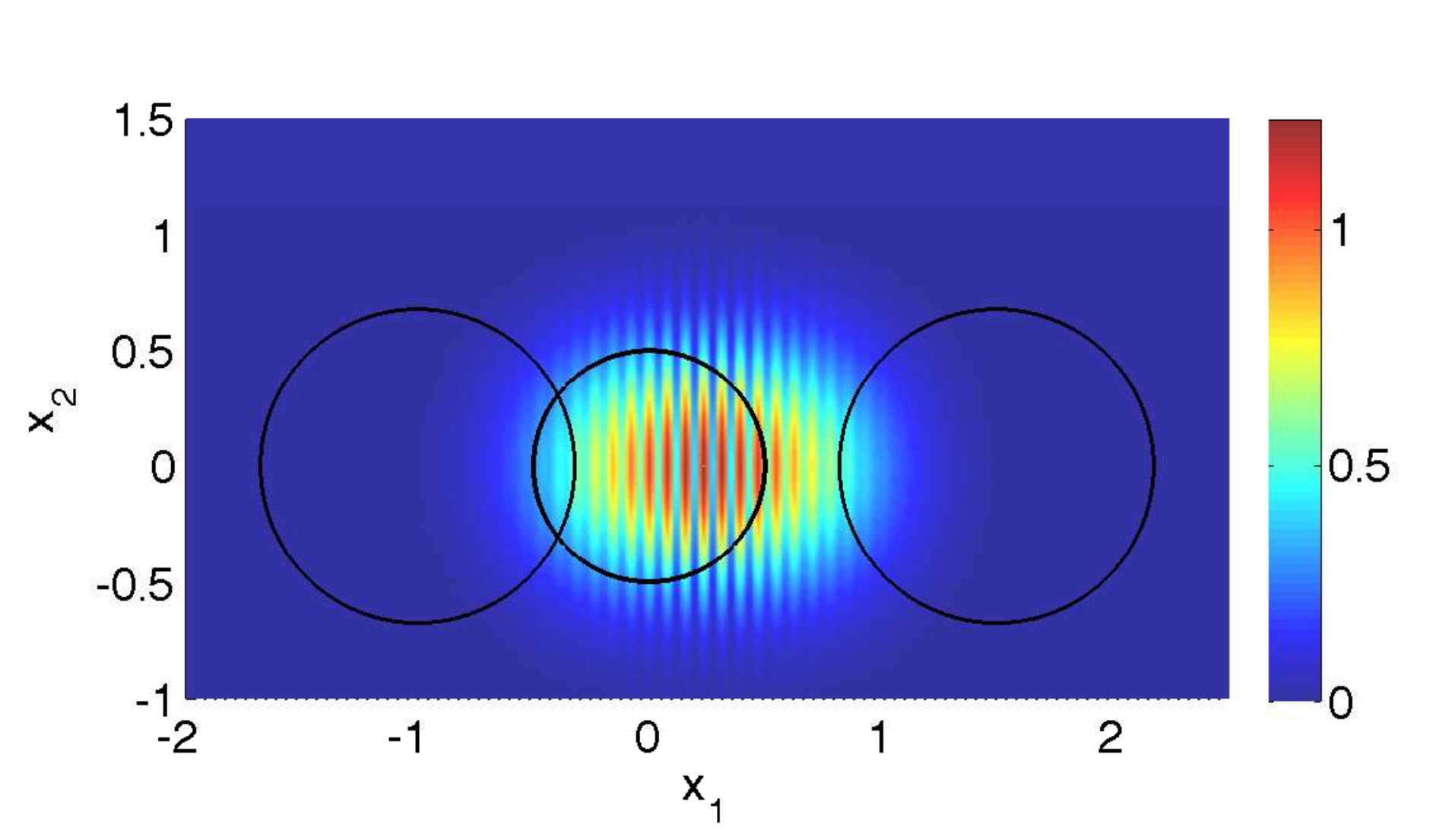}}\\
\end{minipage}
\begin{minipage}{.32\textwidth}
      \subfigure[${\bf y}=(0.8,1.5,0,10,10)$]{\includegraphics[width=\textwidth]{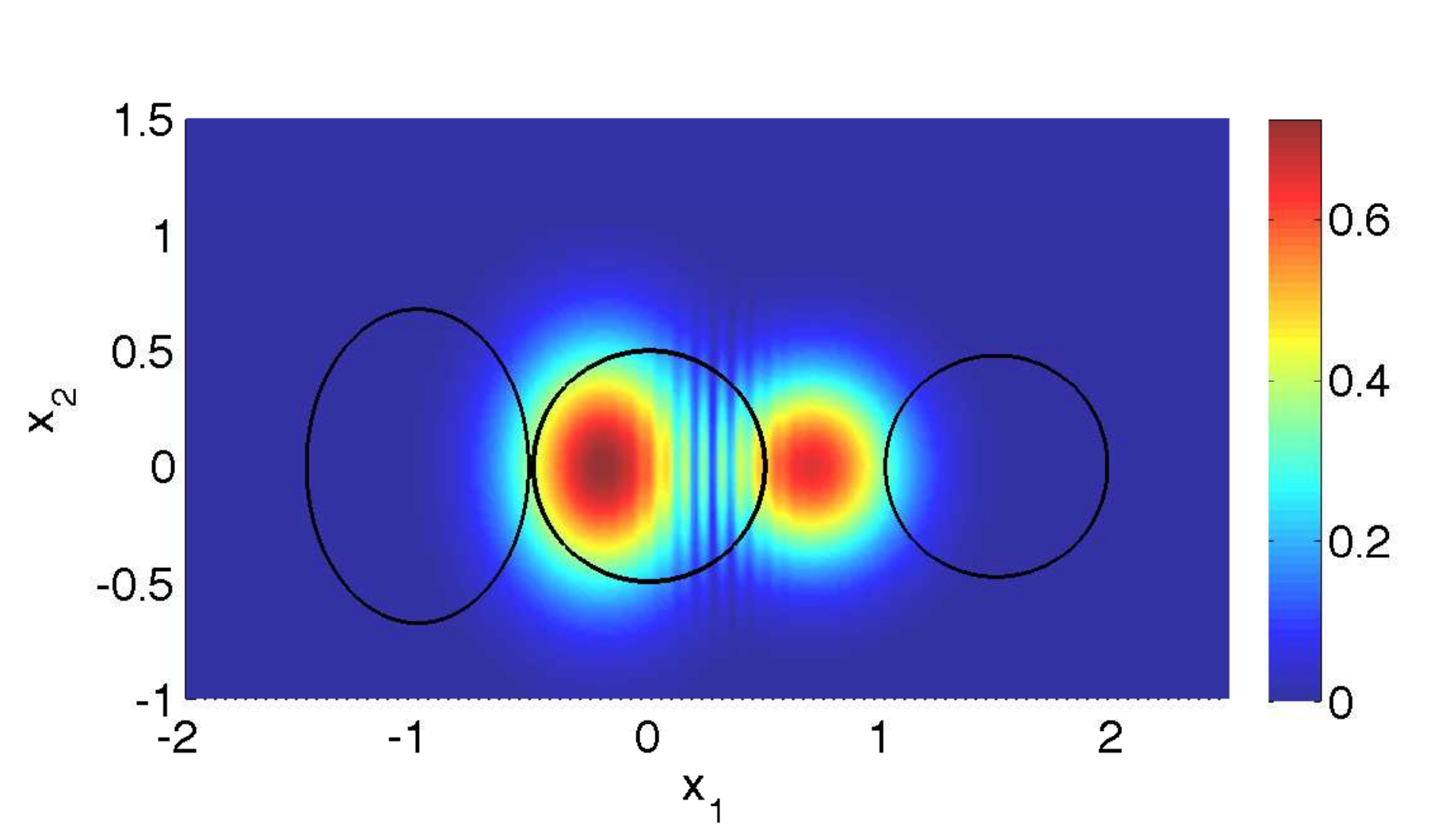}}
\end{minipage}
\begin{minipage}{.32\textwidth}
      \subfigure[${\bf y}=(1,1,0.5,5,5)$]{\includegraphics[width=\textwidth]{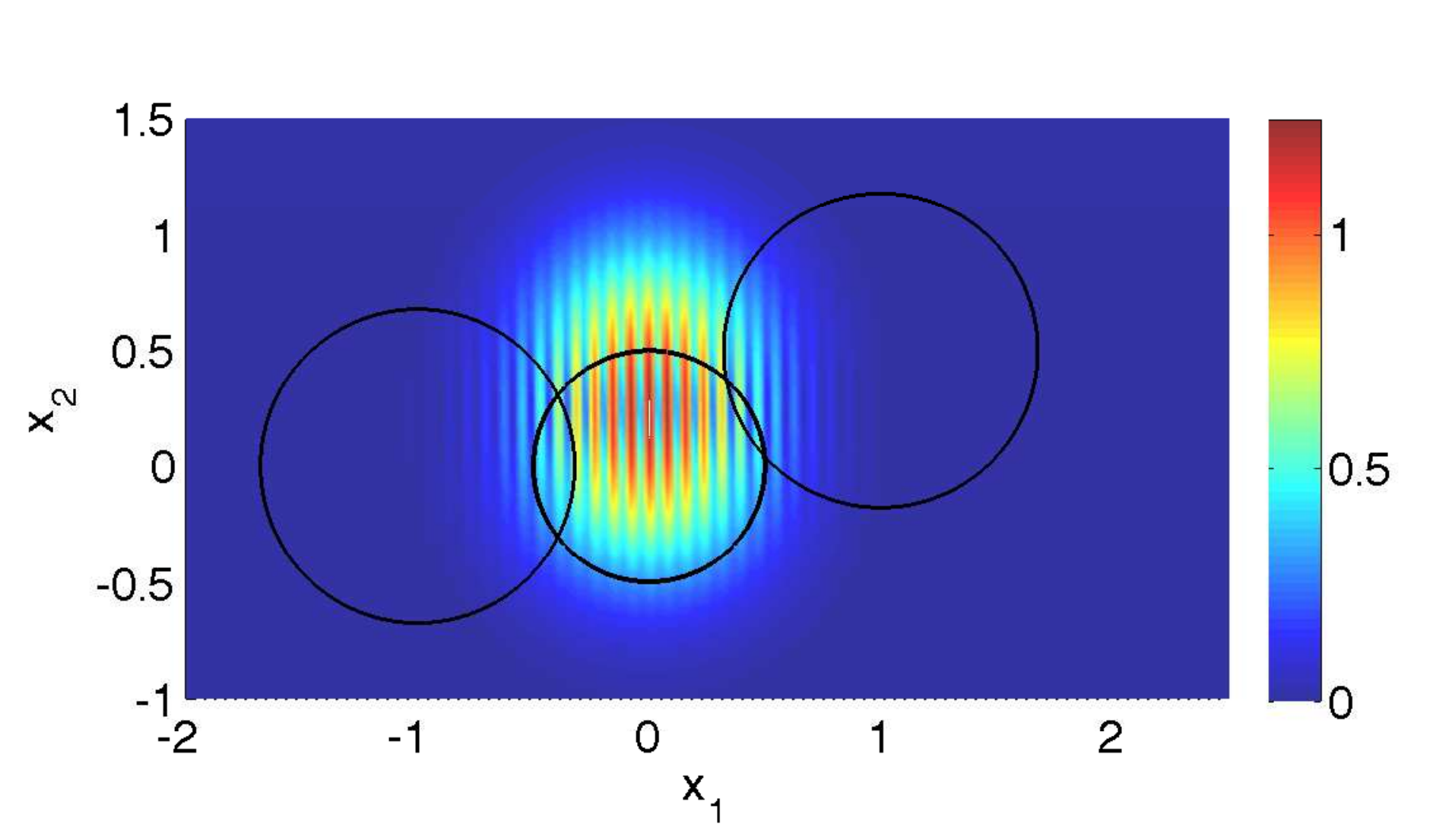}}
\end{minipage}
\begin{minipage}{.32\textwidth}
      \subfigure[${\bf y}=(0.8,1,0,5,10)$]{\includegraphics[width=\textwidth]{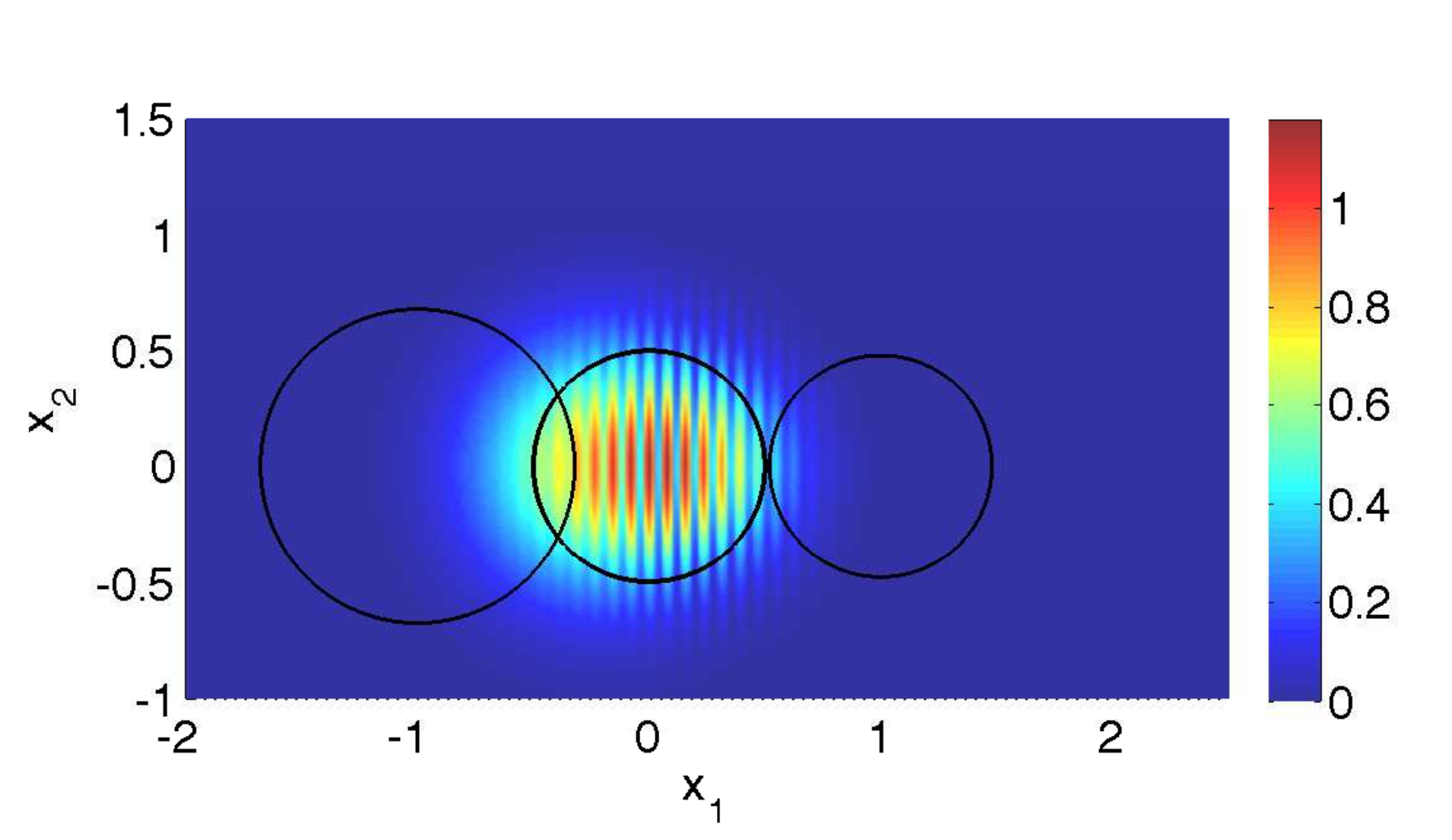}}
\end{minipage}
\caption{Numerical test 1: Six realizations of the magnitude of approximate solution $|u_{\text{GB}}^{\varepsilon}(T,{\bf x},{\bf y})|$ at the fixed time $T = 1$ with wavelength $\varepsilon = 1/40$. In each realization, the central circle indicates the support of the QoI test function and the other two circles/ellipses indicate the supports of the initial solution.}
\label{fig:Test1_sol}
\end{figure}

Figure \ref{fig:Test1_sol} shows six realizations of the magnitude of the approximate solution $|u_{\text{GB}}^{\varepsilon}(T,{\bf x},{\bf y})|$
with wavelength $\varepsilon = 1/40$
at the time $T = 1$. In each realization, the central circle indicates the support of the QoI test function and the other two circles/ellipses indicate the supports of the initial solution that consists of two Gaussian pulses.

We make a convergence study for
a set of wavelengths, $\varepsilon = 1/40, 1/80, 1/160$. For each wavelength, we consider different levels, $\ell \ge 1$, and compute the relative error in the expected value of the QoI in \eqref{Q}:
\begin{equation}\label{abs_error}
{\mathcal E}(\eta(\ell)) :=
\frac{
\Bigl| {\mathbb E}[ {\mathcal S}_{{\mathcal I}(\ell_{\text{ref}})} [{\mathcal Q}^{\varepsilon}]]   -   {\mathbb E}[ {\mathcal S}_{{\mathcal I}(\ell)} [{\mathcal Q}^{\varepsilon}]]   \Bigr|}
{\Bigl|{\mathbb E}[ {\mathcal S}_{{\mathcal I}(\ell_{\text{ref}})} [{\mathcal Q}^{\varepsilon}]]   \Bigr|}.
\end{equation}
Here, for each wavelength, the reference solution is computed with a high level, $\ell_{\text{ref}}$, and with the same Gaussian beam parameters used in all levels, $\ell \ge 1$. The error \eqref{abs_error} therefore reflects only the stochastic collocation error, not the error in the deterministic asymptotic solver.

\begin{figure}
\vskip -.1cm
\centering
                \begin{minipage}{0.8\textwidth}
                    \includegraphics[width=\textwidth]{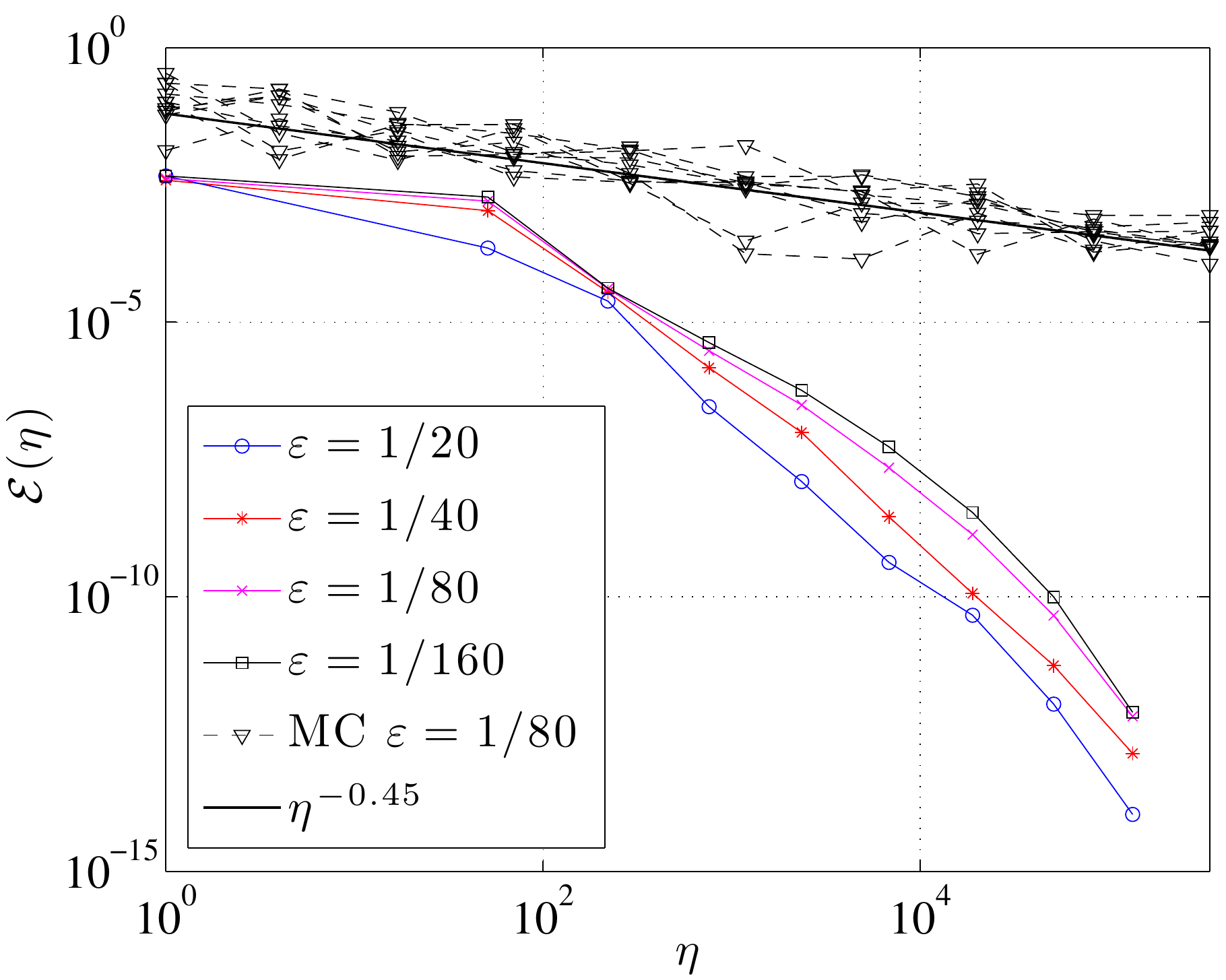}
                \end{minipage}
\vskip -.2cm
\caption{Numerical test 1: Relative error ${\mathcal E}(\eta)$ at time $T=1$ versus the number of collocation points $\eta$ (or the number of samples in the case of Monte Carlo sampling), for various wavelengths. The proposed method performs a fast spectral convergence, while Monte Carlo sampling has a slow algebraic convergence. The rate of convergence of Monte Carlo sampling, obtained by linear regression through the data points, is 0.45.}
\label{fig:Test1_conv}
\end{figure}

Figure \ref{fig:Test1_conv} shows the relative error, ${\mathcal E}(\eta)$, in \eqref{abs_error} at time $T=1$, computed by the proposed method, versus the number of collocation points, $\eta$, for various wavelengths.
It also shows the convergence of the relative error
in ten Monte Carlo runs,
computed
using the same reference solution as above,
 with $\eta$ representing the number
of samples.

We observe a fast spectral convergence rate of the stochastic collocation error in the proposed method, due to the high stochastic regularity of the QoI. A simple linear regression through the data points shows that the rate of convergence of Monte Carlo is 0.45, which is very slow. Consequently, the decay in the stochastic collocation error is much faster than the decay in Monte Carlo error. We also note that as $\varepsilon$ decreases, the
decay rate does not deteriorate.
This points to
the existence of uniform bounds \eqref{est}.

\subsection{Numerical test 2: A Lens}

In this example, the uncertain wave speed is described by a random vector ${\bf y}=(y_1, y_2, y_3)$ containing $N = 3$ independent uniformly distributed random variables, given by
$$
c({\bf x},{\bf y})= 1 - y_1 \, e^{- (y_2 \, x_1^2 - y_3 \,  x_2^2)},
$$
where
$$
y_1 \sim {\mathcal U}(0,0.4), \qquad y_2 \sim {\mathcal U}(0.65,0.85), \qquad y_3 \sim {\mathcal U}(0, 1).
$$
The initial data are assumed to be deterministic and given by \eqref{waveinit} with
$$
\Phi_0({\bf x}) = -x_1, \qquad A_0({\bf x})= e^{-5 \, (x_1+1)^2}.
$$
The problem models a plane wave that is refracted by a lens of uncertain
shape and strength.

\begin{figure}
\centering
    \begin{minipage}{0.49\textwidth}
       \subfigure[${\bf y}=(0.2, 0.75, 0.5)$]
       {\includegraphics[width=\textwidth]{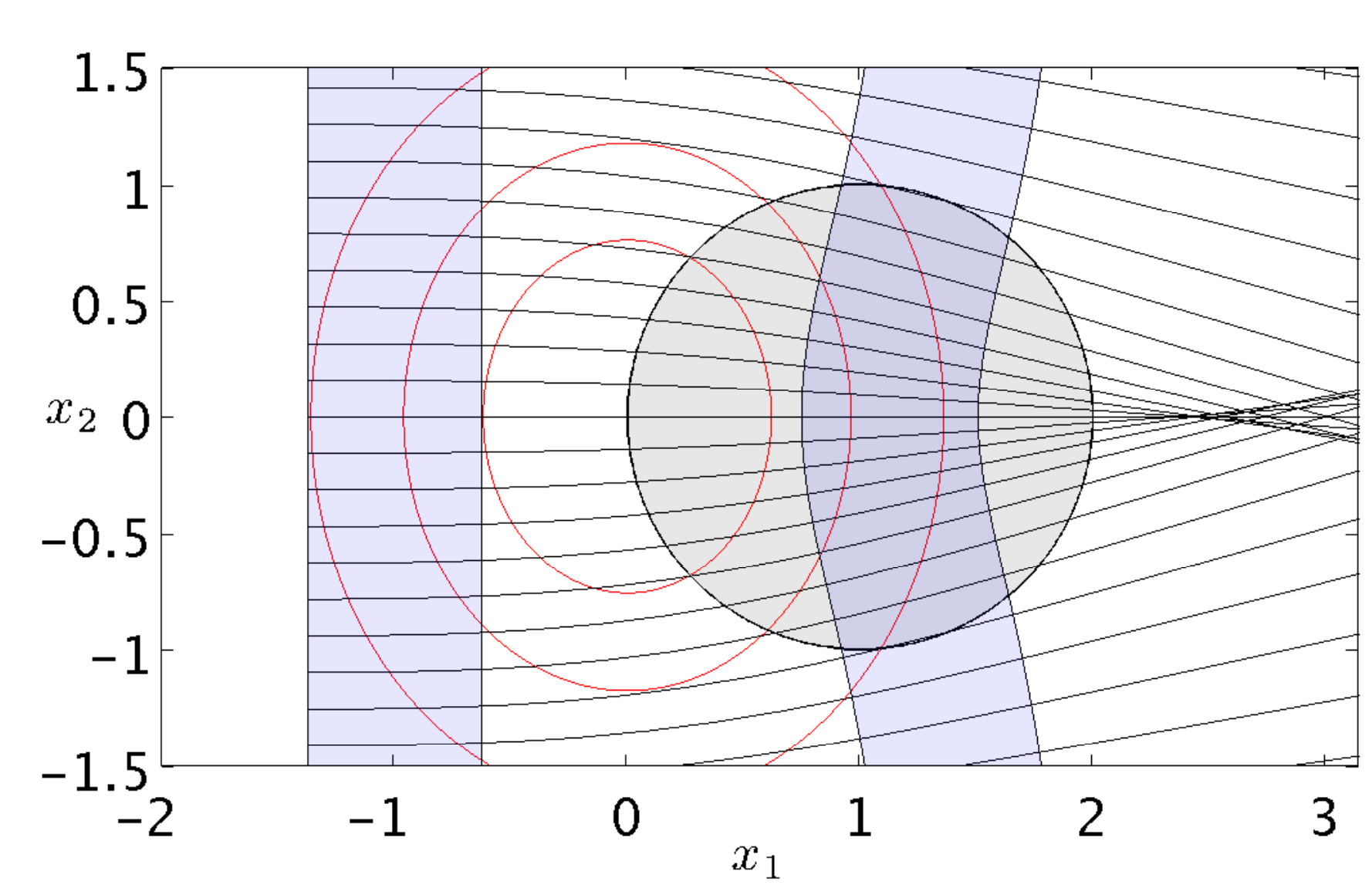}}
    \end{minipage}
    \begin{minipage}{0.49\textwidth}
            \subfigure[${\bf y}=(0.4, 0.75, 0.5)$]
       {\includegraphics[width=\textwidth]{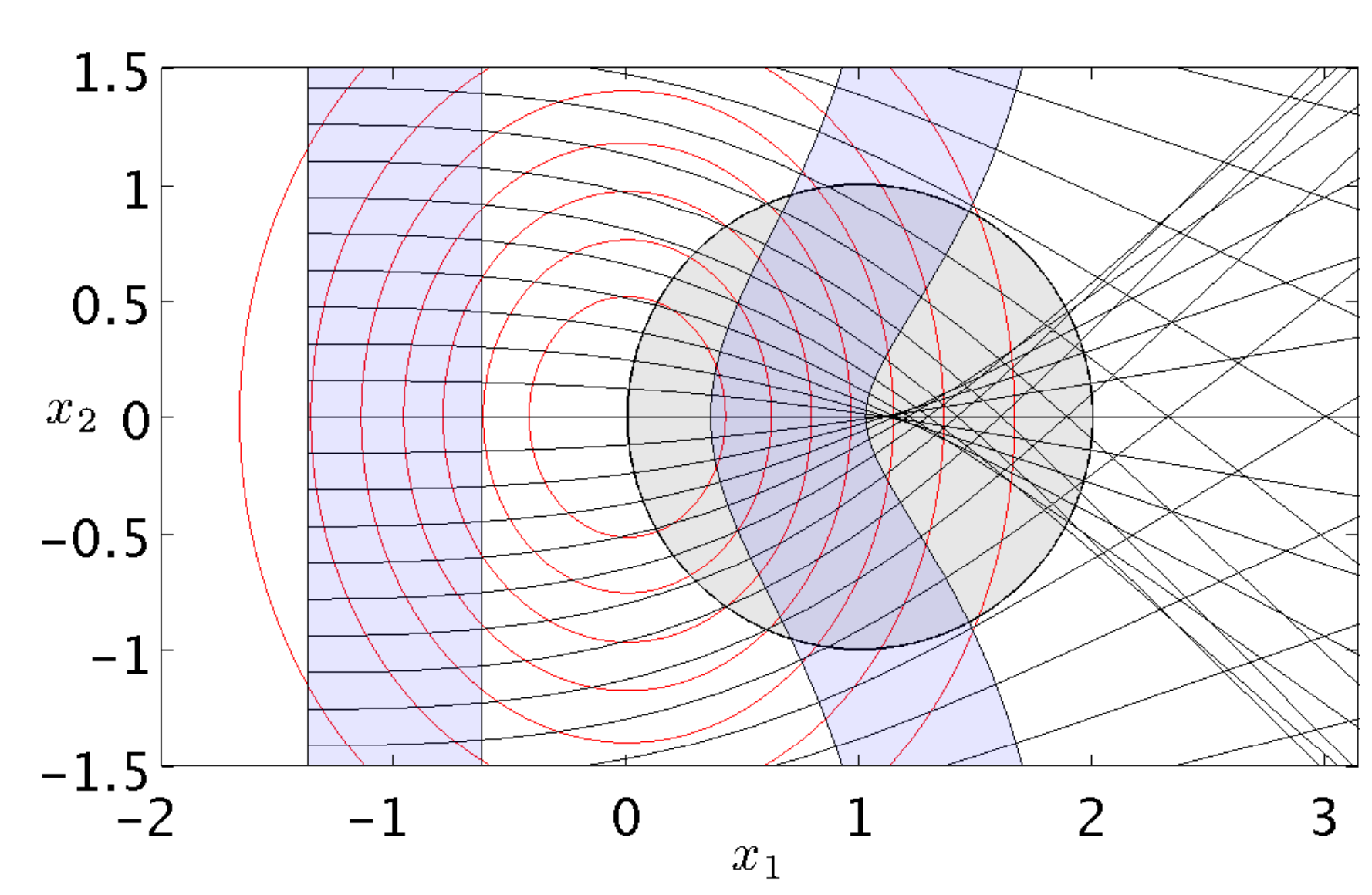}}\\
    \end{minipage}
    \begin{minipage}{0.49\textwidth}
       \subfigure[${\bf y}=(0.2, 0.75, 1.0)$]
       {\includegraphics[width=\textwidth]{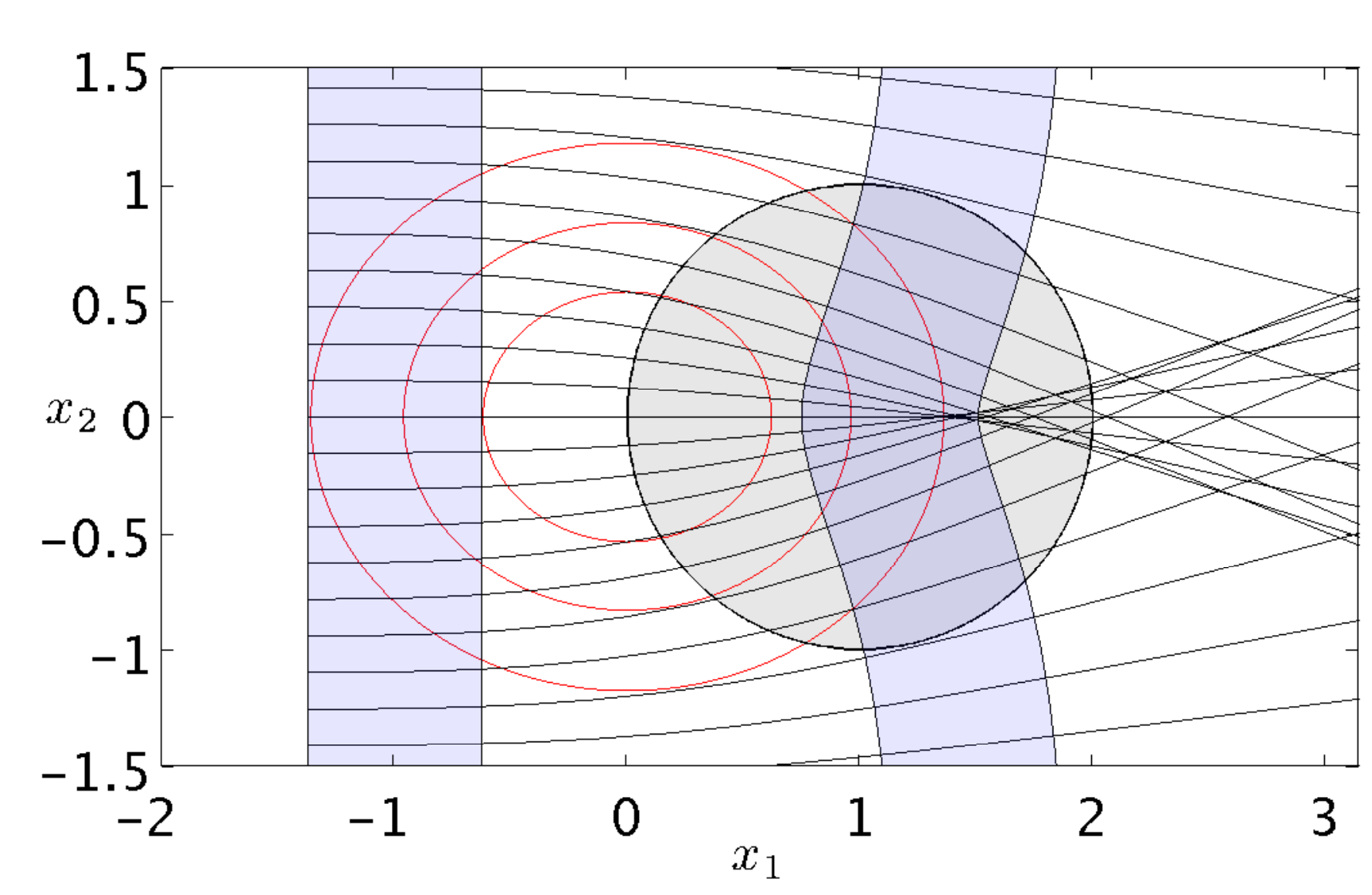}}
    \end{minipage}
    \begin{minipage}{0.49\textwidth}
       \subfigure[${\bf y}=(0.4, 0.75, 1.0)$]
       {\includegraphics[width=\textwidth]{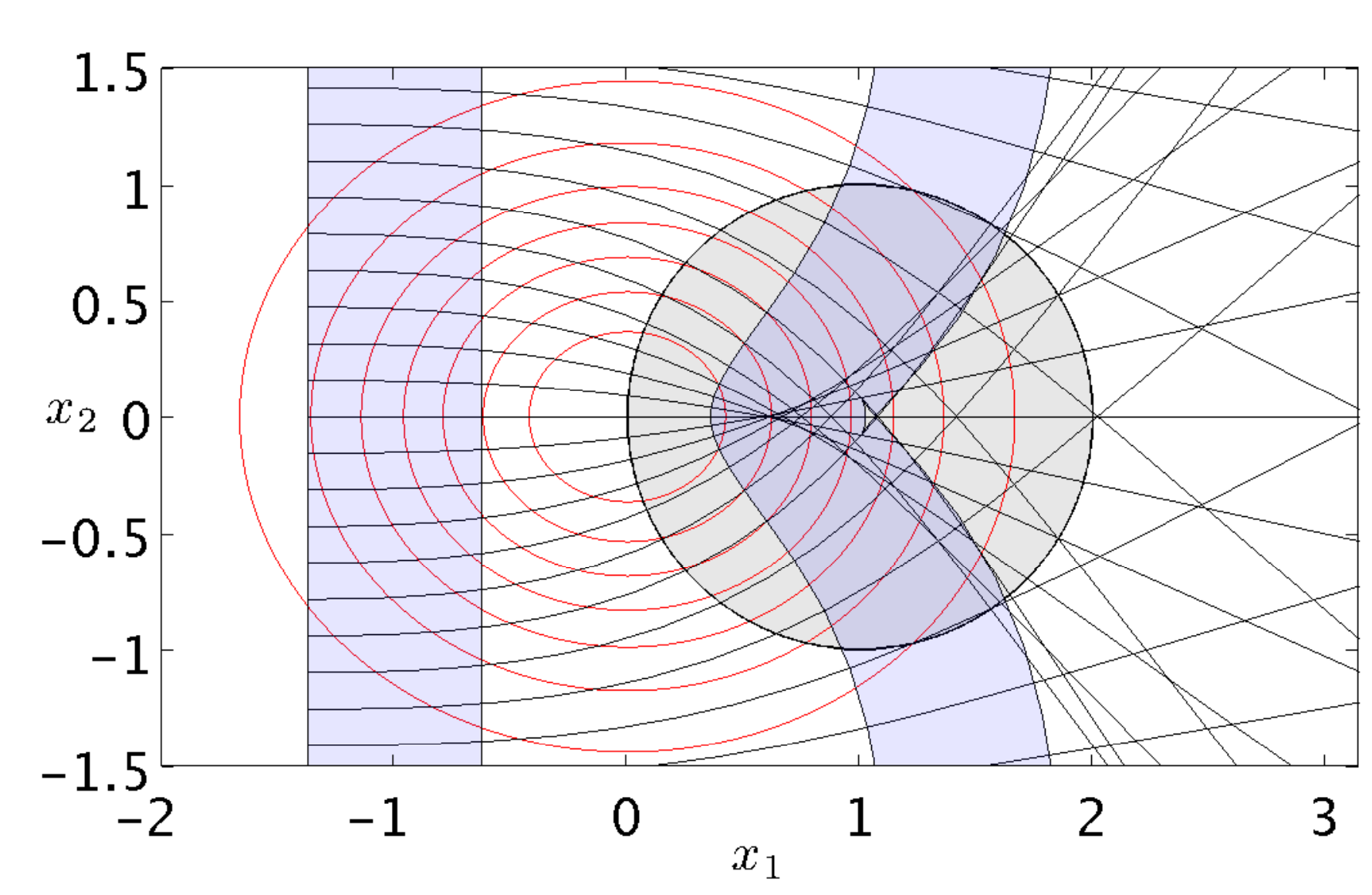}}
    \end{minipage}
\caption{Numerical test 2: Ray tracing solution for four realizations of ${\bold y}$. A cusp caustic and two fold caustics are formed in all realizations. Contour lines of $c({\bold x},{\bold y})$ overlaid in red. The left transparent band shows where the initial amplitude is above $1/2$; the right band shows how this band has been transported by the rays at time $t=2.5$. Circle indicates the support of QoI test function.}
\label{fig:raypic}
\end{figure}

The wave speed varies in the spatial domain and caustics may consequently form. Figure \ref{fig:raypic} shows the ray tracing solution for four different realizations of the random vector ${\bold y}$. A cusp caustic and two fold caustics are formed inside the support of the QoI test function
for the last three realizations, but not in the first one.

\begin{figure}
            \centering
    \begin{minipage}{0.49\textwidth}
            \subfigure[$t=0$]
        {\includegraphics[width=\textwidth]{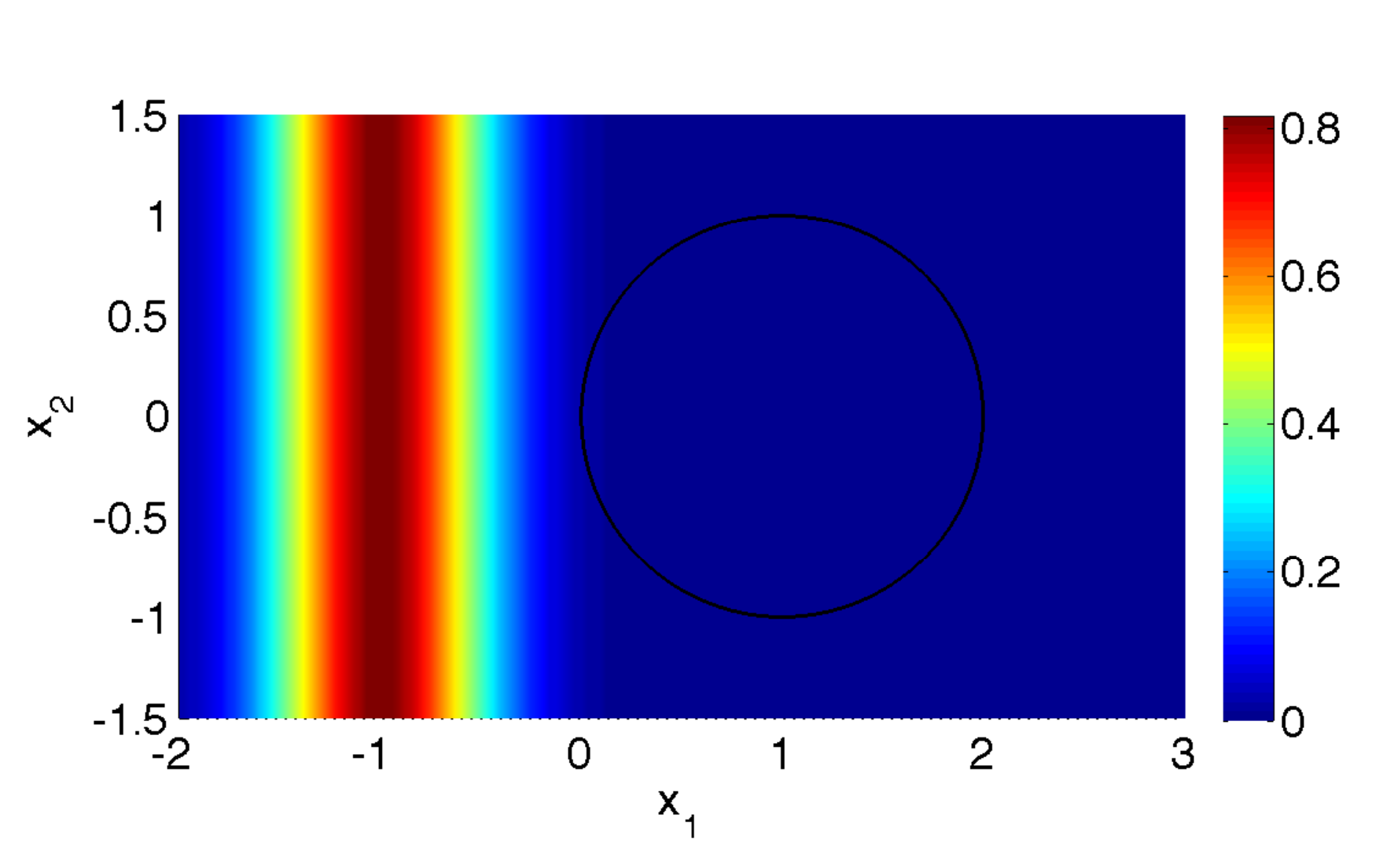}}
    \end{minipage}
    \begin{minipage}{0.49\textwidth}
            \subfigure[$t=1$]
        {\includegraphics[width=\textwidth]{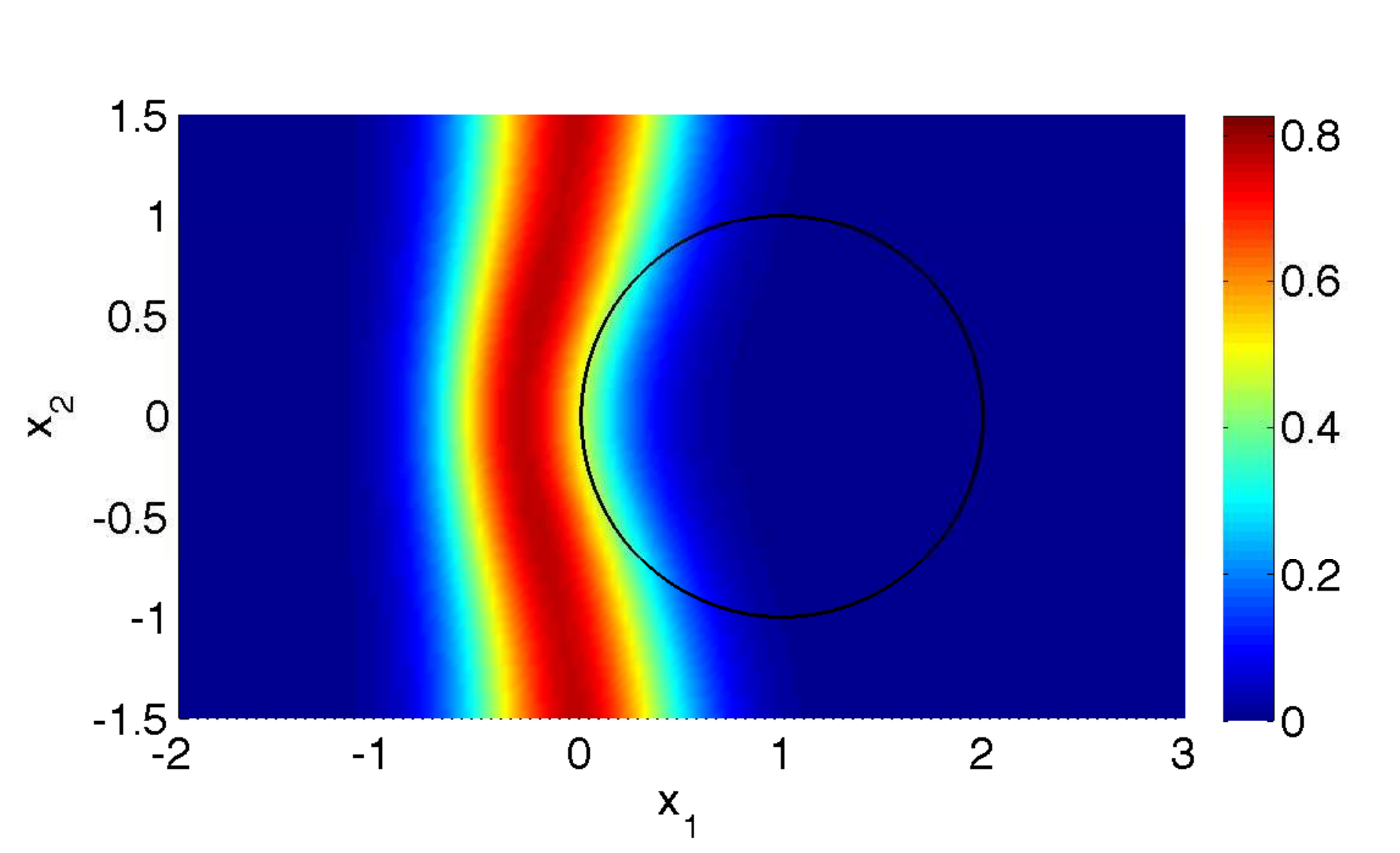}}\\
    \end{minipage}
    \begin{minipage}{0.49\textwidth}
            \subfigure[$t=2$]
        {\includegraphics[width=\textwidth]{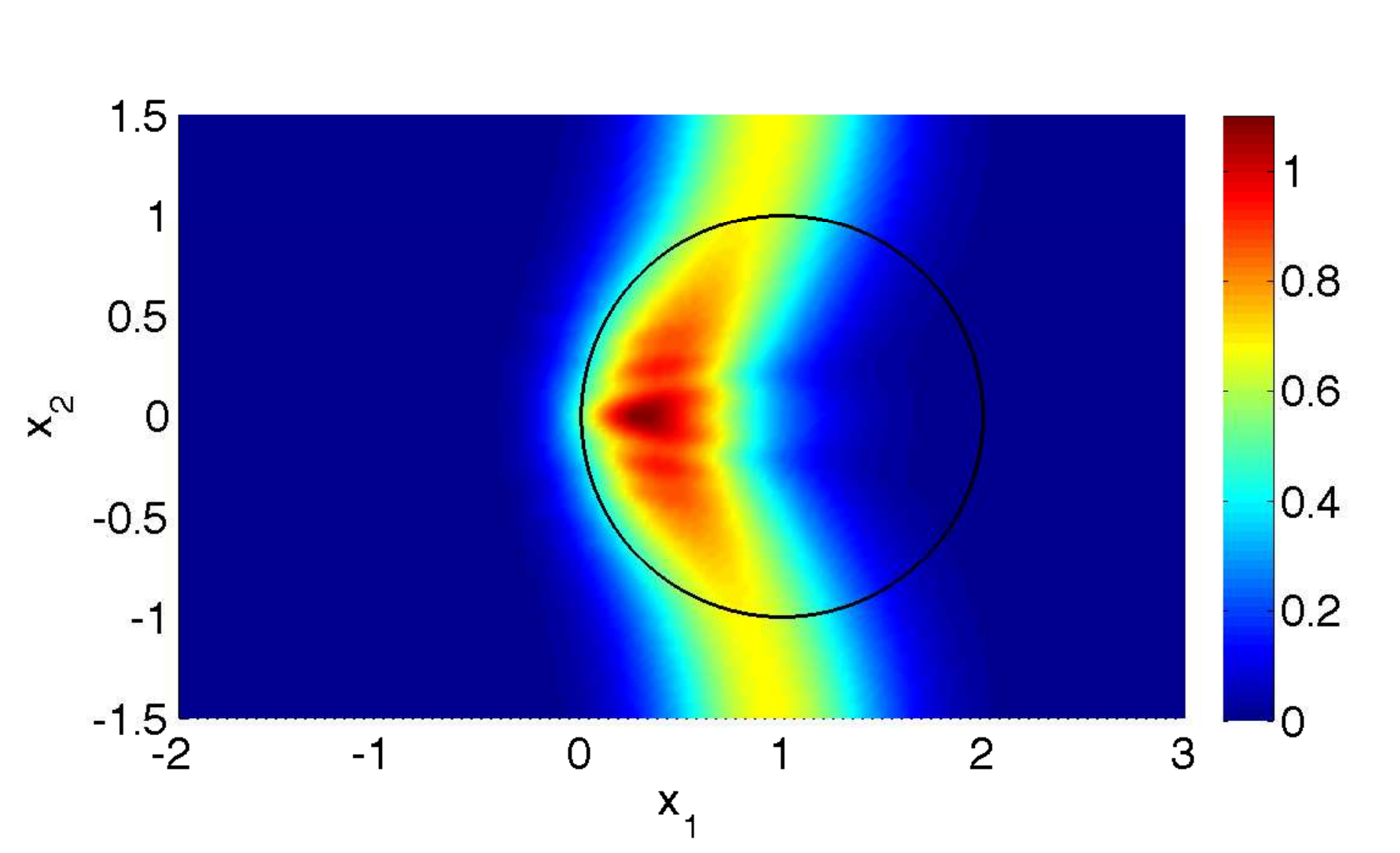}}
    \end{minipage}
    \begin{minipage}{0.49\textwidth}
            \subfigure[$t=2.5$]
        {\includegraphics[width=\textwidth]{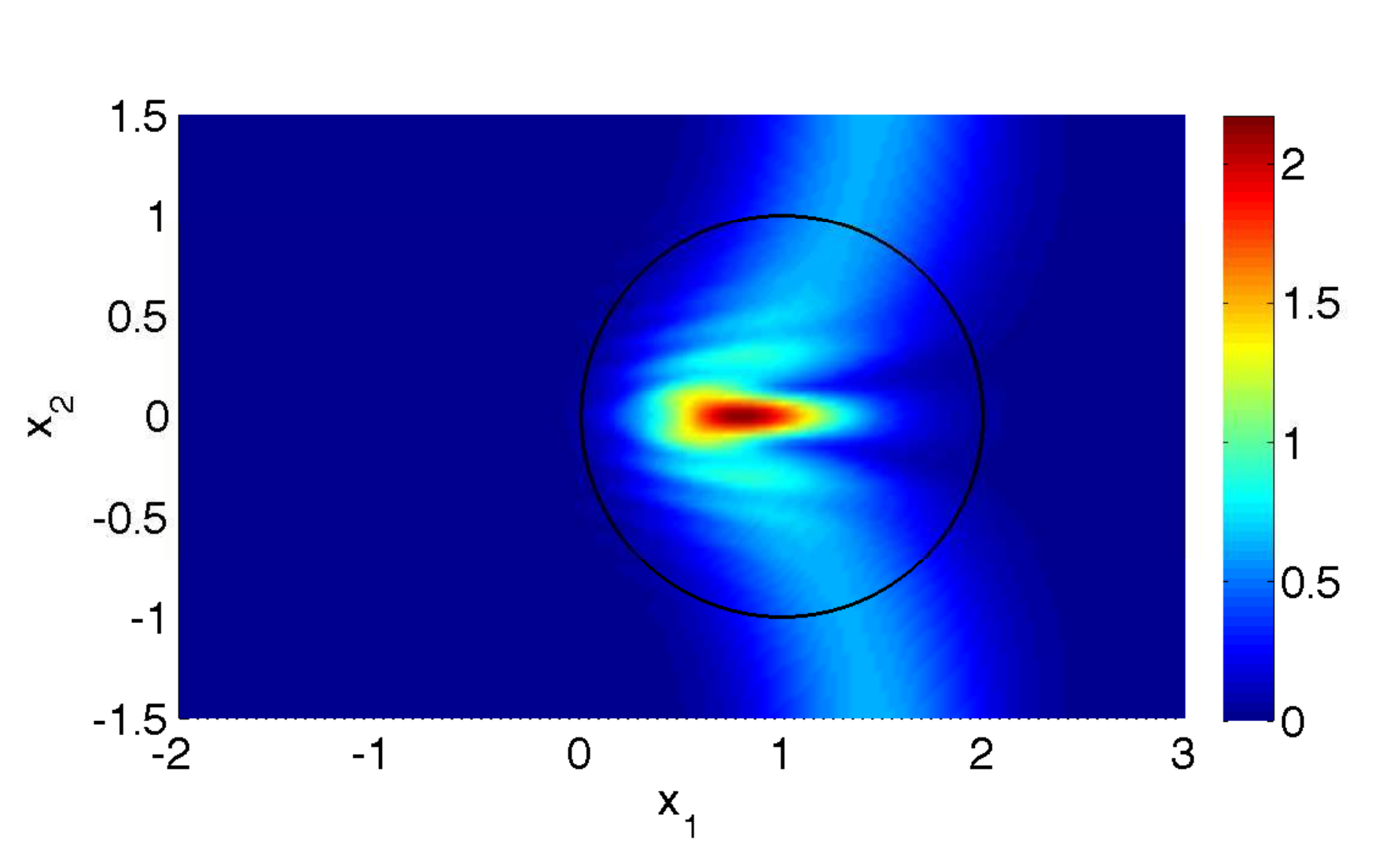}}
    \end{minipage}
            \caption{Numerical test 2: Magnitude of approximate solution $|u_{\text{GB}}^{\varepsilon}(t,{\bf x},{\bf y})|$ with wavelength $\varepsilon = 1/20$ at various time instances and a fixed random vector ${\bf y}=(0.4, 0.75, 1)$, which corresponds to the realization in Figure \ref{fig:raypic}d. Circle indicates the support of QoI test function.}
\label{fig:Test2_sol}
\end{figure}

Figure \ref{fig:Test2_sol} shows the magnitude of the approximate solution $|u_{\text{GB}}^{\varepsilon}(t,{\bf x},{\bf y})|$ with wavelength $\varepsilon = 1/20$ at various time instances and a fixed random vector, ${\bf y}=(0.4, 0.75, 1)$, which corresponds to the realization in Figure \ref{fig:raypic}d. The central circles indicate the support of the QoI test function.

\begin{figure}
\vskip -.1cm
\centering
\begin{minipage}{0.8\textwidth}
    \includegraphics[width=\textwidth]{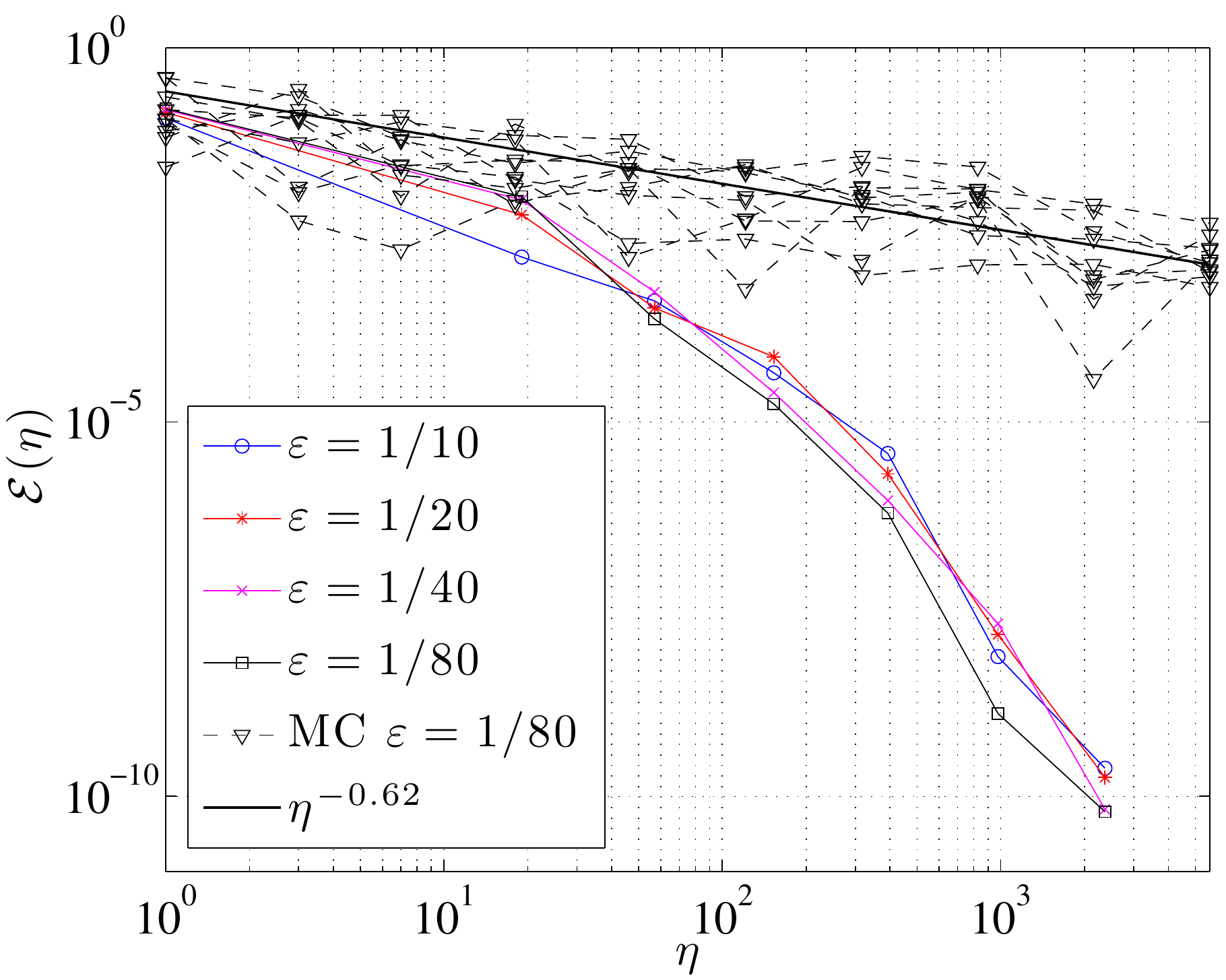}
\end{minipage}
\vskip -.2cm
\caption{Numerical test 2: Relative error ${\mathcal E}(\eta)$ at time $T=2.5$ versus the number of collocation points, $\eta$ (or the number of samples in the case of Monte Carlo sampling), for various wavelengths. The proposed method performs a fast spectral convergence, while Monte Carlo sampling has a slow algebraic convergence. The rate of convergence of Monte Carlo sampling, obtained by linear regression through the data points, is 0.62.}
\label{fig:Test2_conv}
\end{figure}

Figure \ref{fig:Test2_conv} shows the relative error, ${\mathcal E}(\eta)$, in \eqref{abs_error} at time $T=2.5$, computed by the proposed method, versus the number of collocation points, $\eta$ (or the number of samples in the case of Monte Carlo sampling), for various wavelengths.

Similar to the first numerical test, we observe a fast spectral convergence rate of the stochastic collocation error in the proposed method.
The convergence rate of Monte Carlo, given by linear regression, is in
this case 0.62, and the error again decays much more slowly.
Furthermore, the error curves of the proposed method have a rather
unform shape for all the $\varepsilon$ used, suggesting that \eqref{est} holds.

\section{Conclusion}
\label{sec6}

We have proposed a novel stochastic spectral asymptotic method for the forward propagation of uncertainty in high-frequency waves generated by highly oscillatory initial data. The source of uncertainty is the wave speed and/or the initial data, characterized by a finite number of independent random variables with known probability distributions. The proposed method combines a sparse stochastic collocation method for propagating uncertainty and a Gaussian beam summation method for propagating high-frequency waves. Fast error convergence is attained only in the presence of QoIs which are smooth with respect to input random parameters {\it independent of the wave frequency}.

The wave solution is highly oscillatory in both physical and stochastic spaces,
and its derivatives clearly cannot be bounded independently of the frequency.
A priori, QoIs based on the solution would have the same behavior.
However, our main result is that
there are in fact quadratic QoIs that are smooth with uniformly bounded derivatives in the stochastic space. Through both theoretical arguments for simplified problems and numerical experiments for more complicated problems, we have verified the spectral convergence of the proposed method for such a quadratic QoI, which represents the local wave strength. This shows that the proposed method may be a valid alternative to the traditional Monte Carlo method.

Future directions include a rigorous proof of the uniform bounds in Conjecture 1 for the quadratic quantity considered in this work and a regularity analysis of other types of nonlinear QoIs.

\vskip 1cm
\bibliographystyle{plain}
\bibliography{refs_me}

\end{document}